\documentclass[opre,nonblindrev]{informs4}

\OneAndAHalfSpacedXI
\makeatletter
\renewcommand{\AAnormalsizeXI}{%
  \@setfontsize\normalsizeXI{11}{11.5}%
  \Dskips{1.0}{0.3}{0.5}{0.6}%
  \let\@listi\@listI
  \baselineskip17.25pt plus 1pt minus 1pt
  \@bls=17.25pt\@em=10pt\relax}
\normalsize
\makeatother
\usepackage{anyfontsize}
\usepackage[T1]{fontenc}
\RequirePackage{endnotes}
\usepackage{natbib}
\usepackage{amsmath,amssymb,amsfonts}
\usepackage{graphicx}
\usepackage{booktabs}
\usepackage{array}
\usepackage{multirow}
\usepackage{tabularx}
\usepackage{longtable}
\usepackage{float}
\usepackage{placeins}
\usepackage[ruled,vlined]{algorithm2e}
\usepackage{enumitem}
\usepackage[table]{xcolor}
\usepackage[
  colorlinks=true,
  linkcolor=blue!50!black,
  citecolor=blue!50!black,
  urlcolor=blue!50!black
]{hyperref}

\makeatletter
\newcommand{\paperRunningHead}[2]{%
  \vbox{%
    \hbox to\textwidth{\normalfont\normalcolor\fontsize{8}{10}\selectfont
      \strut #1\hfil #2}%
    \kern4pt\hrule height0.4pt}}
\if@BLINDREV
\renewcommand{\setoddRH}{%
  \paperRunningHead{Authors' names not included for peer review\quad\textit{\theRUNTITLE}}{\thepage}}
\renewcommand{\setevenRH}{%
  \paperRunningHead{\thepage}{Authors' names not included for peer review\quad\textit{\theRUNTITLE}}}
\renewcommand{\setoddECRH}{%
  \paperRunningHead{Authors' names not included for peer review\quad\textit{\theECRUNTITLE}\quad E-Companion}{\thepage}}
\renewcommand{\setevenECRH}{%
  \paperRunningHead{\thepage}{Authors' names not included for peer review\quad\textit{\theECRUNTITLE}\quad E-Companion}}
\else
\renewcommand{\setoddRH}{%
  \paperRunningHead{\theRUNAUTHOR\quad\textit{\theRUNTITLE}}{\thepage}}
\renewcommand{\setevenRH}{%
  \paperRunningHead{\thepage}{\theRUNAUTHOR\quad\textit{\theRUNTITLE}}}
\renewcommand{\setoddECRH}{%
  \paperRunningHead{\theECRUNAUTHOR\quad\textit{\theECRUNTITLE}\quad E-Companion}{\thepage}}
\renewcommand{\setevenECRH}{%
  \paperRunningHead{\thepage}{\theECRUNAUTHOR\quad\textit{\theECRUNTITLE}\quad E-Companion}}
\fi
\makeatother

\bibpunct[, ]{(}{)}{,}{a}{}{,}
\def\bibfont{\fontsize{10}{12}\selectfont}

\let\footnote=\endnote

\TheoremsNumberedThrough
\EquationsNumberedThrough
\ECRepeatTheorems

\makeatletter
\let\atlas@classTheoremStyle\th@TH
\def\th@TH{%
  \atlas@classTheoremStyle
  \def\@opargbegintheorem##1##2##3{\normalfont\TheoremTextFont
    \item[\if@OPRE\else\hspace*{1em}\fi\hskip\labelsep
      \theorem@headerfont ##1\ ##2\ (##3).]}}
\@ifundefined{argmax}{}{}
\@ifundefined{argmin}{\DeclareMathOperator*{\argmin}{arg\,min}}{}
\makeatother

\newcommand{\RR}{\mathbb{R}}
\newcommand{\E}{\mathbb{E}}
\newcommand{\cX}{\mathcal X}
\newcommand{\cZ}{\mathcal Z}
\DeclareMathOperator{\softmax}{softmax}

\MANUSCRIPTNO{}

\begin{document}

\RUNAUTHOR{You et al.}
\RUNTITLE{KKT Atlas}
\TITLE{KKT Atlas: Learning Low-Dimensional Local Scenario Models for Contextual Stochastic Optimization}

\ARTICLEAUTHORS{%
\AUTHOR{Zhengzhong Ricky You\textsuperscript{1}, Baichuan Mo\textsuperscript{1,*}, Yuan Zhou\textsuperscript{2}, Jinhua Zhao\textsuperscript{3}, Zuo-Jun Max Shen\textsuperscript{4}}
\AFF{\textsuperscript{1}Department of Civil Engineering, Tsinghua University, Beijing 100084, China}
\AFF{\textsuperscript{2}Yau Mathematical Sciences Center, Tsinghua University, Beijing 100084, China}
\AFF{\textsuperscript{3}Department of Urban Studies and Planning, Massachusetts Institute of Technology, Cambridge, MA 02139, USA}
\AFF{\textsuperscript{4}Department of Industrial Engineering and Operations Research, University of California, Berkeley, Berkeley, CA 94720, USA}
\AFF{\textsuperscript{*}Corresponding author: \EMAIL{bmo@tsinghua.edu.cn}}
}%

\ABSTRACT{%
In contextual stochastic optimization, integrated learning trains conditional
models on operational loss, but repeatedly solving stochastic programs with
many scenarios is costly. We propose \emph{KKT Atlas}, which partitions
contexts into local charts and identifies a low-dimensional subspace per
chart from decision gradients of the scenario cost. It generates shared
scenarios by jointly fitting subspace coordinates and
context-dependent probabilities. We show that one more scenario
than the subspace dimension suffices to approximate the conditional
Karush--Kuhn--Tucker optimality condition under regularity conditions,
yielding decision and risk bounds under local strong convexity, with
generalization error scaling with the compact parameter dimension.
On three synthetic problems, Atlas uses
4--18 scenarios per decision versus 500 for direct-solve full-support integrated
learning. It attains comparable or lower risk with about $3$--$8\times$
construction and $8$--$48\times$ batch deployment speedups. Atlas also
achieves lower risk than clustering-based reduction at matched support
budgets.
In stylized European grid capacity coverage, Atlas compresses 728
scenarios to 8--9, cutting construction and batch deployment
times by about $4\times$ relative to
full-support integrated learning.
Its loss is lowest among primary benchmarks over 11 quarters and, over 17 later quarters with fixed
settings, slightly higher than with full-support integrated learning and
significantly lower than with the remaining benchmarks.
}%


\maketitle

\section{Introduction}
\label{sec:introduction}

Operational decisions often involve uncertainty, while observed operating
conditions can provide information about that uncertainty before a decision
is made. Contextual stochastic optimization incorporates this information
into the uncertainty model used to guide the decision. For example,
inventory replenishment decisions depend on uncertain demand. An
unconditional approach estimates one demand distribution from observations
across operating conditions and uses it in the stochastic program.
Promotions and recent sales provide information about demand
before orders are placed. Conditioning on these features can represent
the demand distribution more accurately than pooling observations across
conditions.
Empirical studies show that incorporating such information can yield
better decisions than ignoring it
\citep{bertsimas2020predictive,kallus2023stochastic,kannan2025datadriven}.
The practical question is how to turn the available observations into
an effective decision rule.

A common approach to contextual stochastic optimization first estimates
conditional uncertainty from data and then solves a stochastic program
based on that estimate. For example, multinomial logistic regression fits
conditional probabilities by maximum likelihood or, equivalently, by minimizing
cross-entropy loss \citep{hastie2009elements}.
With a fixed scenario support, the fitted model gives the probability
of each scenario at a new context.
\citet{kannan2025datadriven} construct scenarios from regression predictions
and fitted residuals, with the prediction model fitted independently of
the stochastic program. In these estimate-then-optimize procedures,
fitting targets predictive or distributional accuracy, so errors
that cause costly decisions may receive little weight in training.

Integrated learning instead uses the operational loss of the induced
decisions as a training criterion for the statistical model, following
the broader principle of decision-focused learning \citep{donti2017task,mandi2024decision}.
\citet{elmachtoub2022smart} develop Smart ``Predict, then Optimize'' (SPO)
for predicting objective coefficients in optimization problems with linear
objectives. They measure prediction quality by the excess cost of the induced
decision and derive the convex SPO+ surrogate loss for training.
\citet{kallus2023stochastic} apply this principle to forests by choosing
splits that improve downstream optimization performance.
\citet{qi2026integrated} develop integrated conditional
estimation--optimization (ICEO), which fits a conditional distribution
model to minimize the operational loss of the resulting decisions.
Their numerical comparisons show that including downstream operational loss
can improve on separate distribution fitting based on cross-entropy loss.
These gains motivate accounting for operational loss when learning
a conditional distribution.

However, training with operational loss can be computationally intensive
because each loss evaluation requires solving the stochastic program induced
by the current conditional distribution model. This burden becomes substantial
when the scenario support is large or recourse is costly. Existing approaches
reduce the training cost by approximating quantities that would otherwise
require repeated optimization. For example, \citet{kallus2023stochastic} use
perturbation-based criteria to avoid reoptimization at every candidate tree
split, while \citet{qi2026integrated} learn a neural approximation of the
solution mapping from solutions computed at sampled probability vectors and use
this approximation to reduce repeated optimization during training. These approximations trade repeated exact optimization for surrogate or splitting error. We take a complementary approach and propose \emph{KKT Atlas}, which uses the scenario gradients
that enter the Karush--Kuhn--Tucker (KKT) optimality conditions to learn compact
local scenario models. Atlas solves the resulting smaller stochastic programs
directly during training and deployment.

Atlas learns from training observations $(x_i,\xi_i)$, each pairing
the context $x_i$ observed before the decision with the realization $\xi_i$
of the uncertain parameter observed afterward. It clusters the observed
contexts to partition the context space into local regions, called charts,
that group similar operating conditions. Each chart uses a shared set of
generated scenarios with context-dependent probabilities. This localization
allows the shared scenarios to focus on uncertainty within their region,
while the collection of local models accommodates differences across regions.
For a new context, Atlas solves one small stochastic program using only the
assigned chart's scenarios and context-dependent probabilities.

Within each chart, Atlas uses scenario gradients to guide the representation
of uncertainty. For a decision $z$, the scenario cost $Q(z,\xi;x)$ is the term in the
operational loss that depends on the realization $\xi$.
Its gradient $\nabla_z Q(z,\xi;x)$, which we call the scenario gradient,
determines the realization's contribution to the KKT first-order
optimality condition. Our analysis considers a local affine-gradient model
in which this gradient is approximately a linear combination of a small
number of features of $\xi$, plus a term independent of $\xi$. This
structure motivates constructing a low-dimensional basis from the
cross-covariance between realizations and scenario gradients.
Each generated scenario is represented by coordinates in this basis and
reconstructed in the original uncertainty space.

\paragraph{Contributions.}
\begin{enumerate}[label=(\roman*),leftmargin=*]
\item
\emph{We develop KKT Atlas, a local scenario framework that keeps
downstream stochastic programs small during training and deployment.}
Within each chart, Atlas constructs a low-dimensional basis from the
cross-covariance between realizations and scenario gradients.
Each chart assigns context-dependent probabilities to a small set of shared
scenarios represented by coordinates in this basis.
To limit computational cost, structure selection first compares candidate
numbers of charts by operational loss on internal assessment folds using a
simplified fitting procedure.
Joint model fitting is then performed once for the selected chart structure
by minimizing regularized empirical operational loss over scenario coordinates
and probability parameters with direct downstream solves.

\item
\emph{We establish theoretical guarantees for chart localization, minimal
shared supports, decision quality, and finite-sample generalization.}
First, without parametric assumptions on the conditional distribution, we
quantify the tradeoff between chart localization and sampling error and show
that, under suitable conditions, the dimension selected from the data recovers
the local dimension relevant to decisions with high probability.
Second, under regularity conditions, one more shared scenario than this
dimension suffices to approximate the population objective gradient and is
minimal among shared supports for matching projected conditional means that affinely span the local
subspace. Local strong convexity then converts gradient error into bounds on
decision error and operational risk.
Third, for a fixed chart structure, we derive an end-to-end excess-risk
bound whose estimation error scales with Atlas's compact parameter
dimension, without explicit dependence on reference support size.

\item
\emph{We conduct extensive computational experiments to evaluate Atlas
on three synthetic applications and in a stylized power grid study on real data.}
Across capacity allocation, conditional value-at-risk (CVaR) portfolio
optimization, and network reservation alongside rolling day-ahead grid
capacity coverage, Atlas compresses supports from 500
to 4--18 scenarios and from 728 to 8--9, respectively.
These compact models accelerate construction and batch deployment over
full-support ICEO with direct downstream solves by about $3$--$8\times$ and
$8$--$48\times$ on the synthetic problems, and by about $4\times$ on the
grid.
In decision quality, Atlas attains comparable or lower risk than ICEO on the
synthetic problems, matches ICEO closely over 11 initial grid quarters, and
incurs $1.3\%$ higher loss than ICEO over 17 later quarters with settings
fixed in advance.
Across all synthetic and grid evaluations, Atlas consistently outperforms
clustering-based reduction at matched support budgets as well as weight-
or neighborhood-based benchmarks, with significantly lower loss over those
later quarters. Ablation studies further illustrate the roles of the $d+1$
support budget, the gradient-informed local basis, and multi-chart
localization.

\end{enumerate}

The rest of the paper is organized as follows.
Section~\ref{sec:literature} reviews contextual stochastic optimization,
scenario reduction, and scenario generation. Section~\ref{sec:problem}
formulates the contextual stochastic program and introduces the training
observations and scenario gradients
used throughout the paper. In Section~\ref{sec:method}, we present KKT Atlas,
from chart construction and gradient-informed scenario coordinates to
structure selection and model fitting. Section~\ref{sec:theory}
establishes the corresponding guarantees for localization, support size,
decision quality, and finite-sample generalization. In Section~\ref{sec:computation}, we
describe the computational study and report the synthetic comparisons,
stylized European grid capacity-coverage study, and ablation and sensitivity
analyses. Section~\ref{sec:conclusion} concludes. The e-companion contains
proofs and detailed computational protocols.

\section{Related Literature}
\label{sec:literature}

In Section~\ref{subsec:contextualLiterature}, we review the use of contextual
data to learn uncertainty models and decision rules, two topics covered
in the broader survey by \citet{sadana2025survey}. We then discuss
scenario reduction in Section~\ref{subsec:reductionLiterature} and methods
for scenario generation and dimension reduction in
Section~\ref{subsec:generationLiterature}.

\subsection{Contextual stochastic optimization}
\label{subsec:contextualLiterature}

Existing contextual stochastic optimization methods either construct a context-dependent uncertainty representation for the downstream stochastic program or map contexts directly to decisions. Methods in the first group differ in both the representation they construct and the criterion used to fit it. \citet{bertsimas2020predictive} represent conditional uncertainty by assigning context-dependent weights to historical observations and solve the corresponding weighted empirical problem. These weights are constructed using predictive methods such as $k$-nearest neighbors, classification and regression trees, and random forests to estimate the conditional distribution without incorporating downstream operational loss. Subsequent tree-based methods incorporate the downstream optimization problem into the construction of these weights. \citet{stratigakos2022prescriptive} train ensembles of decision trees by directly minimizing application-specific operating costs and assign larger weights to historical observations that fall in the same terminal leaves as the current context when constructing a weighted sample average approximation (SAA). In a more general contextual stochastic optimization setting, stochastic optimization forests \citep{kallus2023stochastic} choose tree splits according to the quality of the downstream decisions they induce. \citet{ribeiro2026integrated} consider contextual two-stage linear programs and choose tree splits using an objective that balances conditional distributional fit and downstream decision quality.

Other methods construct context-dependent scenarios or probabilities. \citet{kannan2025datadriven} construct context-dependent SAA scenarios by combining a fitted prediction with its residuals, fitting the prediction model separately from the downstream optimization problem. ICEO \citep{qi2026integrated}, in contrast, uses a prespecified finite support and learns context-dependent probabilities by minimizing downstream operational loss. Different from weighted empirical, residual-based, and fixed-support integrated methods, whose downstream problems remain tied to historical observations or a prespecified reference support, Atlas learns compact chart-level shared scenario supports together with context-dependent probabilities directly under operational loss.

Methods in the second group learn a direct mapping from context to decisions.
\citet{ban2019big} minimize empirical newsvendor cost over linear rules
of the observed features and establish finite-sample performance bounds.
\citet{zhang2026piecewiseaffine} fit piecewise affine rules by minimizing
empirical operational loss, with the number of pieces controlling their
flexibility. These rules provide decisions without solving a new stochastic
program. In their constrained formulation with general nonlinear constraints,
however, feasibility at new contexts is guaranteed only with high probability.
By contrast, learning the uncertainty distribution and determining the
decision through a downstream optimization problem enforces the decision
constraints explicitly at each new context.

\subsection{Scenario reduction}
\label{subsec:reductionLiterature}
Distribution-based scenario reduction seeks a distribution that is
supported on fewer scenarios and remains close to the original distribution.
\citet{dupacova2003scenario} retain a subset of the original scenarios and
transfer the probability mass of each removed scenario to its nearest
retained scenario. Following this approach, \citet{heitsch2003scenario}
select the smaller set iteratively, adding or removing original scenarios
based on the resulting probability distance. In contrast,
\citet{rujeerapaiboon2022scenario} study continuous scenario reduction,
in which new scenarios may be chosen freely without being restricted to
the original scenario set.

Optimization-based scenario reduction uses the objective and structure of
the stochastic program to guide the choice of a smaller scenario set.
\citet{bertsimas2023optimizationbased} define a problem-dependent divergence
between distributions that accounts for the quality of the decisions they
induce and use it to construct new scenarios. \citet{zhang2023optimized}
select and reweight original scenarios to approximate the recourse function
over a pool of candidate decisions, with bounds on solution quality.
For tail-risk problems such as CVaR portfolio selection,
\citet{fairbrother2022problem} generate scenarios by prioritizing
those that can affect the tail-risk objective value. \citet{keutchayan2023problem}
cluster scenarios using the costs of decisions optimized for individual
scenarios and evaluated across the original scenario set. These methods use
the downstream optimization problem to construct a weighted scenario set
for a given distribution. By comparison, Atlas uses scenario gradients to
guide the construction of a small scenario set for each chart and learns
how the scenario probabilities vary with the observed context.

\subsection{Contextual scenario generation and dimension reduction}
\label{subsec:generationLiterature}
Constructing new scenarios adds flexibility even when the number of
scenarios is fixed. For linear two-stage stochastic programs,
\citet{henrion2022problembased} choose scenario values and probabilities to
approximate expected recourse, the expected second-stage cost for a given
first-stage decision. Their criterion minimizes the largest discrepancy
over feasible first-stage decisions between the original expected recourse
and its weighted scenario approximation. In the contextual setting,
\citet{homemdemello2026forecasting} show that, for a class of
two-stage problems with linear second stages, fixed recourse, and fixed
costs, a single suitably constructed scenario, possibly outside the
original support, can recover an optimal decision for each context.
Beyond this linear fixed-recourse setting, Atlas handles nonlinear scenario costs and uses a shared set of generated scenarios within
each chart, with probabilities that vary by context.

Contextual scenario generators learn scenario sets from observed data.
\citet{islip2026contextual} map each context to a prescribed number of
scenarios, while \citet{hornewall2025decision} train a neural scenario
generator through a deterministic equivalent regularized by a log barrier,
using implicit differentiation of the KKT conditions
\citep{donti2017task,agrawal2019differentiable}.
Both formulations specify generated
scenarios in the original uncertainty space. In the latter method, KKT
differentiation supplies gradients for updating the generator parameters.
Atlas combines KKT differentiation during fitting with a preceding
representation step that uses scenario gradients to select a
low-dimensional basis for each chart. The generated scenarios are fitted
through coordinates in this basis, reducing the number of scenario parameters
to learn.

Even with a small support, learning scenarios in the full uncertainty
space can be costly because each scenario has many components.
One way to address this dimensional burden is to identify a few input
directions along which a function varies most and use the resulting
reduced representation. For example, \citet{zahm2020gradient} use gradients
to identify these directions and approximate a vector-valued function through
the corresponding linear combinations of its inputs. \citet{chen2025coupled}
also reduce the output space, jointly selecting input and output subspaces to
account for how the choice in each space affects the relevant directions in the
other. Both studies select reduced spaces using approximation bounds based on
derivatives of the function with respect to its inputs. In contrast,
Atlas uses the cross-covariance between realizations and decision gradients
of the scenario cost, without differentiating with respect to the uncertain
parameter, to select a representation suited to
downstream first-order optimality conditions.
In moment-based distributionally robust optimization,
\citet{jiang2026optimized} optimize a linear projection of the uncertain
parameters within reduced-dimensional approximations of the original
problem. Atlas uses the estimated subspace in each chart to parameterize
a small set of scenarios and jointly fits the scenarios and their
context-dependent probabilities using the operational loss of the
resulting decisions.

\section{Problem Setting and Preliminaries}
\label{sec:problem}

Section~\ref{subsec:contextualProblem} formulates the contextual stochastic
program and its finite-scenario approximation.
Section~\ref{subsec:scenarioGradients} then identifies the
distribution-dependent term in the optimality condition.
Throughout, we write $\RR_+:=[0,\infty)$ and use
$\mathbf 1_s\in\RR^s$ for the all-ones vector, omitting the subscript
when the dimension is clear. The symbols $\|\cdot\|_2$,
$\|\cdot\|_F$, and $\|\cdot\|_{\rm op}$ denote the Euclidean,
Frobenius, and operator norms, respectively.
Table~\ref{tab:notationSummary} in the e-companion summarizes the notation.

\subsection{Contextual stochastic program}
\label{subsec:contextualProblem}

Before making a decision, the decision maker observes a vector $x$ in the
context space $\cX\subseteq\RR^{d_x}$. The observed context provides information about the
uncertain parameter $\xi$, whose value in the uncertainty space
$\Xi\subseteq\RR^{d_\xi}$ remains unknown.
For example, in inventory replenishment, a promotion known before orders
are placed may signal higher future demand. The promotion is part of $x$,
while $\xi$ represents the demand realized after the order is placed.
A decision $z$ chosen from the feasible decision set
$\cZ\subseteq\RR^{d_z}$ under context $x$
incurs the following operational loss when $\xi$ is realized:
\begin{equation}
 \ell(z,\xi;x):=f(z;x)+Q(z,\xi;x),
 \label{eq:operationalLoss}
\end{equation}
where $f(z;x)$ is the deterministic cost and $Q(z,\xi;x)$ is the
scenario cost, which depends on $\xi$.

Given an observed context $x$, let $P_x$ denote the conditional
distribution of $\xi$. If $P_x$ were known, the decision maker would
minimize the conditional operational risk $R_0(z;x)$ by solving
\begin{equation}
 \min_{z\in\cZ}\;\bigl\{R_0(z;x):=\E_{P_x}[\ell(z,\xi;x)]\bigr\}.
 \label{eq:populationObjectives}
\end{equation}
In practice, $P_x$ is unknown, so we cannot solve
Eq.~\eqref{eq:populationObjectives} directly. Instead, we use training
observations $\{(x_i,\xi_i)\}_{i\in I_{\rm tr}}$, each pairing an observed
context with the corresponding realization of the uncertain parameter, to learn a mapping
$x\mapsto\widetilde P_x$. For each new context $x$, the finite
distribution $\widetilde P_x$ assigns probabilities
$p=(p_1,\ldots,p_K)\in\Delta_K:=\{p\in\RR_+^K:\mathbf 1^\top p=1\}$
to $K$ points
$\tilde{\boldsymbol{\xi}}=(\tilde{\xi}_1,\ldots,\tilde{\xi}_K)\in\Xi^K$,
called scenarios. The resulting optimization problem is
$\min_{z\in\cZ}\sum_{k=1}^K p_k\ell(z,\tilde{\xi}_k;x)$.
Our goal is to learn a mapping $x\mapsto\widetilde P_x$ with few scenarios
whose decisions nearly minimize Eq.~\eqref{eq:populationObjectives}.

When the population objective $R_0(z;x)$ in Eq.~\eqref{eq:populationObjectives}
has little curvature near its minimizer, even a small error in approximating
its gradient can substantially change the optimal decision. To stabilize the decision,
we introduce quadratic regularization. Let $z_0:\cX\to\cZ$ be a
prespecified regularization center and $\eta_{\rm reg}\geq0$ the
regularization coefficient. For a probability distribution $P$ on $\Xi$,
we define the regularized objective at context $x$ and, when the minimum
is attained, an optimal decision by
\begin{equation}
 R(z;P,x):=\E_P[\ell(z,\xi;x)]
 +\frac{\eta_{\rm reg}}{2}\bigl\|z-z_0(x)\bigr\|_2^2,\qquad
 z(P;x)\in\argmin_{z\in\cZ}R(z;P,x).
 \label{eq:finiteDecisionLayer}
\end{equation}
The population and finite-scenario decisions are $z(P_x;x)$ and
$z(\widetilde P_x;x)$, respectively. The following assumption specifies
a curvature condition on the population objective $R(\cdot;P_x,x)$.

\begin{assumption}[Strong convexity]
\label{ass:regularizedProblem}
The feasible set $\cZ$ is nonempty, closed, and convex. For every
$x\in\cX$ and $\xi\in\Xi$, the functions $f(\cdot;x)$ and
$Q(\cdot,\xi;x)$ are proper, closed, and convex.
For every $x\in\cX$, $\ell(z,\xi;x)$ is integrable under $P_x$ for all
$z\in\cZ$, and the population objective $R(\cdot;P_x,x)$ attains its minimum on $\cZ$.
Every finite-scenario problem in Eq.~\eqref{eq:finiteDecisionLayer}
has a unique minimizer on $\cZ$.

There exist constants $\eta>0$ and $r_{\rm loc}\in(0,\infty]$
such that, for every $x\in\cX$,
some minimizer $z(P_x;x)$ of $R(\cdot;P_x,x)$ on $\cZ$ satisfies, for all
$z,z'\in\mathcal B_{\cZ}(z(P_x;x),r_{\rm loc})$ and $t\in[0,1]$,
\begin{equation}
 R(tz+(1-t)z';P_x,x)\leq tR(z;P_x,x)+(1-t)R(z';P_x,x)
 -\frac{\eta}{2}t(1-t)\|z-z'\|_2^2.
 \label{eq:strongConvexity}
\end{equation}
Here
$\mathcal B_{\cZ}(z(P_x;x),r_{\rm loc})
:=\{z\in\cZ:\|z-z(P_x;x)\|_2\leq r_{\rm loc}\}$
is the neighborhood of $z(P_x;x)$.
\end{assumption}

Assumption~\ref{ass:regularizedProblem} ensures unique optimal decisions
and provides the strong-convexity constant $\eta$ used to convert
first-order approximation error into decision and risk bounds.
We impose strong convexity on a neighborhood of $z(P_x;x)$ instead
of requiring it globally for two reasons.
First, when $\eta_{\rm reg}$ is zero or small, $R_0(\cdot;x)$ may have
stronger curvature near $z(P_x;x)$ than globally on $\cZ$, allowing
us to use a larger local constant $\eta$.
Second, when $\cZ$ is unbounded, a finite $r_{\rm loc}$ restricts
the required uniform gradient error bound to the bounded neighborhood
$\mathcal B_{\cZ}(z(P_x;x),r_{\rm loc})$.
Theorem~\ref{thm:fixedAtlasRisk} uses these properties to bound the
expected excess operational risk of deployed decisions, with the
gradient error term decreasing as $\eta$ increases for a fixed
gradient error bound.

\subsection{Scenario gradients and optimality}
\label{subsec:scenarioGradients}

To reduce the number of scenarios while preserving the information needed
for decisions, we examine how the uncertain parameter $\xi$ affects the
KKT optimality conditions through the derivative of $Q(z,\xi;x)$ with
respect to the decision $z$. We rely on the following differentiability assumption.

\begin{assumption}[Differentiability]
\label{ass:differentiableScenarioCosts}
For every $x\in\cX$ and $\xi\in\Xi$, both $f(\cdot;x)$ and $Q(\cdot,\xi;x)$
are differentiable on a neighborhood of $\cZ$. For every $z\in\cZ$ and $x\in\cX$, the
gradient of $Q$ is integrable under $P_x$, and differentiation in $z$
commutes with expectation under $P_x$.
\end{assumption}

Assumption~\ref{ass:differentiableScenarioCosts} allows us to express
the gradient of the expected scenario cost as the conditional mean of
scenario gradients. Specifically, define the scenario gradient
$\phi(z,\xi;x):=\nabla_z Q(z,\xi;x)\in\RR^{d_z}$ and its conditional mean
$\bar\phi(z;x):=\E_{P_x}[\phi(z,\xi;x)]$.
Interchanging differentiation and conditional expectation gives the
first-order optimality condition
\citep[Section~4.2.3]{boyd2004convex}:
\begin{equation}
\begin{aligned}
 0&\in \nabla_z R(z(P_x;x);P_x,x)+N_{\cZ}(z(P_x;x)),\\
 \nabla_z R(z;P_x,x)
 &=\nabla_z f(z;x)+\bar\phi(z;x)
 +\eta_{\rm reg}\{z-z_0(x)\},
\end{aligned}
 \label{eq:populationFirstOrderMapping}
\end{equation}
where
$N_{\cZ}(z):=\{v\in\RR^{d_z}:v^\top(u-z)\leq0\text{ for every }u\in\cZ\}$
denotes the normal cone of $\cZ$ at $z$. The mean scenario gradient $\bar\phi(z;x)$ is the only term in the
optimality condition of Eq.~\eqref{eq:populationFirstOrderMapping}
that depends on the conditional distribution. Replacing $P_x$ with a finite scenario
distribution $\widetilde P_x$ leaves the deterministic cost, regularizer, and constraints
unchanged, while replacing $\bar\phi(z;x)$ with
$\sum_{k=1}^K p_k\phi(z,\tilde{\xi}_k;x)$.

In Section~\ref{sec:method}, we use scenario gradients $\phi(z,\xi_i;x)$ evaluated at
training realizations to construct a cross-covariance basis.
Section~\ref{sec:theory} then bounds how errors in approximating
$\bar\phi$ affect the optimal decision of the finite-scenario problem.

\section{Construction of KKT Atlas}
\label{sec:method}

KKT Atlas combines two ideas. First, it partitions the context space into charts so that each chart only needs to represent uncertainty within its local region, allowing a smaller scenario support for each decision. Second, within each chart, Atlas uses scenario gradients $\phi(z,\xi_i;x)$ to identify a low-dimensional representation of uncertainty relevant to downstream decisions, so the chart can use fewer scenarios per decision.

We describe the context charts and local scenario directions in
Sections~\ref{subsec:contextCharts} and \ref{subsec:kktCoordinates},
construct local scenarios and their context-dependent
probabilities in Section~\ref{subsec:localModel}, and detail joint fitting and
structure selection in Sections~\ref{subsec:modelFitting} and
\ref{subsec:integratedTraining}.

\subsection{Context charts from K-means}
\label{subsec:contextCharts}

Let $\{(x_i,\xi_i)\}_{i\in I_{\rm tr}}$ denote the training observations,
and let $\widehat D_x\in\RR^{d_x\times d_x}$ be the diagonal matrix of
sample standard deviations of the training contexts, with all diagonal
entries assumed positive. For $h\in\RR^{d_x}$, define the standardized
norm $\|h\|_x:=\|\widehat D_x^{-1}h\|_2$.

For a given number $B$ of charts, we apply K-means to the standardized
training contexts and transform the resulting cluster centroids back to
the original scale. Let $\bar x_b$ denote the centroid of chart $b$.
We assign each context $x$ to its nearest centroid under the standardized
norm, so that
$\tau_B(x)\in\argmin_{1\leq b\leq B}\|x-\bar x_b\|_x^2$,
with ties broken in favor of the smallest chart index.
Let $I_b:=\{i\in I_{\rm tr}:\tau_B(x_i)=b\}$ denote the indices of
training observations assigned to chart $b$.
Since $\cX$ need not be convex, the chart centroids $\bar x_b$ may lie outside
$\cX$. We therefore enlarge the context space to
$\cX\cup\{\bar x_b:b=1,\ldots,B\}$ and, with a slight abuse of notation,
continue to denote it by $\cX$ throughout the subsequent analysis.
We assume that $P_x$, $f$, $Q$, and the regularization center $z_0$
extend to these centroids, and that Assumptions~\ref{ass:regularizedProblem}
and~\ref{ass:differentiableScenarioCosts} hold there with the same
curvature parameters $\eta$ and $r_{\rm loc}$.

\subsection{Scenario directions from scenario gradients}
\label{subsec:kktCoordinates}

For each chart, we compute the local SAA decision
$z_b^{\rm SAA}$ by solving
\begin{equation}
 \min_{z\in\cZ}\;
 f(z;\bar x_b)
 +\frac{1}{|I_b|}\sum_{i\in I_b}Q(z,\xi_i;\bar x_b)
 +\frac{\eta_{\rm reg}}{2}\|z-z_0(\bar x_b)\|_2^2.
 \label{eq:localSaaProblem}
\end{equation}
At the solution of Eq.~\eqref{eq:localSaaProblem}, we evaluate
$\phi(z_b^{\rm SAA},\xi_i;\bar x_b)$ for each $i\in I_b$.
Fix a transformation $T:\Xi\to\RR^{d_\xi}$ and a decoder
$g:\operatorname{dom}(g)\to\Xi$, where
$\operatorname{dom}(g)\subseteq\RR^{d_\xi}$ denotes the domain of $g$. We need these maps to ensure that generated
scenarios belong to the uncertainty space $\Xi$. The decoder $g$ enforces
this requirement, while $T$ reverses decoding for realizations,
satisfying $g(T(\xi))=\xi$ on the support of $\xi$.
For example, in an inventory problem with nonnegative demand, we may
choose $g(y)=y^2$ and $T(\xi)=\sqrt{\xi}$, applied componentwise.
Decoding then produces nonnegative demands and recovers each observed
demand from its encoded coordinates. Without this decoding step, generated
coordinates may be negative and invalid as inputs to the
stochastic program.

We center and scale the transformed
realizations separately within each chart with $|I_b|\geq 2$. Let $\widehat\mu_b:=|I_b|^{-1}\sum_{i\in I_b}T(\xi_i)$. For a small
numerical tolerance $\varepsilon_{\rm num}>0$, such as $10^{-8}$, define
the diagonal scale matrix $\widehat D_b\in\RR^{d_\xi\times d_\xi}$
componentwise by
\begin{equation}
 (\widehat D_b)_{rr}
 :=
 \max\left\{
 \left(
 \frac{1}{|I_b|-1}\sum_{i\in I_b}
 \bigl([T(\xi_i)]_r-(\widehat\mu_b)_r\bigr)^2
 \right)^{1/2},
 \varepsilon_{\rm num}
 \right\},
 \qquad r=1,\ldots,d_\xi.
\label{eq:chartScaleDefinition}
\end{equation}
Using the mean and scale in Eq.~\eqref{eq:chartScaleDefinition}, define the
standardized encoder
$T_b(\xi):=\widehat D_b^{-1}\{T(\xi)-\widehat\mu_b\}$ and decoder
$g_b(y):=g(\widehat\mu_b+\widehat D_b y)$.
Here $T_b:\Xi\to\RR^{d_\xi}$ and
$\operatorname{dom}(g_b)=
\{y\in\RR^{d_\xi}:\widehat\mu_b+\widehat D_b y\in\operatorname{dom}(g)\}$.
Thus $T_b$ and $g_b$ combine $T$ and $g$ with chart-specific standardization
and its inverse, satisfying $g_b(T_b(\xi))=\xi$.

Since $|I_b|^{-1}\sum_{j\in I_b}T_b(\xi_j)=0$ by construction,
the empirical cross-covariance simplifies to
\begin{align}
 \widehat C_b
 &:=\frac{1}{|I_b|}\sum_{i\in I_b}
 \left(T_b(\xi_i)-\frac{1}{|I_b|}\sum_{j\in I_b}T_b(\xi_j)\right)
 \left(\phi(z_b^{\rm SAA},\xi_i;\bar x_b)
 -\frac{1}{|I_b|}\sum_{j\in I_b}\phi(z_b^{\rm SAA},\xi_j;\bar x_b)\right)^\top
 \notag\\
 &=\frac{1}{|I_b|}\sum_{i\in I_b}
 T_b(\xi_i)\phi(z_b^{\rm SAA},\xi_i;\bar x_b)^\top
 \in\RR^{d_\xi\times d_z}.
 \label{eq:gradientCrossCovariance}
\end{align}
Let $\sigma_j(\widehat C_b)$ denote the $j$th largest singular value
of $\widehat C_b$.
For a fixed retention proportion $\alpha_{\rm ret}\in(0,1)$, choose
\begin{equation}
 d_b:=\min\left\{d\in\{1,\ldots,\min(d_\xi,d_z)\}:
 \sum_{j=1}^d\sigma_j(\widehat C_b)^2
 \geq\alpha_{\rm ret}\|\widehat C_b\|_F^2\right\}.
 \label{eq:spectralWidthRule}
\end{equation}
Note from Eq.~\eqref{eq:spectralWidthRule} that $d_b\leq\min(d_\xi,d_z)$, so the retained dimension is automatically bounded by the decision dimension whenever $d_z<d_\xi$.

\subsection{Local scenarios and probabilities}
\label{subsec:localModel}

Let $\widehat V_b\in\RR^{d_\xi\times d_b}$ be the local basis for chart $b$,
whose columns are the leading $d_b$ left singular vectors of $\widehat C_b$.
These vectors define directions in the standardized coordinate space $\RR^{d_\xi}$.
Intuitively, we retain these directions because they capture the strongest
linear associations between transformed realizations $T_b(\xi_i)$ and scenario gradients $\phi(z_b^{\rm SAA},\xi_i;\bar x_b)$
in the sample, focusing the reduced representation on uncertainty relevant
to decisions. Projecting $T_b(\xi)\in\RR^{d_\xi}$ onto these directions
gives $\widehat V_b^\top T_b(\xi)\in\RR^{d_b}$, representing each realization
with $d_b$ coordinates instead of $d_\xi$, thereby reducing the dimension. Conversely, any coordinate vector $u\in\RR^{d_b}$ with
$\widehat V_b u\in\operatorname{dom}(g_b)$ decodes to a feasible scenario
$\tilde{\xi}=g_b(\widehat V_b u)$.

Motivated by Theorem~\ref{thm:tiltedResponseCapacity}, we use
$K_b=d_b+1$ trainable scenario coordinate vectors
$u_{b,1},\ldots,u_{b,K_b}\in\RR^{d_b}$ in chart $b$ and define
$\tilde{\xi}_{b,k}:=g_b(\widehat V_bu_{b,k})$.
With probability parameters consisting of trainable intercepts $a_b\in\RR^{K_b}$ and context coefficients
$\Gamma_b\in\RR^{K_b\times d_x}$, for a context $x$ in chart $b$, the
probability model for these scenarios is
\begin{equation}
 p_b(x;a_b,\Gamma_b):=
 \softmax\!\left(a_b+\Gamma_b\widehat D_x^{-1}(x-\bar x_b)\right)
 \in\Delta_{K_b},
 \label{eq:localProbabilityModel}
\end{equation}
where \(\softmax(v)_k:=\exp(v_k)/\sum_{j=1}^{K_b}\exp(v_j)\) and
$\Delta_{K_b}$ is the probability simplex. We write
\(p_{b,k}(x;a_b,\Gamma_b)\) for component \(k\), the probability assigned
to $\tilde{\xi}_{b,k}$, and abbreviate it as
\(p_{b,k}(x)\) when the parameters are clear.

Contexts within a chart share the same scenarios, whose probabilities vary
with context. The model has $K_bd_b$ scenario coordinates and $K_b$
scenario terms in its stochastic program. With
$\tilde{\boldsymbol{\xi}}_b:=(\tilde{\xi}_{b,1},\ldots,\tilde{\xi}_{b,K_b})$ and
$\theta_b:=(u_{b,1},\ldots,u_{b,K_b},a_b,\Gamma_b)$,
let $\widetilde P_{b,x}(\theta_b)$ denote the distribution specified by
these scenarios and their probabilities $p_b(x)$. We use
Eq.~\eqref{eq:finiteDecisionLayer} to define the local objective
and its induced decision as
\begin{equation}
 \widehat R_b(z;x,\theta_b):=R\bigl(z;\widetilde P_{b,x}(\theta_b),x\bigr),
 \qquad
 z_b(x;\theta_b):=z\bigl(\widetilde P_{b,x}(\theta_b);x\bigr).
 \label{eq:localChartProblem}
\end{equation}
In the next subsection, we jointly fit the scenario coordinates $u_{b,k}$ and probability
parameters $(a_b,\Gamma_b)$ in $\theta_b$, keeping $\widehat V_b$ and $g_b$ fixed.
We evaluate operational loss using the induced decision $z_b(x;\theta_b)$.

\subsection{Model fitting}
\label{subsec:modelFitting}

For a fixed partition with $B$ charts, consider chart $b$. Let
$I^{\rm fit}\subseteq I_b$ index the observations used for
initialization; model fitting uses $I^{\rm fit}=I_b$, whereas
structure selection in Section~\ref{subsec:integratedTraining} applies
the same initialization to fold subsets $I^{\rm fit}\subset I_b$. We apply
K-means to
$\{\widehat V_b^\top T_b(\xi_i):i\in I^{\rm fit}\}$ and initialize the
scenario coordinates $(u_{b,k})_{k=1}^{K_b}$ at the resulting $K_b$ cluster centroids. For
$i\in I^{\rm fit}$, let $k(i)\in\{1,\ldots,K_b\}$ index the cluster
containing the projected realization $\widehat V_b^\top T_b(\xi_i)$. Holding the scenario
coordinates fixed, we use these cluster assignments to initialize
the affine-softmax probability model $p_b(x;a_b,\Gamma_b)$ by solving
\begin{equation}
 \min_{\substack{a\in\RR^{K_b},\,\Gamma\in\RR^{K_b\times d_x}\\
                 \mathbf 1^\top a=0,\,\mathbf 1^\top\Gamma=0}}\;
 -\frac{1}{|I^{\rm fit}|}
 \sum_{i\in I^{\rm fit}}
 \log p_{b,k(i)}(x_i;a,\Gamma)
 +\lambda_\Gamma\|\Gamma\|_F^2.
 \label{eq:screeningProbabilityFit}
\end{equation}
Here \(\lambda_\Gamma>0\) penalizes the context coefficients $\Gamma$, and $p_{b,k(i)}(x_i;a,\Gamma)$ is the probability that the
model in Eq.~\eqref{eq:localProbabilityModel} assigns to cluster $k(i)$.
The cross-entropy term fits the conditional cluster probabilities,
and the constraints remove logit shift indeterminacy.
Evaluating the objective and its gradient costs
$O(|I^{\rm fit}|K_bd_x)$ and requires no downstream optimization.

With the partition, local bases, dimensions, support sizes, and decoder
fixed, we jointly fit the scenario coordinates $u_{b,k}$ and probability parameters $(a_b,\Gamma_b)$.
For a fixed penalty weight $\lambda_u\geq0$ on the scenario coordinates and
the penalty weight $\lambda_\Gamma$ from
Eq.~\eqref{eq:screeningProbabilityFit}, define
$\Omega_b(\theta_b):=(\lambda_u/K_b)\sum_{k=1}^{K_b}\|u_{b,k}\|_2^2
+\lambda_\Gamma\|\Gamma_b\|_F^2$.
We then minimize the regularized empirical operational loss
over $\theta_b$ with $\widehat V_b u_{b,k}\in\operatorname{dom}(g_b)$
for all $k$ and the zero-sum logit constraints of
Eq.~\eqref{eq:screeningProbabilityFit} by solving
\begin{equation}
 \min_{\substack{\theta_b\\ \mathbf 1^\top a_b=0,\,\mathbf 1^\top\Gamma_b=0}}\;
 \frac{1}{|I_b|}\sum_{i\in I_b}
 \ell(z_b(x_i;\theta_b),\xi_i;x_i)
 +\Omega_b(\theta_b).
 \label{eq:integratedTrainingObjective}
\end{equation}

We solve Eq.~\eqref{eq:integratedTrainingObjective} using first-order
methods \citep{kingma2015adam}. Following standard implicit-differentiation
methods in decision-focused learning
\citep{fiacco1976sensitivity,bonnans2000perturbation,amos2017optnet,donti2017task,agrawal2019differentiable,blondel2022efficient},
we obtain gradients by differentiating through the KKT conditions of the
$K_b$-scenario problem in Eq.~\eqref{eq:localChartProblem}.
We then propagate these gradients through $p_b(x;a_b,\Gamma_b)$ in
Eq.~\eqref{eq:localProbabilityModel} and
$\tilde{\xi}_{b,k}=g_b(\widehat V_bu_{b,k})$ to jointly update the
probability parameters $(a_b,\Gamma_b)$ and scenario coordinates $u_{b,k}$.

\paragraph{Remark.}
The term ``KKT'' in \emph{KKT Atlas} refers to the downstream
optimality condition associated with the first-order mapping $\nabla_z R(z;P_x,x)$ in
Eq.~\eqref{eq:populationFirstOrderMapping}, which links
context-dependent uncertainty to the optimal decision.
This role is
distinct from the KKT-based implicit differentiation described
above, which we use to optimize
Eq.~\eqref{eq:integratedTrainingObjective} with respect to $\theta_b$.
Alternative training approaches, including smoothing, surrogate losses,
and solution-map surrogates
\citep{wilder2019melding,berthet2020learning,elmachtoub2022smart,qi2026integrated},
can be used without changing the KKT Atlas representation.

\subsection{Structure selection}
\label{subsec:integratedTraining}

Structure selection determines the number of charts
$B$. As discussed in Section~\ref{sec:theory}, more charts can reduce
within-chart variation but leave fewer observations for estimating each
local model, so there is no obvious choice of $B$. We therefore compare
a fixed grid of candidate values of $B$ using a simplified fitting
procedure to limit computational cost. For each candidate, the
preceding construction determines a partition and the local dimension
$d_b$ of each chart.

We screen candidate partitions using chart sample sizes $|I_b|$, numerical ranks,
the number of distinct projected realizations $\widehat V_b^\top T_b(\xi_i)$, and the number of fitting
observations per parameter, as specified in
Section~\ref{ec:commonStructureProtocol} of the e-companion.
Let $\mathcal G_B$ contain the candidate values of $B$ that pass these checks.
An empty $\mathcal G_B$ indicates that the training sample is insufficient to
identify a reliable local representation, and no candidate is selected. In all
our experiments, $\mathcal G_B$ is nonempty.

We hold the candidate partition, chart standardizations, and local bases
$\widehat V_b$ fixed as computed on $I_{\rm tr}$, so all assessment folds share
the same chart representations. For each eligible candidate
partition, chart $b$, and assessment fold, let
$I^{\rm fit}\subset I_b$ contain the observations outside that fold;
$I_b\setminus I^{\rm fit}$ is the assessment set. We initialize the
scenario coordinates $u_{b,k}$ using these fitting observations as described
in Section~\ref{subsec:modelFitting}. Holding these coordinates fixed, we fit
the probability parameters $(a_b,\Gamma_b)$ to the resulting cluster assignments using the
cross-entropy objective in Eq.~\eqref{eq:screeningProbabilityFit}.
We solve the downstream problems only when evaluating the initialized model on the
assessment observations of that fold. We aggregate the resulting operational
losses across charts, giving each observation equal weight, and select the
partition with the lowest aggregate loss. Let $\widehat B$ denote the selected
number of charts, which fixes the chart structure
$(\tau_{\widehat B},\widehat D_x,\{\bar x_b,\widehat V_b,g_b,d_b,K_b\}_{b=1}^{\widehat B})$
with $K_b=d_b+1$. For the subsequent model fitting, all chart-specific
quantities refer to that selected chart structure.

With the selected chart structure fixed, we reinitialize each local model using
all training observations in its chart, indexed by $I_b$. We then jointly
fit the scenario coordinates $u_{b,k}$ and probability parameters $(a_b,\Gamma_b)$ by
minimizing the regularized empirical operational loss in
Eq.~\eqref{eq:integratedTrainingObjective}, as described in
Section~\ref{subsec:modelFitting}. This model fitting is performed only
once for the selected chart structure.
At deployment, we assign a new context $x$ to chart
$b=\tau_{\widehat B}(x)$, evaluate its $K_b$ fitted probabilities using
Eq.~\eqref{eq:localProbabilityModel}, and solve one $K_b$-scenario problem
in Eq.~\eqref{eq:localChartProblem}.
Algorithm~\ref{alg:kktAtlas} summarizes the construction.

\begin{algorithm}[htbp]
\caption{Construction of KKT Atlas}
\label{alg:kktAtlas}
\setlength{\linewidth}{\hsize}
\setlength{\columnwidth}{\linewidth}
\small
\linespread{1.}\selectfont
\DontPrintSemicolon
\KwIn{Training and validation observations $\{(x_i,\xi_i)\}_{i\in I_{\rm tr}}$ and $\{(x_i,\xi_i)\}_{i\in I_{\rm val}}$; candidate grid of chart numbers.}

\textit{Step 1. Structure selection.}\;
Standardize the training contexts using their sample means and standard deviations\;
\For{each candidate $B$ in the grid}{
  Partition the standardized contexts into $B$ charts using K-means; if $\min_b |I_b|<n_{\min}$, skip $B$\;
  \For{$b=1,\ldots,B$}{
    Solve the local SAA problem in Eq.~\eqref{eq:localSaaProblem} and compute $\widehat C_b$ using Eq.~\eqref{eq:gradientCrossCovariance}\;
    Select $d_b$ using Eq.~\eqref{eq:spectralWidthRule}, construct the local basis $\widehat V_b$, and set $K_b=d_b+1$\;
  }
  If $B$ passes screening, add $B$ to $\mathcal G_B$; holding the partition and bases $\widehat V_b$ fixed, for each chart $b$ and assessment fold, initialize $u_{b,k}$ by K-means on $I^{\rm fit}$, fit probabilities using Eq.~\eqref{eq:screeningProbabilityFit}, and evaluate operational loss on $I_b\setminus I^{\rm fit}$\;
}
Select $\widehat B\in\mathcal G_B$ with the lowest aggregate assessment loss\;

\textit{Step 2. Model fitting.}\;
\For{$b=1,\ldots,\widehat B$}{
  Initialize scenario coordinates $u_{b,k}$ using all observations in $I_b$ and probabilities using Eq.~\eqref{eq:screeningProbabilityFit}\;
}
Jointly fit $\theta_1,\ldots,\theta_{\widehat B}$ by minimizing Eq.~\eqref{eq:integratedTrainingObjective} and select the model using operational loss on $I_{\rm val}$\;
\KwOut{Selected charts, local bases $\widehat V_b$, and fitted scenario models.}
\end{algorithm}
\FloatBarrier

\section{Theoretical Guarantees}
\label{sec:theory}

We analyze the structural choices introduced in Section~\ref{sec:method}.
Section~\ref{subsec:chartTradeoff} examines the tradeoff between localization
and estimation as the number of charts $B$ varies.
Section~\ref{subsec:sharedScenarioSupport} establishes conditions under which
the retention rule in Eq.~\eqref{eq:spectralWidthRule} recovers the
structural dimension and shows when $K_b=d_b+1$ shared scenarios are
sufficient and minimal. Section~\ref{subsec:operatorGuarantee} bounds
decision error and operational risk for a fixed fitted model.
Finally, Section~\ref{subsec:finiteSampleGeneralization} bounds the
finite-sample generalization error from model fitting and derives
an end-to-end excess-risk bound.

\subsection{Chart localization}
\label{subsec:chartTradeoff}

Increasing $B$ permits finer localization but reduces the average number of
training observations available to each local model. We quantify this
tradeoff by comparing the empirical cross-covariance $\widehat C_b$ in
Eq.~\eqref{eq:gradientCrossCovariance} with the population target signal
$C_b^\star$ in Eq.~\eqref{eq:idealSignalDefinition} below, evaluated at
$(z_b^{\rm SAA},\widehat D_b)$.
Because scenario gradients $\phi(z,\xi;x)$ can be nonlinear, we first decompose the scenario
gradient in each chart into a low-dimensional linear part and an
approximation residual. The linear part defines the ideal target matrix that
$\widehat C_b$ aims to recover, whereas the residual and within-chart moment
differences capture the error from pooling different contexts into one chart.
The tradeoff in $B$ balances these within-chart errors against sampling noise.

\begin{definition}[Affine-gradient decomposition]
\label{def:affineGradientDecomposition}
For chart $b$, let
$\cX_b:=\{x\in\cX:\tau_B(x)=b\}\cup\{\bar x_b\}$
denote the chart context region and
$\cZ_b:=\{z_b^{\rm SAA}\}\cup
\bigcup_{x\in\cX_b}\mathcal B_{\cZ}(z(P_x;x),r_{\rm loc})$
the local decision set. Given an integer $1\leq d_b^\star\leq\min\{d_\xi,d_z\}$,
a matrix $V_b^{\rm rel}\in\RR^{d_\xi\times d_b^\star}$ with orthonormal
columns,
and functions
$\beta_b:\cZ_b\times\cX_b\to\RR^{d_z}$ and
$H_b:\cZ_b\times\cX_b\to\RR^{d_z\times d_b^\star}$,
we define the residual function
$r_b:\cZ_b\times\operatorname{dom}(g_b)\times\cX_b\to\RR^{d_z}$
by
\begin{equation}
 \phi(z,g_b(y);x)
 =\beta_b(z,x)+H_b(z,x)(V_b^{\rm rel})^\top y
 +r_b(z,y,x),
 \label{eq:targetAffineResponse}
\end{equation}
for all $z\in\cZ_b$, $y\in\operatorname{dom}(g_b)$,
and $x\in\cX_b$.
\end{definition}

The following assumption bounds $H_b$ and the residual $r_b$ in
Definition~\ref{def:affineGradientDecomposition}.

\begin{assumption}[Uniform bounds]
\label{ass:affineUniformBounds}
For every $x\in\cX_b$, assume that $T_b(\xi)$ is integrable under $P_x$
and satisfies $g_b(T_b(\xi))=\xi$ almost surely under $P_x$.
There exist a radius $r_{y,b}\in(0,\infty]$ and finite constants
$\overline H_b,\varepsilon_b^{\rm aff}\geq0$ such that
$\{y\in\RR^{d_\xi}:\|y\|_2\leq r_{y,b}\}
\subseteq\operatorname{dom}(g_b)$ and, for all $z\in\cZ_b$, $x\in\cX_b$,
and $\|y\|_2\leq r_{y,b}$,
\begin{equation}
 \|H_b(z,x)\|_{\rm op}\leq \overline H_b,
 \qquad
 \bigl\|\E_{P_x}[r_b(z,T_b(\xi),x)]\bigr\|_2\leq \varepsilon_b^{\rm aff},
 \qquad
 \|r_b(z,y,x)\|_2\leq \varepsilon_b^{\rm aff}.
 \label{eq:affineResidualBound}
\end{equation}
\end{assumption}

To bound the estimation error of $\widehat C_b$, we impose the following
finite-sample conditions. Since $\widehat C_b$ averages observations from
different contexts, we require conditional independence, bounded transformed
realizations, and conditional means and covariances that vary smoothly with
the context.

\begin{assumption}[Finite-sample covariance conditions]
\label{ass:finiteSampleCovariance}
Let $\mathcal F_{x,\rm tr}$ denote the sigma-algebra generated by the
training contexts $\{x_i:i\in I_{\rm tr}\}$. For a candidate partition
with $B$ charts, assume the following conditions hold conditional on
$\mathcal F_{x,\rm tr}$.
\begin{enumerate}[label=(\roman*),ref=\roman*]
\item\label{cond:chartwiseIndependence}
The realizations $\xi_i$ are independent with conditional distributions $P_{x_i}$.

\item\label{cond:chartwiseAffineRepresentation}
Assumption~\ref{ass:affineUniformBounds} holds almost surely for the
decomposition in Definition~\ref{def:affineGradientDecomposition},
with fixed dimensions $d_b^\star$, coordinate radii $r_{y,b}$,
slope bounds $\overline H_b$, and residual bounds $\varepsilon_b^{\rm aff}$.

\item\label{cond:chartwiseSupportEnvelope}
There exists a finite constant $R_T\geq0$ such that, for every chart $b$
and $x\in\{x_i:i\in I_b\}\cup\{\bar x_b\}$, $\|T(\xi)\|_2\leq R_T$ and
$g(T(\xi))=\xi$ almost surely under $P_x$.

\item\label{cond:chartwiseRegularity}
For every chart $b$ and $x\in\cX_b$, the conditional mean
$\mu_T(x):=\E_{P_x}[T(\xi)]$ and covariance
$\Sigma_T(x):=\operatorname{Cov}_{P_x}(T(\xi))$ are defined.
There exist finite nonnegative constants $A_{\rm sc}$, $L_\Sigma$,
and $L_\mu$, fixed conditional on $\mathcal F_{x,\rm tr}$, such that,
for every chart $b$, the first two bounds below hold for all $x\in\cX_b$
and the last two hold almost surely,
\begin{gather}
\|\Sigma_T(x)-\Sigma_T(\bar x_b)\|_{\rm op}
\leq L_\Sigma\|x-\bar x_b\|_x,\qquad
\|\mu_T(x)-\mu_T(\bar x_b)\|_2
\leq L_\mu\|x-\bar x_b\|_x,
\label{eq:chartwiseMomentBounds}\\
\overline H_b\|\widehat D_b^{-1}\|_{\rm op}^2\leq A_{\rm sc},\qquad
2R_T\|\widehat D_b^{-1}\|_{\rm op}\leq\min\{A_{\rm sc},r_{y,b}\}.
\label{eq:chartwiseRegularity}
\end{gather}
\end{enumerate}
\end{assumption}

Let $\Sigma_b:=\widehat D_b^{-1}\Sigma_T(\bar x_b)\widehat D_b^{-1}$
denote the standardized conditional covariance at the chart centroid $\bar x_b$.
We define the corresponding ideal affine-gradient signal as
\begin{equation}
 C_b^\star:=\Sigma_bV_b^{\rm rel}
 H_b(z_b^{\rm SAA},\bar x_b)^\top
 \in\RR^{d_\xi\times d_z}.
 \label{eq:idealSignalDefinition}
\end{equation}
The estimation error of $\widehat C_b$ relative to $C_b^\star$ is
$\varepsilon_b^C:=\|\widehat C_b-C_b^\star\|_{\rm op}$.
The following theorem bounds this error under
Assumption~\ref{ass:finiteSampleCovariance}.

\begin{theorem}[Tradeoff between localization and sampling]
\label{thm:atlasCrossCovarianceTradeoff}
Define the empirical mean squared distance
to chart centroids and the weighted average affine residual as
$\widehat J_B:=|I_{\rm tr}|^{-1}\sum_{i\in I_{\rm tr}}
\|x_i-\bar x_{\tau_B(x_i)}\|_x^2$ and
$\bar\varepsilon_B^{\rm aff}:=\sum_{b=1}^B(|I_b|/|I_{\rm tr}|)
\varepsilon_b^{\rm aff}$.
Fix $\delta\in(0,1)$, and suppose Assumptions~\ref{ass:regularizedProblem},
\ref{ass:differentiableScenarioCosts}, and~\ref{ass:finiteSampleCovariance} hold and
$|I_b|\geq\max\{2,\log\{4(d_\xi+1)B/\delta\}\}$ for every chart.
Then, with probability at least $1-\delta$ conditional on
$\mathcal F_{x,\rm tr}$, we have
\begin{equation}
\begin{aligned}
 \sum_{b=1}^B \frac{|I_b|}{|I_{\rm tr}|}\varepsilon_b^C
 \leq
 A_{\rm sc}
 \Biggl\{
 &\underbrace{L_\Sigma\sqrt{\widehat J_B}}_{\substack{\text{within-chart}\\\text{covariance variation}}}
 +\underbrace{L_\mu^2\widehat J_B}_{\substack{\text{within-chart}\\\text{mean variation}}}
 +\underbrace{\bar\varepsilon_B^{\rm aff}}_{\substack{\text{within-chart}\\\text{affine residual}}}\\
 &+\underbrace{8R_T^2
 \sqrt{\frac{B\log\{4(d_\xi+1)B/\delta\}}{|I_{\rm tr}|}}}_{\substack{\text{fewer observations}\\\text{per chart}}}
 \Biggr\}.
\end{aligned}
 \label{eq:atlasCrossCovarianceTradeoff}
\end{equation}
\end{theorem}
\proof{Proof.} See Section~\ref{ec:proofAtlasCrossCovarianceTradeoff}
of the e-companion. \Halmos \endproof

In Eq.~\eqref{eq:atlasCrossCovarianceTradeoff}, the terms involving
$\widehat J_B$ measure within-chart variation in conditional means $\mu_T(x)$ and
covariances $\Sigma_T(x)$, while $\bar\varepsilon_B^{\rm aff}$ captures departures
from affine-gradient dependence. Increasing $B$ may reduce these errors
but leaves fewer observations per chart, increasing the sampling term
at a rate of $\sqrt{B/|I_{\rm tr}|}$, up to logarithmic factors.
The tradeoff motivates screening values of $B$ using operational loss
on internal assessment observations in Section~\ref{subsec:integratedTraining}.

\subsection{Compact scenario supports}
\label{subsec:responseCompleteCapacity}
\label{subsec:sharedScenarioSupport}

For a fixed chart $b$, we study when the local basis formed from $\widehat C_b$ can represent
the gradient-relevant projections and when $K_b=d_b+1$ shared scenarios
are sufficient and minimal. If the signal is nondegenerate and
$\widehat C_b$ is sufficiently close to $C_b^\star$, the leading
$d_b^\star$ left singular vectors of $\widehat C_b$ can represent each
relevant projection. Under the
local affine-gradient model, approximating conditional mean gradients $\bar\phi(z;x)$
then reduces, up to the residual error, to matching conditional means
in these directions. With unrestricted probabilities, the $d_b^\star+1$
vertices of a feasible enclosing simplex can match these means. When the
means have affine dimension $d_b^\star$, exact matching requires at least
as many shared scenarios.

\begin{assumption}[Signal nondegeneracy]
\label{ass:localSignalNondegeneracy}
The matrix $H_b(z_b^{\rm SAA},\bar x_b)$ has full column rank
$d_b^\star$, and $(V_b^{\rm rel})^\top\Sigma_bV_b^{\rm rel}$
is positive definite.
\end{assumption}

Under Assumption~\ref{ass:localSignalNondegeneracy},
Proposition~\ref{thm:spectralDimensionRecovery} in the e-companion
shows that the retention rule in Eq.~\eqref{eq:spectralWidthRule}
recovers the structural dimension $d_b^\star$ with high probability
when $\alpha_{\rm ret}$ is suitably chosen and the high-probability
bound on the cross-covariance estimation error $\varepsilon_b^C$ is
sufficiently small.

For $x\in\cX_b$, define the projected conditional mean
$\mu_b^{\rm rel}(x):=\E_{P_x}[(V_b^{\rm rel})^\top T_b(\xi)]$
and let $\mathcal M_b:=\{\mu_b^{\rm rel}(x):x\in\cX_b\}$ collect these means.
If $\mathcal M_b$ is bounded, it can be enclosed in a simplex in
$\RR^{d_b^\star}$ with vertices
$u_{b,1}^{\rm rel},\ldots,u_{b,d_b^\star+1}^{\rm rel}$,
as constructed in Section~\ref{ec:proofTiltedResponseCapacity} of the
e-companion. For these vertices, we define the best uniform approximation
error of the affine-softmax model in Eq.~\eqref{eq:localProbabilityModel} as
\begin{equation}
 \varepsilon_b^{\rm mean}
 :=\inf_{a_b\in\RR^{d_b^\star+1},\,\Gamma_b\in\RR^{(d_b^\star+1)\times d_x}}\;
 \sup_{x\in\cX_b}
 \left\|
 \mu_b^{\rm rel}(x)-
 \sum_{k=1}^{d_b^\star+1}
 p_{b,k}(x;a_b,\Gamma_b)u_{b,k}^{\rm rel}
 \right\|_2.
 \label{eq:conditionalMeanApproximationError}
\end{equation}

Let $\widehat V_{b,d_b^\star}$ contain the first $d_b^\star$ left singular
vectors of $\widehat C_b$, and define
$\sigma_b:=\sigma_{d_b^\star}(C_b^\star)$.
Under Assumption~\ref{ass:localSignalNondegeneracy},
Lemma~\ref{lem:idealTiltedSpace} in the e-companion establishes that
$\operatorname{rank}(C_b^\star)=d_b^\star$ and that
$(V_b^{\rm rel})^\top V_b^\star$ is invertible for any matrix
$V_b^\star\in\RR^{d_\xi\times d_b^\star}$ whose columns form an
orthonormal basis of $\operatorname{range}(C_b^\star)$.
We can therefore define
$\kappa_b:=\|((V_b^{\rm rel})^\top V_b^\star)^{-1}\|_{\rm op}$.

\begin{theorem}[Sufficiency and minimality of shared supports]
\label{thm:tiltedResponseCapacity}
\leavevmode\par
\begin{enumerate}[label=(\roman*),ref=\roman*,align=left]

\item\label{part:sharedSupportInvertibility}
Suppose Assumption~\ref{ass:localSignalNondegeneracy} holds and
$\varepsilon_b^C<\sigma_b/(1+\sqrt{2}\kappa_b)$.
Then $(V_b^{\rm rel})^\top\widehat V_{b,d_b^\star}$ is invertible.

\item\label{part:sharedSupportSufficiency}
Under the conditions of (i), suppose Assumptions~\ref{ass:differentiableScenarioCosts}
and~\ref{ass:affineUniformBounds} hold and $\mathcal M_b$ is bounded.
For the simplex vertices $u_{b,k}^{\rm rel}$, define
$u_{b,k}^\circ:=((V_b^{\rm rel})^\top\widehat V_{b,d_b^\star})^{-1}
u_{b,k}^{\rm rel}$ for $k=1,\ldots,d_b^\star+1$.
If $\|u_{b,k}^\circ\|_2\leq r_{y,b}$ for all $k$, the generated scenarios
$\tilde{\boldsymbol{\xi}}_b^\circ:=
(g_b(\widehat V_{b,d_b^\star}u_{b,k}^\circ))_{k=1}^{d_b^\star+1}$
form a fixed shared support.
For every $\varepsilon>0$, let
$\varepsilon_b^\circ(\varepsilon):=
\overline H_b(\varepsilon_b^{\rm mean}+\varepsilon)
+2\varepsilon_b^{\rm aff}$.
Then there exists an affine-softmax probability model $p_b^\circ$ on this
support whose induced distributions $P_{b,x}^\circ$ satisfy
$\sup_{x\in\cX_b,\,z\in\cZ_b}
\|\nabla_z R(z;P_x,x)-\nabla_z R(z;P_{b,x}^\circ,x)\|_2
\leq\varepsilon_b^\circ(\varepsilon)$.

\item\label{part:sharedSupportMinimality}
Let $u_k=(V_b^{\rm rel})^\top y_k$ ($k=1,\ldots,K$) be the gradient-relevant
projections of $K$ shared scenarios $\tilde{\xi}_k=g_b(y_k)$ in chart $b$.
If, for every $x\in\cX_b$, there exist probabilities
$p(x)\in\Delta_K$ such that
$\sum_{k=1}^Kp_k(x)u_k=\mu_b^{\rm rel}(x)$, then
$K\geq\dim\operatorname{aff}(\mathcal M_b)+1$,
where $\operatorname{aff}(\mathcal M_b)$ denotes the affine hull
of $\mathcal M_b$. In particular, $K\geq d_b^\star+1$
when $\mathcal M_b$ has nonempty interior in $\RR^{d_b^\star}$.

\end{enumerate}
\end{theorem}

\proof{Proof.} See Section~\ref{ec:proofTiltedResponseCapacity} of the
e-companion. \Halmos \endproof

Parts~(\ref{part:sharedSupportSufficiency}) and~(\ref{part:sharedSupportMinimality})
of Theorem~\ref{thm:tiltedResponseCapacity} motivate setting $K_b=d_b+1$
when Eq.~\eqref{eq:spectralWidthRule} recovers $d_b=d_b^\star$.
In part~(\ref{part:sharedSupportSufficiency}), the simplex determines
fixed shared scenarios $\tilde{\boldsymbol{\xi}}_b^\circ$. During model fitting, we fit the scenario coordinates $u_{b,k}$
jointly with the context-dependent probability parameters $(a_b,\Gamma_b)$.

\subsection{Quality of fitted models}
\label{subsec:operatorGuarantee}

Sections~\ref{subsec:chartTradeoff} and~\ref{subsec:sharedScenarioSupport}
motivate the chart partition and the number of scenarios in each chart.
We now study the out-of-sample operational risk of fixed fitted models.
We first compare the fitted objective gradient with the population
gradient $\nabla_z R(z;P_{\bar x_b},\bar x_b)$ at each chart centroid $\bar x_b$. If their difference is
Lipschitz in the context, the centroid error and the distance to the
centroid bound the gradient error throughout the chart. Local strong convexity then
converts this error into bounds on the decision error and regularized
objective gap (Lemma~\ref{thm:localDecisionGuarantee} in the e-companion).
Finally, accounting for $\eta_{\rm reg}$ then bounds excess
operational risk at $x$.

For chart $b$, the fitted scenarios $\tilde{\xi}_{b,1},\ldots,\tilde{\xi}_{b,K_b}$
and probability map $\widehat p_b:\cX_b\to\Delta_{K_b}$
define the distribution $\widehat P_{b,x}$ at context $x$.
Eq.~\eqref{eq:finiteDecisionLayer} gives the fitted objective
$\widehat R_b(z;x):=R(z;\widehat P_{b,x},x)$
and its induced decision
$\widehat z_b(x):=z(\widehat P_{b,x};x)$.
The gradient approximation error at the centroid $\bar x_b$ is
\begin{equation}
\varepsilon_b^{\rm ctr}
 :=
 \sup_{z\in\cZ_b}
 \bigl\|\nabla_z R(z;P_{\bar x_b},\bar x_b)-\nabla_z\widehat R_b(z;\bar x_b)\bigr\|_2.
 \label{eq:centroidOperatorError}
\end{equation}
To extend the centroid error bound $\varepsilon_b^{\rm ctr}$ to other contexts in $\cX_b$,
we impose the following Lipschitz continuity condition on
the gradient error.

\begin{assumption}[Context regularity]
\label{ass:operatorContextRegularity}
There exists a finite constant $L_B\geq0$ such that, for every chart $b$,
$z\in\cZ_b$, and $x\in\cX_b$,
\begin{equation}
 \bigl\|\{\nabla_z R(z;P_x,x)-\nabla_z\widehat R_b(z;x)\}
 -\{\nabla_z R(z;P_{\bar x_b},\bar x_b)-\nabla_z\widehat R_b(z;\bar x_b)\}\bigr\|_2
 \leq L_B\|x-\bar x_b\|_x.
 \label{eq:contextRegularity}
\end{equation}
\end{assumption}

We next consider an independent deployment context $x$, conditional
on the fitted model.
To state the out-of-sample risk bound, define
$J_B:=\E[\|x-\bar x_{\tau_B(x)}\|_x^2]$,
$\varepsilon_B^{\rm ctr}:=\{\E[(\varepsilon_{\tau_B(x)}^{\rm ctr})^2]\}^{1/2}$,
and
$R_0^\star(x):=\inf_{z\in\cZ}R_0(z;x)$.

\begin{theorem}[Out-of-sample operational risk]
\label{thm:fixedAtlasRisk}
Suppose Assumptions~\ref{ass:regularizedProblem},
\ref{ass:differentiableScenarioCosts}, and~\ref{ass:operatorContextRegularity}
hold, with $\varepsilon_{\tau_B(x)}^{\rm ctr}+L_B\|x-\bar x_{\tau_B(x)}\|_x<\eta r_{\rm loc}$ almost surely,
$J_B<\infty$, and $\varepsilon_B^{\rm ctr}<\infty$.
Assume there exists a finite radius $\operatorname{rad}_{\cZ,0}\geq0$
such that, for every $x\in\cX$,
$\inf_{z\in\cZ:\|z-z_0(x)\|_2\leq\operatorname{rad}_{\cZ,0}}R_0(z;x)=R_0^\star(x)$.
If $R_0(\widehat z_{\tau_B(x)}(x);x)$ and $R_0^\star(x)$ are integrable,
then the deployed decision satisfies
\begin{equation}
 \E[R_0(\widehat z_{\tau_B(x)}(x);x)]-\E[R_0^\star(x)]
 \leq
 \frac{\bigl(\varepsilon_B^{\rm ctr}+L_B\sqrt{J_B}\bigr)^2}{2\eta}
 +
 \frac{\eta_{\rm reg} \operatorname{rad}_{\cZ,0}^2}{2}.
 \label{eq:fixedAtlasRisk}
\end{equation}
\end{theorem}

\proof{Proof.} See Section~\ref{ec:unregularizedRisk} of the e-companion.
\Halmos \endproof

In Eq.~\eqref{eq:fixedAtlasRisk}, the first term accounts for centroid
gradient error $\varepsilon_B^{\rm ctr}$ and within-chart variation, while the second bounds the
bias introduced by quadratic regularization \citep[as in][]{qi2026integrated}.
If $R_0(\cdot;x)$ already satisfies Eq.~\eqref{eq:strongConvexity},
setting $\eta_{\rm reg}=0$ eliminates this bias.

\subsection{Finite-sample generalization}
\label{subsec:finiteSampleGeneralization}

We analyze the finite-sample generalization error of model fitting in
Eq.~\eqref{eq:integratedTrainingObjective}, with the chart structure fixed.
Assumption~\ref{ass:stageTwoRegularity} requires the induced operational
loss to be Lipschitz in the model parameters $\theta_b$.
Lemma~\ref{lem:ecChartParameterLipschitz} in the e-companion derives this
property from strong convexity of the local problems, bounded chart radii,
bounded and Lipschitz scenario gradients $\phi$, and Lipschitz operational
loss $\ell$. We first use concentration inequalities to bound the gap between
empirical and population losses at a finite set of points covering the
bounded parameter space. A union bound controls this gap simultaneously
at these points, and Lipschitz continuity extends the bound to the entire
parameter space. Finally, we use the shared-support model constructed in
Theorem~\ref{thm:tiltedResponseCapacity}(\ref{part:sharedSupportSufficiency})
as a comparator to derive an end-to-end excess-risk bound for the fitted
model.

For chart $b$, let
$\theta_b=(u_{b,1},\ldots,u_{b,K_b},a_b,\Gamma_b)$
collect its scenario coordinates and probability parameters, with norm
$\|\theta_b\|_2:=(\sum_{k=1}^{K_b}\|u_{b,k}\|_2^2+\|a_b\|_2^2+\|\Gamma_b\|_F^2)^{1/2}$.
For the full parameter vector $\theta=(\theta_1,\ldots,\theta_B)$,
the decision rule is
$\pi_{\rm Atlas}(x;\theta):=z_{\tau_B(x)}(x;\theta_{\tau_B(x)})$,
which uses the model associated with the assigned chart.

Let $\{(x_i,\xi_i)\}_{i=1}^n$ be $n\geq1$ independent observations drawn
from the joint distribution of $(x,\xi)$.
The chart structure
$(\tau_B,\widehat D_x,\{\bar x_b,\widehat V_b,g_b,d_b,K_b\}_{b=1}^B)$
and nonempty parameter sets $\Theta_b$ are independent of these observations,
as when Steps~1 and~2 of Algorithm~\ref{alg:kktAtlas} use disjoint subsets of
the training data.
We define the full parameter space as
$\Theta_{\rm Atlas}:=\Theta_1\times\cdots\times\Theta_B$.
Throughout this subsection, all probabilities and expectations are
conditional on the chart structure and parameter sets.
We impose the following regularity conditions on the parameter space
and the induced losses.

\begin{assumption}[Parameter regularity]
\label{ass:stageTwoRegularity}
There exist finite positive constants $R_\Theta$, $L_\Theta$, and $M_\ell$
such that the following conditions hold.
\begin{enumerate}[label=(\roman*),ref=\roman*]
\item\label{cond:stageTwoDomain}
For every chart $b$, each $\theta_b\in\Theta_b$ satisfies
$\mathbf 1^\top a_b=0$, $\mathbf 1^\top\Gamma_b=0$,
$\|\theta_b\|_2\leq R_\Theta$, and
$\widehat V_bu_{b,k}\in\operatorname{dom}(g_b)$ for $k=1,\ldots,K_b$.

\item\label{cond:stageTwoLipschitz}
For every chart $b$, $x\in\cX_b$, $\xi\in\Xi$, and
$\theta_b,\theta_b'\in\Theta_b$, we have
$\bigl|\ell(z_b(x;\theta_b),\xi;x)-\ell(z_b(x;\theta_b'),\xi;x)\bigr|
\leq L_\Theta\|\theta_b-\theta_b'\|_2$.

\item\label{cond:stageTwoLoss}
For every $\theta\in\Theta_{\rm Atlas}$, the induced loss
$(x,\xi)\mapsto\ell(\pi_{\rm Atlas}(x;\theta),\xi;x)$ is measurable
and takes values in an interval of length $M_\ell$ over
$\cX\times\Xi$.
\end{enumerate}
\end{assumption}

In chart $b$, the $K_b$ scenario
coordinate vectors have $d_b$ entries each, giving $K_bd_b$ parameters, and the
zero-sum constraints in Eq.~\eqref{eq:screeningProbabilityFit} leave
$(K_b-1)(d_x+1)$ free probability parameters.
The total number of free parameters is
$P_{\rm Atlas}:=\sum_{b=1}^B\{K_bd_b+(K_b-1)(d_x+1)\}$.

The empirical operational loss is
$\widehat{\mathcal L}_n(\theta):=n^{-1}\sum_{i=1}^n
\ell(\pi_{\rm Atlas}(x_i;\theta),\xi_i;x_i)$,
and the corresponding population loss is
$\mathcal L(\theta):=\E[\ell(\pi_{\rm Atlas}(x;\theta),\xi;x)]$.
For each chart $b$, let
$I_{b,n}:=\{1\leq i\leq n:\tau_B(x_i)=b\}$
index its training observations. When $I_{b,n}$ is nonempty, we fit
$\theta_b\in\Theta_b$ by minimizing
$|I_{b,n}|^{-1}\sum_{i\in I_{b,n}}
\ell(z_b(x_i;\theta_b),\xi_i;x_i)+\Omega_b(\theta_b)$,
as in Eq.~\eqref{eq:integratedTrainingObjective}.

\begin{theorem}[Generalization bound]
\label{thm:stageTwoGeneralization}
Fix $\delta\in(0,1)$, and suppose Assumptions~\ref{ass:regularizedProblem}
and~\ref{ass:stageTwoRegularity} hold. Then, with probability at least
$1-\delta$, the following bounds hold.
\begin{enumerate}[label=(\roman*),ref=\roman*,align=left]
\item\label{part:uniformGeneralization}
We have
\begin{equation}
\begin{aligned}
 \sup_{\theta\in\Theta_{\rm Atlas}}
 \bigl|\widehat{\mathcal L}_n(\theta)-\mathcal L(\theta)\bigr|
 &\leq
 \varepsilon_n^{\rm est}(\delta)\\
 &:=\frac{M_\ell}{\sqrt{n}}
 \left[
 1+\sqrt{\frac{P_{\rm Atlas}\log(1+4R_\Theta L_\Theta\sqrt{n}/M_\ell)+\log(2/\delta)}{2}}
 \right].
\end{aligned}
 \label{eq:stageTwoUniformDeviation}
\end{equation}

\item\label{part:oracleInequality}
If $\widehat\theta=(\widehat\theta_1,\ldots,\widehat\theta_B)\in\Theta_{\rm Atlas}$
minimizes each chart objective over $\Theta_b$ up to additive error
$\varepsilon_{\rm opt}\geq0$ whenever $I_{b,n}$ is nonempty, then
\begin{equation}
 \mathcal L(\widehat\theta)
 \leq
 \inf_{\theta\in\Theta_{\rm Atlas}}
 \Bigl\{\mathcal L(\theta)+\max_{1\leq b\leq B}\Omega_b(\theta_b)\Bigr\}
 +2\varepsilon_n^{\rm est}(\delta)+\varepsilon_{\rm opt}.
 \label{eq:stageTwoOracleInequality}
\end{equation}
\end{enumerate}
\end{theorem}
\proof{Proof.} See Section~\ref{ec:proofStageTwoGeneralization} of the
e-companion. \Halmos \endproof

In Eq.~\eqref{eq:stageTwoUniformDeviation}, $P_{\rm Atlas}$ measures the
dimension of the model being fitted. The size of any reference support
does not enter the bound explicitly. The proof of
Corollary~\ref{cor:endToEndRisk} in Section~\ref{ec:proofEndToEndRisk} of the
e-companion uses the shared-support model from
Theorem~\ref{thm:tiltedResponseCapacity}(\ref{part:sharedSupportSufficiency})
as the comparator in Eq.~\eqref{eq:stageTwoOracleInequality}.

\begin{corollary}[Excess risk of the fitted model]
\label{cor:endToEndRisk}
Suppose the conditions of
Theorem~\ref{thm:stageTwoGeneralization}(\ref{part:oracleInequality})
hold with $d_b=d_b^\star$ and $K_b=d_b^\star+1$ for every chart.
For some $\varepsilon>0$, assume the conditions of
Theorem~\ref{thm:tiltedResponseCapacity}(\ref{part:sharedSupportSufficiency})
also hold in every chart; let
$\theta_b^\circ:=(u_{b,1}^\circ,\ldots,u_{b,K_b}^\circ,a_b^\circ,\Gamma_b^\circ)$
collect the resulting scenario coordinates and zero-sum-centered parameters of
$p_b^\circ$, and assume $\theta_b^\circ\in\Theta_b$ and
$\varepsilon_b^\circ(\varepsilon)<\eta r_{\rm loc}$ for every chart $b$.
If $R_0^\star(x)$ is integrable and a finite
$\operatorname{rad}_{\cZ,0}$ satisfies the condition in
Theorem~\ref{thm:fixedAtlasRisk}, then, with probability at least
$1-\delta$,
\begin{equation}
 \mathcal L(\widehat\theta)-\E[R_0^\star(x)]
 \leq
 \max_{1\leq b\leq B}\frac{\{\varepsilon_b^\circ(\varepsilon)\}^2}{2\eta}
 +\frac{\eta_{\rm reg}\operatorname{rad}_{\cZ,0}^2}{2}
 +\max_{1\leq b\leq B}\Omega_b(\theta_b^\circ)
 +2\varepsilon_n^{\rm est}(\delta)+\varepsilon_{\rm opt}.
 \label{eq:endToEndRisk}
\end{equation}
\end{corollary}

Eq.~\eqref{eq:endToEndRisk} shows the benefit of compactness when the
shared support approximates the population gradient accurately.
Holding other constants fixed, a smaller $P_{\rm Atlas}$ reduces the
estimation term $\varepsilon_n^{\rm est}(\delta)$ for a given sample size. When shared-support approximation
dominates the bound, more observations leave the approximation
term unchanged for a fixed chart structure and comparator.
In this analysis, ensuring that the chart structure is independent of the
fitting sample requires sample splitting, using separate samples for structure
selection (Step~1) and model fitting (Step~2). In practice,
Algorithm~\ref{alg:kktAtlas} uses the full training sample for both steps to
preserve sample efficiency, and Section~\ref{subsubsec:round7Robustness} shows
empirically that reusing the training sample does not degrade out-of-sample
performance relative to sample splitting.

\section{Computational Study}
\label{sec:computation}
\setlength{\textfloatsep}{12pt plus 3pt minus 2pt}
\setlength{\intextsep}{12pt plus 3pt minus 2pt}
\renewcommand{\topfraction}{0.88}
\renewcommand{\bottomfraction}{0.65}
\renewcommand{\textfraction}{0.10}
\renewcommand{\floatpagefraction}{0.75}

We introduce the decision problems and experimental design in
Sections~\ref{subsec:decisionProblems} and~\ref{subsec:experimentalDesign},
and report the computational results for three decision problems using
synthetic data in Section~\ref{subsec:applicationResults}.
Section~\ref{subsec:grid64Empirical} then evaluates KKT Atlas on
a stylized day-ahead capacity-coverage benchmark built from real European
power-grid data \citep{wiese2019open,hirth2018entsoe}, and Section~\ref{subsec:ablationStudies}
groups the controlled experiments, ablations, and sensitivity analyses
by design choice.

\subsection{Decision problems}
\label{subsec:decisionProblems}

We consider three applications in capacity allocation, CVaR
portfolio optimization, and network reservation.
We first introduce the common notation and evaluation conventions.

In each application, a fitted method maps an observed context to a finite
scenario distribution $\widetilde P_x$ specified by $\{(p_k,\tilde{\xi}_k)\}_{k=1}^K$ and computes a feasible
decision by minimizing the regularized downstream objective
$R(z;\widetilde P_x,x)$ in Eq.~\eqref{eq:finiteDecisionLayer} with the
regularization center $z_0(x)$ and coefficient
$\eta_{\rm reg}\geq0$. Here $\eta_{\rm reg}$ supplies curvature when missing
from the operational loss $\ell$, whereas out-of-sample operational risk is
evaluated using $\ell$ alone.
Because the conditional distribution is known in each synthetic experiment, we compute
the decision that solves the unregularized contextual stochastic program in
Eq.~\eqref{eq:populationObjectives} under that
distribution and use its risk as the reference for evaluating the true raw
and relative excess risks.
All methods use the same evaluator for operational loss and the same feasible
decision set, with different finite distributions defining their stochastic
programs.

\paragraph{Smooth capacity allocation (CAP).}
The CAP problem allocates total capacity $C>0$ across $d_z$ resources that
jointly serve $d_\xi$ demand classes. We adapt the model of demand uncertainty
across multiple resources from \citet{harrison1999multiresource} to smooth
capacity allocation. Let $A\in\RR_+^{d_\xi\times d_z}$ have row sums equal
to one and define
\(\cZ_{\rm CAP}(C):=\{z\in\RR_+^{d_z}:\mathbf 1^\top z=C\}\). With
$h_\tau(a):=\tau\log\{1+\exp(a/\tau)\}$ ($\tau>0$) and overage and underage cost
coefficients $c_j^o,c_j^u>0$, the operational loss under demand realization
$\xi\in(0,\infty)^{d_\xi}$ is
\begin{equation}
 \ell_{\rm CAP}(z,\xi)
 :=
 \frac1{d_\xi}\sum_{j=1}^{d_\xi}
 \left[
 c_j^o h_\tau((Az)_j-\xi_j)
 +c_j^u h_\tau(\xi_j-(Az)_j)
 \right].
 \label{eq:capacityComputationalProblem}
\end{equation}
In CAP, $d_z<d_\xi$, so the uncertain demand vector has more components than
the capacity allocation decision. CAP therefore tests whether compact local
scenario models remain effective when the uncertainty dimension exceeds the
decision dimension. Section~\ref{ec:capacityProblem} of the e-companion gives the parameter
values, the contextual demand model with four states, and the strong
convexity proof.

\paragraph{Smoothed CVaR portfolio optimization (PORT).}
We use a smoothed CVaR portfolio formulation following
\citet{rockafellar2000optimization}. The PORT problem selects portfolio
weights $w\in\RR^{d_\xi}$ and a threshold $t\in\RR$, subject to
$\mathbf 1^\top w=1$. For a return realization
$\xi\in(-r_{\max},r_{\max})^{d_\xi}$, tail probability $\alpha\in(0,1)$,
expected-return coefficient $\lambda_\mu\geq0$, and leverage coefficient
$\lambda_w\geq0$, the operational loss is
\begin{equation}
 \ell_{\rm PORT}((w,t),\xi)
 :=
 -\lambda_\mu\xi^\top w+t
 +\frac{1}{\alpha}h_\tau(-\xi^\top w-t)
 +\frac{\lambda_w}{2}\bigl\|w\bigr\|_2^2.
 \label{eq:portfolioComputationalProblem}
\end{equation}
Although $(w,t)\in\RR^{d_\xi+1}$, the equality constraint
$\mathbf 1^\top w=1$ removes one degree of freedom. Thus, the feasible
decision set has dimension $d_\xi$, matching that of the uncertain
parameter. PORT therefore tests KKT Atlas when the decision and uncertainty
dimensions match. We specify the bounded factor model, including a
market-wide crash state that jointly reduces asset returns, and the numerical
parameters in Section~\ref{ec:portfolioProblem} of the e-companion.

\paragraph{Multicommodity network reservation (NET).}
The NET problem models advance capacity reservation under uncertain demand
between origins and destinations \citep{riis2002capacitated}. The decision
\(z\in\cZ_{\rm NET}=[0,1.5]^{d_z}\), made before demand is observed,
reserves link capacity. For a demand realization
\(\xi\in(0,\infty)^{d_\xi}\), each of the \(d_\xi\) commodities can use two
fixed paths. Let \(P^{(1)},P^{(2)}\in\{0,1\}^{d_z\times d_\xi}\) denote the
link incidence matrices for these paths, and let $q\in\RR^{d_\xi}$ denote
the routing vector, where $q_j\in[0,\xi_j]$ is the amount of commodity $j$
routed on its first path. The resulting link load is
\(L(q,\xi):=P^{(2)}\xi+\{P^{(1)}-P^{(2)}\}q\). For reservation costs
\(c\in\RR_+^{d_z}\), congestion coefficient \(\lambda>0\), and routing
coefficient \(\rho>0\), the operational loss is
\begin{equation}
 \ell_{\rm NET}(z,\xi)
 :=
 c^\top z+
 \min_{0\leq q\leq\xi}
 \left\{
 \lambda\bigl\|\max\{L(q,\xi)-z,\mathbf 0\}\bigr\|_2^2
 +\rho\bigl\|q-\xi/2\bigr\|_2^2
 \right\},
 \label{eq:reservationComputationalProblem}
\end{equation}
where $\mathbf 0\in\RR^{d_z}$ denotes the zero vector and the maximum is
taken componentwise.
In NET, $d_z$ remains fixed while $d_\xi$ determines the number of
commodities and routing recourse variables per scenario. NET therefore
tests whether compact local scenario supports retain decision quality while
reducing the routing recourse work required in downstream solves. We specify
the network, demand model, path construction, and numerical parameters in
Section~\ref{ec:reservationProblem} of the e-companion.

We compare KKT Atlas with six benchmarks, summarized in
Table~\ref{tab:benchmarkMethods}. Let $K_{\rm ref}$ denote the size of
the full reference support
$\{\tilde{\xi}^{(s)}:s=1,\ldots,K_{\rm ref}\}\subseteq\Xi$. The benchmarks are ICEO,
stochastic optimization forests (SOF), the $k$-nearest neighbor method (kNN),
their scenario-reduction variants ICEO-red and kNN-red, and unconditional SAA. Since neither
scenario-reduction benchmark
uses context charts, we set the support budget for both to
$K_{\rm match}:=\max_b K_b$.

\begin{table}[htbp]
\centering
\caption{Summary of the six benchmarks.}
\label{tab:benchmarkMethods}
\small
\renewcommand{\arraystretch}{1.04}
\begin{tabularx}{\textwidth}{@{}lX@{}}
\toprule
Method & Main idea \\
\midrule
ICEO
& Learns context-dependent probabilities over all $K_{\rm ref}$ reference scenarios by minimizing operational loss \citep{qi2026integrated}. \\\addlinespace[2pt]
ICEO-red
& Applies K-means clustering in standardized $T$-coordinates and decodes the centroids via $g$ to construct $K_{\rm match}$ fixed scenarios from training realizations, then learns context-dependent probabilities by minimizing operational loss \citep{qi2026integrated,rujeerapaiboon2022scenario}. \\\addlinespace[2pt]
SOF
& Uses optimization-aware tree splits to assign context-dependent weights to historical observations \citep{kallus2023stochastic}. \\\addlinespace[2pt]
kNN
& Selects the $k$ training observations with contexts closest to the new context and forms a local empirical SAA using their realizations, each with probability $1/k$ \citep{hastie2009elements}. \\\addlinespace[2pt]
kNN-red
& Independently selects $k$ by validation on compressed supports, applying weighted K-means in standardized $T$-coordinates when $k>K_{\rm match}$ to reduce the $k$ neighbor realizations to $K_{\rm match}$ decoded centroids with summed cluster weights \citep{rujeerapaiboon2022scenario}. \\\addlinespace[2pt]
SAA
& Assigns uniform weights to all training observations, ignoring context. \\
\bottomrule
\end{tabularx}
\end{table}

\subsection{Experimental design and comparisons}
\label{subsec:experimentalDesign}

We first compare Atlas with the six benchmarks on CAP, PORT, and NET using
10 paired replications for CAP and PORT and three for NET.
We use fewer replications for NET because its two-stage routing recourse is computationally expensive (Table~\ref{tab:capPortConstructionCost}), and the resulting standard errors (SE) for Atlas and ICEO are already at most 0.02 percentage points. We further examine the effects of scenario support size $K$, local basis choice,
retention threshold $\alpha_{\rm ret}$, and training sample size $|I_{\rm tr}|$. Two controlled quadratic
experiments isolate the contributions of chart localization and the
cross-covariance basis $\widehat V_b$.

Table~\ref{tab:confirmationDesigns} summarizes the application settings and
sample sizes, with $I_{\rm tr}$, $I_{\rm val}$, and $I_{\rm eval}$ indexing
training observations, validation observations, and evaluation contexts,
respectively.
Sections~\ref{ec:testProblems} and~\ref{ec:samplingEvaluation} of the
e-companion give the data generation and sampling protocols.

We evaluate decision quality on independent contexts by computing exact
conditional operational risk $R_0(z;x)$, summing losses over all $K_{\rm ref}$ reference
scenarios with their known probabilities $p_s^\star(x)$ at context $x$, for $s=1,\ldots,K_{\rm ref}$.
Relative excess risk (abbreviated Excess in tables) is the percentage increase in operational risk over the
optimum under the known conditional distribution. We also report the selected
Atlas structure, construction time, and deployment time.
Construction time covers structure or hyperparameter
selection and model fitting. Deployment time includes distribution prediction
and computation of the downstream decision. Section~\ref{ec:timingProtocol} of the e-companion defines
these quality measures and the reported construction and
deployment times.

Atlas and ICEO are implemented in PyTorch 2.12.1 \citep{paszke2019pytorch}.
For the main synthetic comparisons, we train both with Adam \citep{kingma2015adam} for
1,600 updates, using a learning rate of \(0.001\),
a minibatch size of 64, and gradient-norm clipping at 10. We select the final
model only from the saved checkpoints \(\{0,400,800,1200,1600\}\) by minimizing
operational loss on the validation sample; checkpoint 0 denotes the initialized
model. CAP uses sequential least squares programming in SciPy 1.18.0 \citep{virtanen2020scipy}, PORT uses
Clarabel 0.11.1 \citep{goulart2026clarabel}, and NET uses Gurobi 13.0.2
\citep{gurobi}.

The experiments run on a Windows server with AMD EPYC 9654 processors
(2.4~GHz, 192 physical cores) and approximately 1.5~TiB of memory.
For each replication, all methods use the same allocation of 32 physical cores
for construction, with up to 32 single-threaded workers.
Atlas and ICEO parallelize decision solves and implicit differentiation
across contexts before each sequential Adam update, with Atlas operating
separately within each chart. Neither uses the solution-mapping surrogate of
\citet{qi2026integrated}, so both are evaluated as direct-solve
implementations without surrogate approximation. kNN evaluates neighborhood
sizes sequentially but parallelizes validation solves. SOF distributes its
100 trees across workers.

Tables~\ref{tab:capPortConstructionCost} and~\ref{tab:deploymentTiming}
report timing results for six methods, excluding SAA, which requires no
model fitting. Deployment is measured on
1 and 32 physical cores; with 32 cores, each method distributes
128 contexts across 32 workers, each processing four contexts.
Section~\ref{ec:protocol} of the e-companion provides the remaining
training settings. The source code and scripts for reproducing the
results will be released as open source upon publication.

\begin{table}[htbp]
\centering
\caption{Core application dimensions and settings.}
\label{tab:confirmationDesigns}
\small
\renewcommand{\arraystretch}{1.03}
\setlength{\tabcolsep}{6pt}
\begin{tabular*}{\textwidth}{@{\extracolsep{\fill}}lrrrrrrr@{}}
\toprule
Study & $K_{\rm ref}$ & $d_\xi$ & $d_z$ &
$\eta_{\rm reg}$ & $|I_{\rm tr}|$ & $|I_{\rm val}|$ & $|I_{\rm eval}|$ \\
\midrule
CAP  & 500 & 64 & 8  & $0$       & 12,800 & 1,000 & 200 \\
PORT & 500 & 32 & 33 & $10^{-6}$ & 12,800 & 1,000 & 200 \\
NET  & 500 & 16 & 48 & $10^{-6}$ & 12,800 & 1,000 & 200 \\
\bottomrule
\end{tabular*}
\end{table}

\subsection{Results across applications}
\label{subsec:applicationResults}

\subsubsection{Decision quality.}
Table~\ref{tab:capPortAnchorQuality} reports out-of-sample decision quality
across the three applications. Atlas achieves the lowest mean operational risk
in CAP and PORT, attaining mean relative excess risks of $0.07\%$ and $2.94\%$,
respectively. Across all 23 paired replications, Atlas outperforms SOF, kNN,
and SAA in every replication and attains lower operational risk than 500-scenario
ICEO in 20 (all 10 in CAP, 9 of 10 in PORT, and 1 of 3 in NET). In NET, ICEO achieves a slightly lower mean relative excess risk ($0.09\%$ versus
$0.12\%$), whereas single-chart Atlas attains $0.07\%$
(Table~\ref{tab:ecAtlasDiagnostics} of the e-companion). In NET, Atlas fits $371$--$800$ parameters across charts,
compared with $17{,}844$ in ICEO. Neither
method achieves lower risk in all three applications. Atlas, however,
requires substantially less construction and deployment time than ICEO
(Tables~\ref{tab:capPortConstructionCost} and~\ref{tab:deploymentTiming}).

Atlas has lower mean operational risk than ICEO-red and kNN-red in all three
applications. ICEO-red trails ICEO throughout, whereas kNN-red slightly improves
on kNN in CAP but worsens in PORT and NET---most sharply in PORT, where mean
excess risk rises from $5.20\%$ to $31.15\%$ after scenario reduction, as cluster averaging pulls scenarios toward
the center and understates tail losses that drive CVaR.

\begin{table}[htbp]
\centering
\caption{Out-of-sample decision quality under known conditional distributions.}
\label{tab:capPortAnchorQuality}
\footnotesize
\renewcommand{\arraystretch}{1.00}
\setlength{\tabcolsep}{0pt}
\begin{tabular*}{\textwidth}{@{\extracolsep{\fill}}l*{6}{r}@{}}
\toprule
 & \multicolumn{2}{c}{CAP} & \multicolumn{2}{c}{PORT} & \multicolumn{2}{c}{NET} \\
\cmidrule(lr){2-3}\cmidrule(lr){4-5}\cmidrule(lr){6-7}
Method & Risk (SE) $\downarrow$ & Excess\,/\% (SE) $\downarrow$ & Risk (SE) $\downarrow$ & Excess\,/\% (SE) $\downarrow$ & Risk (SE) $\downarrow$ & Excess\,/\% (SE) $\downarrow$ \\
\midrule
Atlas & \textbf{1.5957} (0.021) & \textbf{0.07} (0.01) & \textbf{0.0380} (0.0005) & \textbf{2.94} (0.17) & 1.7253 (0.100) & 0.12 (0.01) \\
ICEO & 1.6009 (0.022) & 0.40 (0.03) & 0.0382 (0.0005) & 3.66 (0.27) & \textbf{1.7247} (0.100) & \textbf{0.09} (0.02) \\
ICEO-red & 1.6015 (0.021) & 0.44 (0.03) & 0.0393 (0.0005) & 6.69 (2.08) & 1.7259 (0.100) & 0.16 (0.03) \\
SOF & 1.6331 (0.023) & 2.41 (0.13) & 0.0408 (0.0005) & 10.66 (0.40) & 1.7384 (0.102) & 0.87 (0.11) \\
kNN & 1.6275 (0.023) & 2.06 (0.11) & 0.0388 (0.0005) & 5.20 (0.21) & 1.7343 (0.100) & 0.65 (0.06) \\
kNN-red & 1.6263 (0.023) & 1.99 (0.11) & 0.0482 (0.0033) & 31.15 (10.15) & 1.7391 (0.102) & 0.92 (0.07) \\
SAA & 1.9495 (0.031) & 22.23 (0.52) & 0.0431 (0.0005) & 16.91 (0.64) & 1.8022 (0.109) & 4.56 (0.31) \\
\midrule
Reference & 1.5945 (0.021) & 0.00 (0.00) & 0.0369 (0.0005) & 0.00 (0.00) & 1.7232 (0.100) & 0.00 (0.00) \\
\bottomrule
\end{tabular*}
\par\vspace{2pt}\begin{minipage}{\textwidth}\renewcommand{\baselinestretch}{1.02}\footnotesize
\textit{Notes.} Entries are means. Parentheses report the standard error of the mean (SE). Reference reports the mean optimal risk under the known conditional distribution.
\end{minipage}
\end{table}

\subsubsection{Construction and deployment costs.}
Table~\ref{tab:capPortConstructionCost} reports mean construction wall times
over replications 1--3, of 0.03 versus 0.11 hours for Atlas and ICEO in CAP,
0.02 versus 0.08 hours in PORT, and 0.50 versus 3.86 hours in NET,
corresponding to $3$--$8\times$ speedups based on the unrounded times. Both
methods repeatedly solve stochastic programs during fitting, but Atlas solves
smaller programs on its compact supports. ICEO-red similarly cuts training cost
but incurs higher risk.

Compared with SOF and kNN, Atlas selects charts, estimates local bases, and
jointly fits scenarios and probabilities, whereas SOF splits via approximate-risk
scores and kNN avoids iterative fitting.

\begin{table}[htbp]
\centering
\caption{Model construction time.}
\label{tab:capPortConstructionCost}
\small
\renewcommand{\arraystretch}{1.00}
\setlength{\tabcolsep}{4pt}
\begin{tabular*}{\textwidth}{@{\extracolsep{\fill}}l*{6}{r}@{}}
\toprule
 & \multicolumn{2}{c}{CAP} & \multicolumn{2}{c}{PORT} & \multicolumn{2}{c}{NET} \\
\cmidrule(lr){2-3} \cmidrule(lr){4-5} \cmidrule(lr){6-7}
Method & Wall\,/h $\downarrow$ & CPU\,/h $\downarrow$ & Wall\,/h $\downarrow$ & CPU\,/h $\downarrow$ & Wall\,/h $\downarrow$ & CPU\,/h $\downarrow$ \\
\midrule
Atlas & 0.03 & 0.18 & 0.02 & 0.08 & 0.50 & 2.49 \\
ICEO & 0.11 & 2.30 & 0.08 & 1.18 & 3.86 & 69.11 \\
ICEO-red & 0.02 & $<0.01^\ast$ & 0.02 & $<0.01^\ast$ & 0.13 & 1.19 \\
SOF & 0.03 & 0.69 & 0.02 & 0.37 & 0.66 & 15.80 \\
kNN & 0.01 & 0.16 & $<0.01$ & 0.08 & 0.12 & 3.69 \\
kNN-red & 0.02 & 0.04 & 0.03 & 0.03 & 0.03 & 0.16 \\
\bottomrule
\end{tabular*}
\par\vspace{2pt}\begin{minipage}{\textwidth}\renewcommand{\baselinestretch}{1.02}\footnotesize
\textit{Notes.} Entries are mean construction times over replications 1--3 on 32 physical cores. $^\ast$On fast CAP and PORT subproblems, individual solves are shorter than the $1/64$~s operating-system CPU accounting tick, so measured CPU time can be smaller than wall time.
\end{minipage}
\end{table}

Table~\ref{tab:deploymentTiming} reports deployment wall times, measured from
predicting a scenario distribution to returning the decision.
For batches of 128 contexts, Atlas achieves approximately
$15\times$, $19\times$, and $49\times$ speedups over ICEO on one core in
CAP, PORT, and NET, respectively, and $8$--$48\times$ speedups on 32 cores.
The advantage is largest in NET, where fewer scenarios also reduce the size
of the recourse problem. Atlas also deploys faster than SOF and kNN in every
setting.

Among compact-support methods, Atlas and ICEO-red deploy at similar speeds
because neither clusters online, with Atlas achieving lower risk throughout
(Table~\ref{tab:capPortAnchorQuality}). In contrast, kNN-red clusters neighbor
realizations online for each new context, making it slower than Atlas across
all settings.

\begin{table}[htbp]
\centering
\caption{Model deployment time.}
\label{tab:deploymentTiming}
\small
\renewcommand{\arraystretch}{1.00}
\setlength{\tabcolsep}{3pt}
\begin{tabular*}{\textwidth}{@{\extracolsep{\fill}}l*{9}{r}@{}}
\toprule
& \multicolumn{3}{c}{1 context, 1 core} & \multicolumn{6}{c}{128 contexts} \\
\cmidrule(lr){2-4}\cmidrule(lr){5-10}
& \multicolumn{3}{c}{Wall time\,/ms $\downarrow$} & \multicolumn{6}{c}{Wall time\,/s $\downarrow$} \\
\cmidrule(lr){2-4}\cmidrule(lr){5-10}
Method & CAP & PORT & NET & \multicolumn{2}{c}{CAP} & \multicolumn{2}{c}{PORT} & \multicolumn{2}{c}{NET} \\
& & & & 1 core & 32 cores & 1 core & 32 cores & 1 core & 32 cores \\
\midrule
Atlas & \textbf{5.63} & \textbf{2.73} & \textbf{38.99} & \textbf{0.74} & \textbf{0.04} & \textbf{0.33} & \textbf{0.03} & \textbf{4.97} & \textbf{0.21} \\
ICEO & 85.08 & 48.80 & 2,898.92 & 11.31 & 0.41 & 6.26 & 0.25 & 244.54 & 10.22 \\
ICEO-red & 5.66 & 2.79 & 39.78 & \textbf{0.74} & \textbf{0.04} & \textbf{0.33} & \textbf{0.03} & 4.99 & 0.22 \\
SOF & 33.50 & 30.80 & 314.15 & 4.82 & 0.47 & 2.91 & 0.41 & 73.11 & 3.56 \\
kNN & 16.34 & 33.22 & 871.77 & 2.11 & 0.23 & 3.86 & 0.30 & 123.50 & 5.45 \\
kNN-red & 9.03 & 10.66 & 43.39 & 1.12 & 0.21 & 1.31 & 0.24 & 5.45 & 0.40 \\
\bottomrule
\end{tabular*}
\par\vspace{2pt}\begin{minipage}{\textwidth}\renewcommand{\baselinestretch}{1.02}\footnotesize
\textit{Notes.} Entries are medians of repeated runs on replication~1.
\end{minipage}
\end{table}

\subsubsection{Selected local structures and model sizes.}
\nopagebreak

Table~\ref{tab:capPortStructures} reports the selected local structures
and model sizes. Across the three applications, Atlas uses only 4--18
scenarios per decision, reducing the 500-scenario reference support by
$96.4\%$--$99.2\%$. Atlas fits 39--800 parameters across all charts, compared with
17,844 for ICEO, a $22$--$458\times$ reduction.

\begin{table}[htbp]
\centering
\caption{Selected Atlas structures, parameter dimensions, and scenario reduction.}
\label{tab:capPortStructures}
\small
\renewcommand{\arraystretch}{1.00}
\setlength{\tabcolsep}{5pt}
\begin{tabular*}{\textwidth}{@{\extracolsep{\fill}}llrrllr@{}}
\toprule
Application & Replications & $\widehat B$ & $P_{\rm Atlas}$ & \shortstack[l]{Dimensions\\$d_b$} & \shortstack[l]{Scenarios\\$K_b$} & \shortstack[r]{Reduction\\/\%} \\
\midrule
\multirow{2}{*}{CAP} & 1, 3--10 & 1 & 39 & (3) & (4) & 99.2 \\
 & 2 & 2 & 78 & (3, 3) & (4, 4) & 99.2 \\
\midrule
\multirow{5}{*}{PORT} & 1, 2, 5, 7, 9 & 1 & 299 & (13) & (14) & 97.2 \\
 & 10 & 1 & 144 & (8) & (9) & 98.2 \\
 & 6 & 1 & 264 & (12) & (13) & 97.4 \\
 & 4, 8 & 1 & 416 & (16) & (17) & 96.6 \\
 & 3 & 1 & 459 & (17) & (18) & 96.4 \\
\midrule
\multirow{3}{*}{NET} & 1 & 4 & 684 & (9, 9, 9, 9) & (10, 10, 10, 10) & 98.0 \\
 & 2 & 2 & 371 & (9, 10) & (10, 11) & 97.8 \\
 & 3 & 4 & 800 & (10, 10, 10, 10) & (11, 11, 11, 11) & 97.8 \\
\bottomrule
\end{tabular*}
\par\vspace{2pt}\begin{minipage}{\textwidth}\renewcommand{\baselinestretch}{1.02}\footnotesize
\textit{Notes.} Parentheses list the per-chart dimensions $(d_b)_{b=1}^{\widehat B}$ and numbers of scenarios $(K_b)_{b=1}^{\widehat B}$ across the $\widehat B$ local charts. Reduction\,/\% $=100(1-\max_b K_b/K_{\rm ref})$.
\end{minipage}
\end{table}

\subsection{Stylized day-ahead grid capacity coverage}
\label{subsec:grid64Empirical}

In this section, we consider a stylized real-data benchmark similar to CAP,
motivated by day-ahead capacity planning in power systems.
Grid operators commit generation capacity in advance to meet uncertain net
demand, and forecast errors in wind and solar generation make these capacity
requirements depend on operating conditions \citep{dumas2023dynamic}.
Committing too little capacity leaves shortages uncovered, while committing too
much increases capacity costs.

We study rolling day-ahead capacity coverage across eight interconnected
European zones using hourly electricity load, wind
generation, and solar generation from Open Power System Data
(version 2020-10-06; \citealt{wiese2019open,opsd2020timeseries})
and the ENTSO-E Transparency Platform \citep{hirth2018entsoe}. Each zone's
net load (electricity demand minus wind and solar generation) is normalized by
its 2015--2016 median load. The operator chooses capacities
$z\in[0,2.5]^8$, with neighboring zones sharing capacity across
interconnecting links. The uncertain net load $\xi\in\RR^{64}$ spans
eight zones and eight three-hour delivery blocks, and the operator observes
a 16-dimensional context $x$ of recent net loads and calendar
features before the day-ahead deadline. Under a specified sharing
rule and shortage-cost objective, the decision problem is
\begin{equation}
 \min_{z\in[0,2.5]^8}\;
 c_0\mathbf 1_8^\top z
+\frac{1}{8}\sum_{r=1}^8\sum_{h=1}^8
c_h^u\,
\E_{P_x}\!\left[
h_\tau\!\left(\xi_{r,h}-[M_{\rm grid}z]_r\right)
\right]+\frac{\eta_{\rm reg}}{2}\|z-z_0(x)\|_2^2.
\label{eq:grid64MainOptimization}
\end{equation}
Here $r$ and $h$ index zones and delivery blocks, respectively, and
$\xi_{r,h}$ is the normalized net load in zone $r$ during block $h$. The
coefficients $c_0>0$ and $c_h^u>0$ are the unit capacity cost and shortage
penalty, respectively. The matrix $M_{\rm grid}\in\RR^{8\times8}$ maps
committed capacities to available capacity after cross-zone sharing.
Section~\ref{ec:grid64Protocol} of the e-companion provides the data-processing
details.

We evaluate the methods over 11 quarterly test windows (2018Q1--2020Q3,
1,004 days) using data from Open Power System Data and 17 subsequent windows
(2020Q4--2024Q4, 1,553 days) using data from the ENTSO-E Transparency Platform.
For each evaluation window, we train on 728 days, validate on the next 364
days, and test on the following quarter.

Table~\ref{tab:grid64Rolling} shows that Atlas uses a single chart with 8--9
scenarios per decision, compared with 728 for ICEO. It achieves comparable
decision quality to ICEO in the first 11 quarters and retains approximately
$90\%$ of the loss reduction of ICEO relative to SAA in the later 17 quarters,
outperforming all other benchmarks in both periods. Atlas also outperforms
ICEO-red and kNN-red at matched support budgets, showing the advantage of
jointly learned local scenarios over geometric scenario reduction. Atlas
achieves approximately $4\times$ speedups in both construction and 32-core
batch deployment relative to ICEO. Sections~\ref{ec:grid64Protocol}
and~\ref{ec:grid64Later} of the e-companion provide detailed results and
experimental protocols.

\begin{table}[htbp]
\centering
\caption{Out-of-sample performance and computational cost in rolling European grid capacity coverage.}
\label{tab:grid64Rolling}
\small
\renewcommand{\arraystretch}{1.00}
\setlength{\tabcolsep}{4pt}
\begin{tabular*}{\textwidth}{@{\extracolsep{\fill}}l*{5}{r}@{}}
\toprule
 & & \multicolumn{2}{c}{Pooled loss (SE) $\downarrow$} & & \\
\cmidrule(lr){3-4}
Method & Scenarios & 2018Q1--2020Q3 & 2020Q4--2024Q4 & Build\,/s $\downarrow$ & \shortstack{128 contexts,\\32 cores\,/ms $\downarrow$} \\
\midrule
Atlas & 8.4 & \textbf{25.10} (0.32) & 24.62 (0.28) & 87.40 & 65.02 \\
ICEO & 728.0 & 25.15 (0.35) & \textbf{24.30} (0.28) & 355.72 & 263.48 \\
ICEO-red & 8.4 & 26.95 (0.35) & 27.22 (0.55) & 81.78 & \textbf{50.29} \\
SOF & 106.4 & 25.38 (0.30) & 25.10 (0.27) & 10.47 & 181.04 \\
kNN & 19.6 & 25.51 (0.31) & 25.07 (0.29) & 2.78 & 68.05 \\
kNN-red & 8.3 & 25.85 (0.37) & 25.28 (0.35) & 60.12 & 77.39 \\
SAA & 728.0 & 27.96 (0.17) & 27.64 (0.17) & $<0.01$ & 132.85 \\
\bottomrule
\end{tabular*}
\par\vspace{2pt}
\begin{minipage}{\textwidth}
\renewcommand{\baselinestretch}{1.02}\footnotesize
\textit{Notes.} Losses report pooled daily means with 14-day Newey--West standard errors (SE). Scenarios average over the first 11 windows; build times and 5-run median deployment times are measured on 2019Q1.
\end{minipage}
\end{table}

\subsection{Ablation and sensitivity studies}
\label{subsec:ablationStudies}

\subsubsection{Number of scenarios.}
\label{subsubsec:supportSizeStudy}
\begin{figure}[htbp]
\centering
\caption{Decision quality versus local support size in CAP.}
\label{fig:fixedKCapCurve}
\includegraphics[width=0.83\textwidth,trim=0 4pt 0 2pt,clip]{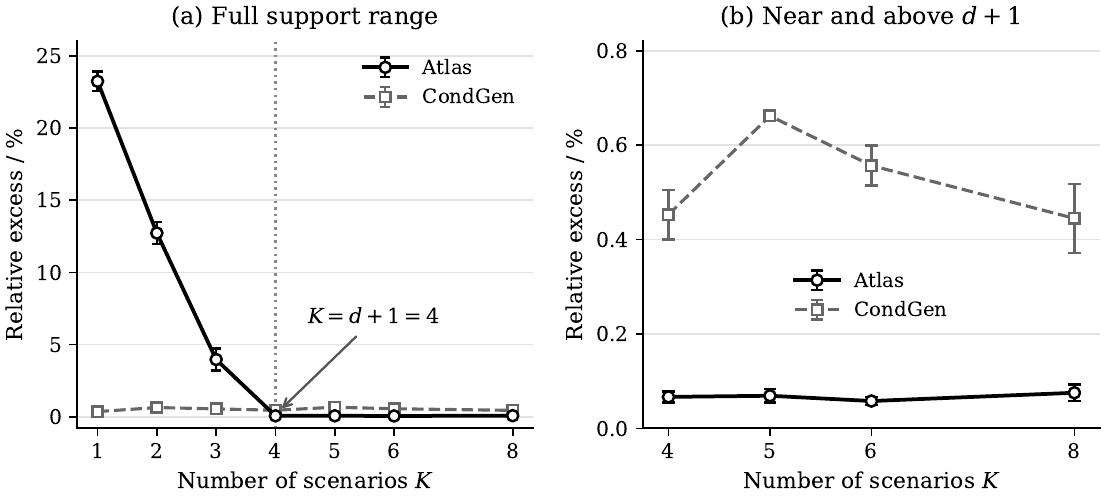}
\par\vspace{6pt}\begin{minipage}{\textwidth}\renewcommand{\baselinestretch}{1.02}\footnotesize
\textit{Notes.} The experiment fixes \(B=1\) and \(d=3\) for Atlas, with three paired training repetitions. The left plot shows mean relative excess risks and standard errors across all support sizes; the right plot provides a closer view of the results near and above \(K=d+1=4\).
\end{minipage}
\end{figure}

To examine how many local scenarios are needed for good decisions, we vary
\(K\in\{1,2,3,4,5,6,8\}\) in CAP while fixing the number of charts at one and
the representation dimension at three, using the same 500 reference
scenarios as in Section~\ref{subsec:applicationResults} and three paired
training repetitions for each support size. We also compare Atlas's shared
support with CondGen, a neural generator that produces context-dependent
scenarios and probabilities. CondGen is trained using the same operational
loss and scenario budget, as detailed in Section~\ref{ec:methodTrainingRules}
of the e-companion.
Figure~\ref{fig:fixedKCapCurve} shows that Atlas's mean relative excess risk
falls from \(23.24\%\) at \(K=1\) to \(12.73\%\) at \(K=2\), \(3.97\%\) at
\(K=3\), and \(0.07\%\) at \(K=d+1=4\), with little change at \(K=5,6,8\). The
plateau begins at \(K=d+1=4\), supporting the choice motivated by
Theorem~\ref{thm:tiltedResponseCapacity}. CondGen achieves \(0.36\%\)
excess risk with only one scenario, but its risk remains between \(0.44\%\)
and \(0.66\%\) for \(K=4\)--\(8\), compared with approximately \(0.07\%\)
for Atlas. CondGen uses 9,891--18,471 parameters over these support sizes,
versus 39--87 for Atlas. Its greater flexibility allows good decisions
with fewer scenarios but also increases the number of parameters to estimate
from limited training data.

\subsubsection{Number of context charts.}
\label{subsec:chartMechanismResults}

We use a controlled quadratic decision problem with four operating modes (each
spanning two scenario directions) to examine when local charts help: we rotate
the mode-specific directions from identical to mutually orthogonal while fixing
contexts and scenario probabilities.

The four modes are represented by vectors
$c_1,\ldots,c_4\in\mathbb R^6$. For an observation from mode $r$, the context
is $x=(c_r+c_e,t)$, where $c_e\in\mathbb R^6$ adds noise to the mode vector and
$t\in\mathbb R^2$ determines the scenario probabilities. Each mode has four
possible scenarios, corresponding to the corners of a square,
$v\in\{-1,1\}^2$. Their probabilities are
$p(v\mid t)=\exp(v^\top t)/\sum_{w\in\{-1,1\}^2}\exp(w^\top t)$.
We map this square into the ten-dimensional uncertainty space using
$\tilde{\xi}_{r,v}=V_r(\omega)v$, where
$V_r(\omega)=\cos(\omega)V_{\rm c}+\sin(\omega)V_r^{\rm aux}$.
The two columns of $V_{\rm c}\in\mathbb R^{10\times2}$ give the shared
directions, and those of $V_r^{\rm aux}\in\mathbb R^{10\times2}$ give the
directions specific to mode $r$. The columns of these five matrices are
mutually orthonormal. At $\omega=0^\circ$, all four squares lie in the same plane, whereas
at $90^\circ$, their planes are mutually orthogonal.

The decision minimizes the expected quadratic loss
$\ell_{\rm chart}(z,\xi)=\frac12\|z-\xi\|_2^2$, so its optimum is the
conditional mean. For an observation from mode $r$,
\begin{equation}
z(P_x;x)=\mathbb E[\xi\mid x]=V_r(\omega)\tanh(t).
\label{eq:chartMechanismMain}
\end{equation}
Here $\tanh$ is applied componentwise, so
$\tanh(t)=(\tanh(t_1),\tanh(t_2))^\top$ for $t=(t_1,t_2)^\top$.
To obtain Eq.~\eqref{eq:chartMechanismMain}, note that
$v^\top t=v_1t_1+v_2t_2$, so the probability of each coordinate factors as
$p(v_j\mid t)=e^{v_jt_j}/(e^{t_j}+e^{-t_j})$ for $j=1,2$.
Since $v_j$ takes values $-1$ and $1$, its conditional mean is
$\mathbb E[v_j\mid t]=(e^{t_j}-e^{-t_j})/(e^{t_j}+e^{-t_j})=\tanh(t_j)$.
Thus, $\mathbb E[\xi\mid x]=V_r(\omega)\mathbb E[v\mid t]
=V_r(\omega)\tanh(t)$, which gives the exact unregularized optimal decision for evaluation.

For each angle and three replications, we compare $B\in\{1,2,4,8,16,32\}$ on
the same 12,800 training observations and 1,600-update budget, selecting $d_b$
from training data with $K_b=d_b+1$ (full protocol in
Section~\ref{ec:chartMechanismProtocol}).

\begin{figure}[htbp]
\centering
\caption{Decision quality as mode-specific scenario planes separate.}
\label{fig:chartSensitivity}
\includegraphics[width=0.77\textwidth,trim=0 8pt 0 4pt,clip]{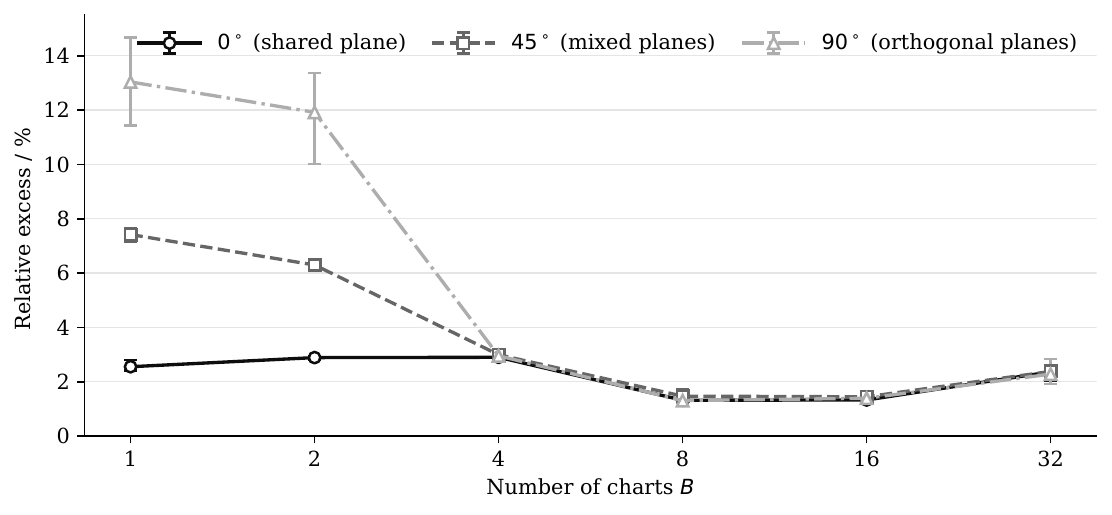}
\par\vspace{2pt}\begin{minipage}{\textwidth}\renewcommand{\baselinestretch}{1.02}\footnotesize
\textit{Notes.} Points show means, and error bars show the range across three replications.
\end{minipage}
\end{figure}

Figure~\ref{fig:chartSensitivity} shows that additional charts become
more useful as the mode-specific planes separate. When $\omega=90^\circ$,
increasing $B$ from 1 to 4 reduces mean relative excess risk from
$13.04\%$ to $2.95\%$ and the selected dimension per chart from
8 to 2, reducing the number of scenarios per decision from 9 to 3.
Increasing $B$ to 8 further reduces relative excess risk to $1.31\%$,
with three scenarios per decision and 24 scenarios in total.
The improvement is smaller when $\omega=0^\circ$, where all modes share the
same plane. Across all three angles, performance levels off between 8 and 16 charts
and worsens at 32, consistent with the localization--estimation tradeoff
in Theorem~\ref{thm:atlasCrossCovarianceTradeoff} and the increased parameter
estimation burden with more charts in
Theorem~\ref{thm:stageTwoGeneralization} ($P_{\rm Atlas}\propto B$).

\subsubsection{Local basis.}
\label{subsubsec:basisAblation}
\label{subsec:controlledBasisResults}

We compare Atlas's cross-covariance basis $\widehat V_b$ against PCA (which
retains largest-variance directions regardless of decision impact) on a
200-dimensional quadratic problem where four directions determine the optimal
decision and eight orthogonal directions add irrelevant variation of scale $\sigma$.

 Let $V_{\rm rel}\in\mathbb R^{200\times4}$
and $V_{\rm noise}\in\mathbb R^{200\times8}$ have mutually orthonormal columns,
with decision-relevant coordinates $u^{\rm rel}\in\mathbb R^4$ and unrelated
coordinates $e^\perp\in\mathbb R^8$. We define
\begin{equation}
 \ell_{\rm quad}(z,\xi)
 =\frac12\|z-V_{\rm rel}^\top\xi\|_2^2,
 \qquad
 \xi^{(\sigma)}
 =V_{\rm rel}u^{\rm rel}+\sigma V_{\rm noise}e^\perp,
 \qquad \sigma\in\{1/4,4\}.
 \label{eq:controlledBasisMain}
\end{equation}
Because $V_{\rm rel}^\top V_{\rm noise}=0$, we have
$V_{\rm rel}^\top\xi^{(\sigma)}=u^{\rm rel}$ and, with $f(z)=\frac12\|z\|_2^2$,
$\phi_{\rm quad}(z,\xi)=-V_{\rm rel}^\top\xi$; increasing $\sigma$ inflates
variance along $V_{\rm noise}$ without changing scenario gradients or optimal
decisions.

Table~\ref{tab:controlledBasisMeans} shows that when unrelated variation
is small ($\sigma=1/4$), both methods recover the decision-relevant subspace
and achieve identical normalized error. When unrelated variation is large
($\sigma=4$), PCA selects high-variance noise directions, reducing subspace
alignment to $0.006$ and increasing normalized error to $0.986$. In contrast,
the cross-covariance basis $\widehat V_b$ maintains an alignment of $0.884$ and
a normalized error of $0.062$, a $93.7\%$ reduction relative to PCA, while
keeping operational risk near its low-noise level.

\begin{table}[htbp]
\centering
\caption{Cross-covariance basis and PCA under decision-irrelevant variation.}
\label{tab:controlledBasisMeans}
\footnotesize
\renewcommand{\arraystretch}{0.95}
\setlength{\tabcolsep}{6pt}
\begin{tabular*}{\textwidth}{@{\extracolsep{\fill}}llrrr@{}}
\toprule
Noise scale $\sigma$ & Basis & Normalized error $\downarrow$ & Alignment $\uparrow$ & Risk $\downarrow$ \\
\midrule
\multirow{2}{*}{$1/4$} & Cross-covariance & 0.041 & 0.999 & 3.4620 \\
 & PCA & 0.041 & 0.999 & 3.4620 \\
\midrule
\multirow{2}{*}{$4$} & Cross-covariance & 0.062 & 0.884 & 3.4843 \\
 & PCA & 0.986 & 0.006 & 4.5003 \\
\bottomrule
\end{tabular*}
\par\vspace{2pt}\begin{minipage}{\textwidth}\renewcommand{\baselinestretch}{1.02}\footnotesize
\textit{Notes.} Entries are means ($\mathcal Z=\RR^4$, identity $T,g$).
Normalized error is $\sum_i\|\hat z(x_i)-z^\star(x_i)\|_2^2/\max(\sum_i\|z^\star(x_i)\|_2^2,10^{-12})$.
Alignment is $\sum_b |I_b|\rho_b/|I_{\rm tr}|$ with
$\rho_b=\sigma_{\min}(V_{\rm rel}^{\top}\widehat V_b^{\rm orig})$, where
$\widehat V_b^{\rm orig}$ orthonormalizes $\widehat D_b\widehat V_b$.
\end{minipage}
\end{table}

Beyond the controlled experiment, Table~\ref{tab:ecBasisAblation} in
Section~\ref{ec:additionalAblations} of the e-companion compares the
scenario bases in the main applications and European grid study.
In CAP and NET, $\widehat V_b$ performs similarly to PCA and full-space
fitting. In PORT, it outperforms PCA ($2.94\%$ versus $3.17\%$) and
approaches full-space fitting ($2.75\%$), which uses $1.6$--$2.5\times$
as many parameters. Across 28 grid windows, $\widehat V_b$ achieves
lower loss than PCA in 20 and full-space fitting in 18.

\subsubsection{Retention threshold.}
\label{subsubsec:retentionSensitivity}
We vary $\alpha_{\rm ret}$ over $\{0.9,0.99,0.999,0.9999\}$.
Table~\ref{tab:ecAlphaSensitivity} shows that lower thresholds reduce
$K_b=d_b+1$ but increase risk. At $0.9$, relative excess risk is
$2.50\%$ in CAP and $11.40\%$ in PORT. At $0.999$, PORT achieves
$3.21\%$, below ICEO's $3.66\%$. The default $0.9999$ achieves
the lowest risk in CAP and PORT and the lowest pooled grid loss
over 28 quarters.

\begin{table}[!htbp]
\centering
\caption{Sensitivity to the spectral retention threshold $\alpha_{\rm ret}$.}
\label{tab:ecAlphaSensitivity}
\footnotesize
\renewcommand{\arraystretch}{0.90}
\setlength{\tabcolsep}{0pt}
\begin{tabular*}{\textwidth}{@{\extracolsep{\fill}}l*{9}{r}@{}}
\toprule
 & \multicolumn{2}{c}{CAP (10 reps.)} & \multicolumn{2}{c}{PORT (10 reps.)} & \multicolumn{5}{c}{European grid (28 quarters, 2,557 days)} \\
\cmidrule(lr){2-3}\cmidrule(lr){4-5}\cmidrule(lr){6-10}
Setting & $\overline K$ & Excess\,/\% & $\overline K$ & Excess\,/\% & $\overline K$ & 18Q1--20Q3 & 20Q4--24Q4 & All 28q & Diff.\ (SE) \\
\midrule
$\alpha_{\rm ret}=0.9$ & 3.0 & 2.50 (0.21) & 4.0 & 11.40 (0.56) & 2.6 & 27.40 & 28.09 & 27.82 & $-3.01$ (0.17) \\
$\alpha_{\rm ret}=0.99$ & 3.1 & 2.33 (0.23) & 7.4 & 4.54 (0.52) & 4.0 & 25.60 & 25.28 & 25.41 & $-0.60$ (0.10) \\
$\alpha_{\rm ret}=0.999$ & 3.9 & 0.52 (0.31) & 8.7 & 3.21 (0.22) & 6.5 & 25.25 & 24.92 & 25.05 & $-0.24$ (0.10) \\
$\alpha_{\rm ret}=0.9999$ (default) & 4.0 & \textbf{0.07} (0.01) & 14.4 & \textbf{2.94} (0.17) & 8.4 & \textbf{25.10} & \textbf{24.62} & \textbf{24.81} &  \\
\bottomrule
\end{tabular*}
\par\vspace{2pt}\begin{minipage}{\textwidth}\renewcommand{\baselinestretch}{1.02}\footnotesize
\textit{Notes.} $\overline K$ is the mean $\max_b K_b$; grid columns report pooled daily means; Diff.\ is default minus row loss over all 28 quarters. Parentheses give cross-replication SEs (CAP, PORT) and 14-day Newey--West SEs (Diff.).
\end{minipage}
\end{table}

\subsubsection{Training sample size.}
\label{subsubsec:trainingSampleSensitivity}
Table~\ref{tab:ecSampleScaling} shows that Atlas achieves lower mean
relative excess risk than all benchmarks at every tested sample size in
CAP and PORT, averaged over three replications. The largest gap to ICEO
occurs in PORT at $1{,}600$ observations ($5.90\%$ versus $18.17\%$),
consistent with the parameter-dimension dependence of the estimation
bound in Theorem~\ref{thm:stageTwoGeneralization}.
\begin{table}[!htbp]
\centering
\caption{Mean relative excess risk (\%) by training sample size $|I_{\rm tr}|$ ($|I_{\rm val}|=1{,}000$, replications 1--3).}
\label{tab:ecSampleScaling}
\footnotesize
\renewcommand{\arraystretch}{0.88}
\setlength{\tabcolsep}{0pt}
\begin{tabular*}{\textwidth}{@{\extracolsep{\fill}}l*{8}{r}@{}}
\toprule
 & \multicolumn{4}{c}{CAP (mean Excess\,/\% and SE)} & \multicolumn{4}{c}{PORT (mean Excess\,/\% and SE)} \\
\cmidrule(lr){2-5}\cmidrule(lr){6-9}
Method & $1{,}600$ & $3{,}200$ & $6{,}400$ & $12{,}800$ & $1{,}600$ & $3{,}200$ & $6{,}400$ & $12{,}800$ \\
\midrule
Atlas & \textbf{0.30} (0.09) & \textbf{0.20} (0.03) & \textbf{0.19} (0.06) & \textbf{0.09} (0.01) & \textbf{5.90} (0.39) & \textbf{3.73} (0.16) & \textbf{3.23} (0.19) & \textbf{2.88} (0.23) \\
ICEO & 1.29 (0.30) & 0.64 (0.15) & 0.48 (0.07) & 0.38 (0.07) & 18.17 (1.24) & 8.09 (0.69) & 5.50 (0.89) & 3.08 (0.09) \\
ICEO-red & 1.12 (0.25) & 0.70 (0.13) & 0.52 (0.09) & 0.42 (0.06) & 13.85 (1.21) & 7.10 (0.66) & 5.39 (0.54) & 4.49 (0.76) \\
SOF & 4.07 (0.17) & 3.35 (0.14) & 2.93 (0.07) & 2.22 (0.16) & 16.89 (1.53) & 14.01 (0.62) & 13.49 (0.71) & 11.18 (0.98) \\
kNN & 3.72 (0.05) & 3.00 (0.15) & 2.40 (0.05) & 1.83 (0.07) & 9.72 (0.90) & 7.37 (0.60) & 6.37 (0.37) & 5.46 (0.44) \\
kNN-red & 3.64 (0.06) & 2.80 (0.11) & 2.33 (0.05) & 1.77 (0.06) & 23.10 (6.94) & 18.60 (7.27) & 17.33 (5.38) & 17.15 (6.01) \\
SAA & 21.64 (0.37) & 21.65 (0.36) & 21.62 (0.40) & 21.59 (0.44) & 17.92 (1.05) & 17.63 (1.15) & 17.59 (1.16) & 17.52 (1.15) \\
\bottomrule
\end{tabular*}
\end{table}

\FloatBarrier
\subsection{Additional comparisons and robustness checks}
\label{subsubsec:round7Robustness}

We conduct additional comparisons and robustness checks, with detailed
results reported in Tables~\ref{tab:ecCondGenAll}--\ref{tab:ecHeadMatched}
of the e-companion.

First, we extend the CondGen comparison from
Section~\ref{subsubsec:supportSizeStudy} to all four applications at the same
scenario budget. Atlas achieves lower relative excess risk in CAP ($0.07\%$
versus $0.49\%$) and PORT ($2.94\%$ versus $5.17\%$), while using substantially
fewer parameters. CondGen achieves slightly lower risk in NET ($0.10\%$
versus $0.12\%$) and lower loss on the European grid, particularly over the
later 17 quarters ($24.22$ versus $24.62$). These results illustrate the
tradeoff between compact shared supports and more flexible context-dependent
scenarios.

Second, we examine whether the probability model explains Atlas's
performance. Replacing the probability models of ICEO and ICEO-red
with the affine-softmax model used by Atlas does not eliminate
Atlas's advantage in the reported comparisons. Replacing Atlas's
affine-softmax model with a neural probability model also increases
operational loss across the tested applications. These results indicate
that Atlas's performance reflects both the learned local scenarios and the
parsimonious probability model, not the probability model alone.

Finally, we examine robustness to training and structure selection.
Using separate samples for chart selection and parameter fitting increases
operational loss, particularly on the grid, where fewer observations are
available for each step; using the full training sample in
Algorithm~\ref{alg:kktAtlas} avoids this loss in sample efficiency.
Constructing local bases from five reference decisions instead of the
single local SAA decision $z_b^{\rm SAA}$ gives similar results
on the synthetic applications and slightly higher loss on the grid,
suggesting little benefit from additional reference decisions in these
settings. Further checks of chart selection, ICEO tuning, and statistical
inference generally support the main findings.

\section{Conclusion}
\label{sec:conclusion}

We develop KKT Atlas, a framework for learning compact local scenario
models for repeated contextual stochastic optimization. Atlas partitions the
context space into charts, identifies local scenario directions from
scenario gradients $\phi(z,\xi;x)$, and jointly fits
scenarios $\tilde{\xi}_{b,k}$ and context-dependent probabilities $p_{b,k}(x)$ by minimizing operational loss.
The spectral rule selects the local representation dimension $d_b$, and the shared-support result
motivates \(K_b=d_b+1\). For a fixed chart structure, the estimation error in the generalization bound scales with the
resulting compact parameter dimension $P_{\rm Atlas}$, without explicit dependence on the reference
support size. For each new context, we then solve only the small
stochastic program associated with its chart.

Atlas attains lower mean risk than all six benchmarks in CAP and PORT and nearly
matches ICEO in NET, using only 4--18 scenarios per decision. It achieves
approximately $3$--$8\times$ construction speedups and $8$--$48\times$
32-core batch deployment speedups over ICEO, while outperforming
clustering-based reduction at matched support budgets in every application.
On the stylized eight-zone European grid, it compresses 728 historical
scenarios to 8--9 learned scenarios per decision, achieves similar
operational loss to ICEO over the first 11 quarterly windows, and achieves approximately
$4\times$ speedups in construction and batch deployment
relative to ICEO. On 17 later quarters evaluated with every setting
fixed in advance, its advantage over SOF, kNN, and the clustering-based
reduction methods persists, and its loss is 1.3\% above that of ICEO.
Additional comparisons with CondGen illustrate the tradeoff between
compact shared supports and more flexible context-dependent scenarios.

These computational savings are particularly useful when downstream
problems are costly and decisions must be made repeatedly, as in network
reservation and rolling operations. Our analysis rests on smoothness
and local curvature. Nonsmooth losses, discrete decisions, and nonconvex downstream models
remain outside this analysis. We also keep the chart partition fixed after training, so the
current construction does not adapt to changes in the contextual
distribution over time. Extending the method to these settings and
updating charts and scenario supports as contexts evolve are directions
for future work. In each case, the central question remains whether
uncertainty can be compressed along the directions that govern the
downstream decision.

\ifnum\value{endnote}>0
  \THEEndNotes
\fi

\FloatBarrier
\bibliographystyle{informs2014}
\bibliography{references}

@inproceedings{agrawal2019differentiable,
  title = {Differentiable {{Convex Optimization Layers}}},
  booktitle = {Advances in {{Neural Information Processing Systems}}},
  author = {Agrawal, Akshay and Amos, Brandon and Barratt, Shane and Boyd, Stephen and Diamond, Steven and Kolter, J. Zico},
  year = 2019,
  volume = {32},
  url = {https://proceedings.neurips.cc/paper/2019/hash/9ce3c52fc54362e22053399d3181c638-Abstract.html}
}

@inproceedings{amos2017optnet,
  title = {{{OptNet}}: {{Differentiable Optimization}} as a {{Layer}} in {{Neural Networks}}},
  booktitle = {Proceedings of the 34th {{International Conference}} on {{Machine Learning}}},
  author = {Amos, Brandon and Kolter, J. Zico},
  year = 2017,
  volume = {70},
  pages = {136--145},
  publisher = {PMLR},
  url = {https://proceedings.mlr.press/v70/amos17a.html}
}

@article{ban2019big,
  title = {The {{Big Data Newsvendor}}: {{Practical Insights}} from {{Machine Learning}}},
  author = {Ban, Gah-Yi and Rudin, Cynthia},
  year = 2019,
  journal = {Operations Research},
  volume = {67},
  number = {1},
  pages = {90--108},
  doi = {10.1287/opre.2018.1757},
  url = {https://doi.org/10.1287/opre.2018.1757}
}

@inproceedings{berthet2020learning,
  title = {Learning with {{Differentiable Perturbed Optimizers}}},
  booktitle = {Advances in {{Neural Information Processing Systems}}},
  author = {Berthet, Quentin and Blondel, Mathieu and Teboul, Olivier and Cuturi, Marco and Vert, Jean-Philippe and Bach, Francis},
  year = 2020,
  volume = {33},
  pages = {9508--9519},
  url = {https://proceedings.neurips.cc/paper/2020/file/6bb56208f672af0dd65451f869fedfd9-Paper.pdf}
}

@article{bertsimas2020predictive,
  title = {From {{Predictive}} to {{Prescriptive Analytics}}},
  author = {Bertsimas, Dimitris and Kallus, Nathan},
  year = 2020,
  journal = {Management Science},
  volume = {66},
  number = {3},
  pages = {1025--1044},
  doi = {10.1287/mnsc.2018.3253},
  url = {https://doi.org/10.1287/mnsc.2018.3253}
}

@article{bertsimas2023optimizationbased,
  title = {Optimization-{{Based Scenario Reduction}} for {{Data-Driven Two-Stage Stochastic Optimization}}},
  author = {Bertsimas, Dimitris and Mundru, Nishanth},
  year = 2023,
  month = jul,
  journal = {Operations Research},
  volume = {71},
  number = {4},
  pages = {1343--1361},
  publisher = {{Institute for Operations Research and the Management Sciences (INFORMS)}},
  issn = {0030-364X, 1526-5463},
  doi = {10.1287/opre.2022.2265},
  url = {https://doi.org/10.1287/opre.2022.2265}
}

@inproceedings{blondel2022efficient,
  title = {Efficient and {{Modular Implicit Differentiation}}},
  booktitle = {Advances in {{Neural Information Processing Systems}}},
  author = {Blondel, Mathieu and Berthet, Quentin and Cuturi, Marco and Frostig, Roy and Hoyer, Stephan and {Llinares-L{\'o}pez}, Felipe and Pedregosa, Fabian and Vert, Jean-Philippe},
  year = 2022,
  volume = {35},
  pages = {5230--5242},
  doi = {10.52202/068431-0378},
  url = {https://proceedings.neurips.cc/paper\_files/paper/2022/hash/228b9279ecf9bbafe582406850c57115-Abstract.html}
}

@book{bonnans2000perturbation,
  title = {Perturbation {{Analysis}} of {{Optimization Problems}}},
  author = {Bonnans, J. Fr{\'e}d{\'e}ric and Shapiro, Alexander},
  year = 2000,
  publisher = {Springer},
  address = {New York},
  doi = {10.1007/978-1-4612-1394-9},
  url = {https://doi.org/10.1007/978-1-4612-1394-9}
}

@book{boyd2004convex,
  title = {Convex Optimization},
  author = {Boyd, Stephen and Vandenberghe, Lieven},
  year = 2004,
  publisher = {Cambridge University Press},
  doi = {10.1017/CBO9780511804441},
  url = {https://web.stanford.edu/\textasciitilde boyd/cvxbook/bv\_cvxbook.pdf}
}

@article{chen2025coupled,
  title = {Coupled {{Input-Output Dimension Reduction}}: {{Application}} to {{Goal}}-Oriented {{Bayesian Experimental Design}} and {{Global Sensitivity Analysis}}},
  author = {Chen, Qiao and Arnaud, {\'E}lise and Baptista, Ricardo and Zahm, Olivier},
  year = 2025,
  month = sep,
  journal = {SIAM Journal on Scientific Computing},
  volume = {47},
  number = {5},
  pages = {A2403-A2430},
  publisher = {Society for Industrial \& Applied Mathematics (SIAM)},
  issn = {1064-8275, 1095-7197},
  doi = {10.1137/24M1670421},
  url = {https://doi.org/10.1137/24m1670421}
}

@inproceedings{donti2017task,
  title = {Task-{{Based End}}-to-{{End Model Learning}} in {{Stochastic Optimization}}},
  booktitle = {Advances in {{Neural Information Processing Systems}}},
  author = {Donti, Priya and Amos, Brandon and Kolter, J. Zico},
  year = 2017,
  volume = {30},
  url = {https://proceedings.neurips.cc/paper/2017/hash/3fc2c60b5782f641f76bcefc39fb2392-Abstract.html}
}

@inproceedings{dumas2023dynamic,
  title = {Dynamic Sizing of Required Balancing Capacities: The Operational Approach in {{France}}},
  booktitle = {2023 19th {{International Conference}} on the {{European Energy Market}} ({{EEM}})},
  author = {Dumas, Jonathan and Terrier, Viktor and Bienvenu, Fr{\'e}d{\'e}ric and Finet, S{\'e}bastien and Grisey, Nathalie},
  year = 2023,
  publisher = {IEEE},
  doi = {10.1109/EEM58374.2023.10161932},
  url = {https://doi.org/10.1109/EEM58374.2023.10161932}
}

@article{dupacova2003scenario,
  title = {Scenario Reduction in Stochastic Programming},
  author = {Dupa{\v c}ov{\'a}, J. and {Gr{\"o}we-Kuska}, N. and R{\"o}misch, W.},
  year = 2003,
  month = mar,
  journal = {Mathematical Programming},
  volume = {95},
  number = {3},
  pages = {493--511},
  issn = {1436-4646},
  doi = {10.1007/s10107-002-0331-0},
  url = {https://doi.org/10.1007/s10107-002-0331-0}
}

@article{elmachtoub2022smart,
  title = {Smart ``Predict, Then Optimize''},
  author = {Elmachtoub, Adam N. and Grigas, Paul},
  year = 2022,
  journal = {Management Science},
  volume = {68},
  number = {1},
  pages = {9--26},
  doi = {10.1287/mnsc.2020.3922},
  url = {https://doi.org/10.1287/mnsc.2020.3922}
}

@article{fairbrother2022problem,
  title = {Problem-Driven Scenario Generation: An Analytical Approach for Stochastic Programs with Tail Risk Measure},
  author = {Fairbrother, Jamie and Turner, Amanda and Wallace, Stein W.},
  year = 2022,
  journal = {Mathematical Programming},
  volume = {191},
  number = {1},
  pages = {141--182},
  publisher = {{Springer Science and Business Media LLC}},
  issn = {0025-5610, 1436-4646},
  doi = {10.1007/s10107-019-01451-7},
  url = {https://doi.org/10.1007/s10107-019-01451-7}
}

@article{fiacco1976sensitivity,
  title = {Sensitivity Analysis for Nonlinear Programming Using Penalty Methods},
  author = {Fiacco, Anthony V.},
  year = 1976,
  journal = {Mathematical Programming},
  volume = {10},
  number = {1},
  pages = {287--311},
  doi = {10.1007/BF01580677},
  url = {https://doi.org/10.1007/BF01580677}
}

@article{goulart2026clarabel,
  title = {Clarabel: {{An}} Interior-Point Solver for Conic Programs with Quadratic Objectives},
  author = {Goulart, Paul J. and Chen, Yuwen},
  year = 2026,
  journal = {Mathematical Programming Computation},
  volume = {18},
  number = {4},
  pages = {1279--1361},
  doi = {10.1007/s12532-026-00320-7}
}

@misc{gurobi,
  title = {Gurobi {{Optimizer Reference Manual}}, {{Version}} 13.0},
  author = {{Gurobi Optimization, LLC}},
  year = 2026,
  url = {https://docs.gurobi.com/\_/downloads/optimizer/en/13.0/pdf/}
}

@article{harrison1999multiresource,
  title = {Multi-Resource Investment Strategies: {{Operational}} Hedging under Demand Uncertainty},
  author = {Harrison, J. Michael and Van Mieghem, Jan A.},
  year = 1999,
  journal = {European Journal of Operational Research},
  volume = {113},
  number = {1},
  pages = {17--29},
  doi = {10.1016/S0377-2217(97)00425-6},
  url = {https://doi.org/10.1016/S0377-2217(97)00425-6}
}

@book{hastie2009elements,
  title = {The Elements of Statistical Learning: {{Data}} Mining, Inference, and Prediction},
  author = {Hastie, Trevor and Tibshirani, Robert and Friedman, Jerome},
  year = 2009,
  edition = {2nd},
  publisher = {Springer},
  address = {New York},
  doi = {10.1007/978-0-387-84858-7},
  url = {https://hastie.su.domains/ElemStatLearn/main.html}
}

@article{heitsch2003scenario,
  title = {Scenario {{Reduction Algorithms}} in {{Stochastic Programming}}},
  author = {Heitsch, Holger and R{\"o}misch, Werner},
  year = 2003,
  month = feb,
  journal = {Computational Optimization and Applications},
  volume = {24},
  number = {2-3},
  pages = {187--206},
  issn = {0926-6003, 1573-2894},
  doi = {10.1023/A:1021805924152},
  url = {https://link.springer.com/10.1023/A:1021805924152}
}

@article{henrion2022problembased,
  title = {Problem-Based Optimal Scenario Generation and Reduction in Stochastic Programming},
  author = {Henrion, R. and R{\"o}misch, W.},
  year = 2022,
  journal = {Mathematical Programming},
  volume = {191},
  number = {1},
  pages = {183--205},
  publisher = {{Springer Science and Business Media LLC}},
  issn = {0025-5610, 1436-4646},
  doi = {10.1007/s10107-018-1337-6},
  url = {https://doi.org/10.1007/s10107-018-1337-6}
}

@article{hirth2018entsoe,
  title = {The {{ENTSO-E Transparency Platform}} -- {{A}} Review of {{Europe}}'s Most Ambitious Electricity Data Platform},
  author = {Hirth, Lion and M{\"u}hlenpfordt, Jonathan and Bulkeley, Marisa},
  year = 2018,
  journal = {Applied Energy},
  volume = {225},
  pages = {1054--1067},
  doi = {10.1016/j.apenergy.2018.04.048},
  url = {https://doi.org/10.1016/j.apenergy.2018.04.048}
}

@article{hoeffding1963probability,
  title = {Probability {{Inequalities}} for {{Sums}} of {{Bounded Random Variables}}},
  author = {Hoeffding, Wassily},
  year = 1963,
  journal = {Journal of the American Statistical Association},
  volume = {58},
  number = {301},
  pages = {13--30},
  doi = {10.1080/01621459.1963.10500830},
  url = {https://doi.org/10.1080/01621459.1963.10500830}
}

@article{homemdemello2026forecasting,
  title = {Forecasting Outside the {{Box}}: Application-Driven Optimal Pointwise Forecasts for a Class of Two-Stage Stochastic Programs},
  shorttitle = {Forecasting Outside the {{Box}}},
  author = {{Homem-de-Mello}, Tito and Valencia, Juan and Lagos, Felipe and Lagos, Guido},
  year = 2026,
  month = may,
  journal = {Mathematical Programming},
  issn = {0025-5610, 1436-4646},
  doi = {10.1007/s10107-026-02359-9},
  url = {https://link.springer.com/10.1007/s10107-026-02359-9}
}

@inproceedings{hornewall2025decision,
  title = {Decision {{Focused Scenario Generation}} for {{Contextual Two-Stage Stochastic Linear Programming}}},
  booktitle = {{{NeurIPS Workshop}} on {{ML}} x {{OR}}},
  author = {Hornewall, Jonathan and {Delannoy-Pavy}, Sol{\`e}ne and Lecl{\`e}re, Vincent and {Homem-de-Mello}, Tito},
  year = 2025,
  url = {https://openreview.net/forum?id=j8BVbta5lg}
}

@article{islip2026contextual,
  title = {Contextual Scenario Generation for Two-Stage Stochastic Programming},
  author = {Islip, David and Kwon, Roy H. and Bae, Sanghyeon and Kim, Woo Chang},
  year = 2026,
  month = nov,
  journal = {Computational Optimization and Applications},
  publisher = {{Springer Science and Business Media LLC}},
  issn = {0926-6003},
  doi = {10.1007/s10589-026-00856-3},
  url = {http://dx.doi.org/10.1007/s10589-026-00856-3}
}

@article{jiang2026optimized,
  title = {Optimized {{Dimensionality Reduction}} for {{Moment-Based Distributionally Robust Optimization}}},
  author = {Jiang, Shiyi and Cheng, Jianqiang and Pan, Kai and Shen, Zuo-Jun Max},
  year = 2026,
  month = may,
  journal = {Operations Research},
  volume = {74},
  number = {3},
  pages = {1648--1669},
  publisher = {{Institute for Operations Research and the Management Sciences (INFORMS)}},
  issn = {0030-364X, 1526-5463},
  doi = {10.1287/opre.2023.0645},
  url = {https://doi.org/10.1287/opre.2023.0645}
}

@article{kallus2023stochastic,
  title = {Stochastic {{Optimization Forests}}},
  author = {Kallus, Nathan and Mao, Xiaojie},
  year = 2023,
  month = apr,
  journal = {Management Science},
  volume = {69},
  number = {4},
  pages = {1975--1994},
  publisher = {{Institute for Operations Research and the Management Sciences (INFORMS)}},
  issn = {0025-1909, 1526-5501},
  doi = {10.1287/mnsc.2022.4458},
  url = {https://doi.org/10.1287/mnsc.2022.4458}
}

@article{kannan2025datadriven,
  title = {Technical {{Note}}---{{Data-Driven Sample Average Approximation}} with {{Covariate Information}}},
  author = {Kannan, Rohit and Bayraksan, G{\"u}zin and Luedtke, James R.},
  year = 2025,
  month = nov,
  journal = {Operations Research},
  volume = {73},
  number = {6},
  pages = {3245--3259},
  publisher = {{Institute for Operations Research and the Management Sciences (INFORMS)}},
  issn = {0030-364X, 1526-5463},
  doi = {10.1287/opre.2020.0533},
  url = {https://doi.org/10.1287/opre.2020.0533}
}

@article{keutchayan2023problem,
  title = {Problem-Driven Scenario Clustering in Stochastic Optimization},
  author = {Keutchayan, Julien and Ortmann, Janosch and Rei, Walter},
  year = 2023,
  journal = {Computational Management Science},
  volume = {20},
  number = {1},
  pages = {13},
  doi = {10.1007/s10287-023-00446-2},
  url = {https://doi.org/10.1007/s10287-023-00446-2}
}

@inproceedings{kingma2015adam,
  title = {Adam: A Method for Stochastic Optimization},
  booktitle = {International Conference on Learning Representations},
  author = {Kingma, Diederik P. and Ba, Jimmy},
  year = 2015,
  url = {https://arxiv.org/abs/1412.6980}
}

@article{mandi2024decision,
  title = {Decision-{{Focused Learning}}: {{Foundations}}, {{State}} of the {{Art}}, {{Benchmark}} and {{Future Opportunities}}},
  author = {Mandi, Jayanta and Kotary, James and Berden, Senne and Mulamba, Maxime and Bucarey, Victor and Guns, Tias and Fioretto, Ferdinando},
  year = 2024,
  journal = {Journal of Artificial Intelligence Research},
  volume = {80},
  pages = {1623--1701},
  doi = {10.1613/jair.1.15320},
  url = {https://doi.org/10.1613/jair.1.15320}
}

@article{newey1987simple,
  title = {A {{Simple}}, {{Positive Semi-Definite}}, {{Heteroskedasticity}} and {{Autocorrelation Consistent Covariance Matrix}}},
  author = {Newey, Whitney K. and West, Kenneth D.},
  year = 1987,
  journal = {Econometrica},
  volume = {55},
  number = {3},
  pages = {703--708},
  doi = {10.2307/1913610},
  url = {https://doi.org/10.2307/1913610}
}

@misc{opsd2020timeseries,
  title = {Data {{Package Time Series}}, {{Version}} 2020-10-06},
  author = {{Open Power System Data}},
  year = 2020,
  doi = {10.25832/time_series/2020-10-06},
  url = {https://data.open-power-system-data.org/time\_series/2020-10-06/}
}

@inproceedings{paszke2019pytorch,
  title = {{{PyTorch}}: {{An}} Imperative Style, High-Performance Deep Learning Library},
  booktitle = {Advances in Neural Information Processing Systems},
  author = {Paszke, Adam and Gross, Sam and Massa, Francisco and Lerer, Adam and Bradbury, James and Chanan, Gregory and Killeen, Trevor and Lin, Zeming and Gimelshein, Natalia and Antiga, Luca and Desmaison, Alban and K{\"o}pf, Andreas and Yang, Edward and DeVito, Zachary and Raison, Martin and Tejani, Alykhan and Chilamkurthy, Sasank and Steiner, Benoit and Fang, Lu and Bai, Junjie and Chintala, Soumith},
  year = 2019,
  volume = {32},
  url = {https://proceedings.neurips.cc/paper/2019/hash/bdbca288fee7f92f2bfa9f7012727740-Abstract.html}
}

@article{qi2026integrated,
  title = {Integrated {{Conditional Estimation-Optimization}}},
  author = {Qi, Meng and Grigas, Paul and Shen, Zuo-Jun (Max)},
  year = 2026,
  month = may,
  journal = {Operations Research},
  volume = {74},
  number = {3},
  pages = {1604--1625},
  publisher = {{Institute for Operations Research and the Management Sciences (INFORMS)}},
  issn = {0030-364X, 1526-5463},
  doi = {10.1287/opre.2023.0427},
  url = {https://doi.org/10.1287/opre.2023.0427}
}

@article{ribeiro2026integrated,
  title = {Integrated Estimate-and-Optimize Decision Trees Learning for Two-Stage Linear Decision-Making Problems},
  author = {Ribeiro, Rafaela and Fanzeres, Bruno},
  year = 2026,
  month = mar,
  journal = {European Journal of Operational Research},
  volume = {329},
  number = {2},
  pages = {607--628},
  publisher = {Elsevier BV},
  issn = {0377-2217},
  doi = {10.1016/j.ejor.2025.08.048},
  url = {https://doi.org/10.1016/j.ejor.2025.08.048}
}

@article{riis2002capacitated,
  title = {Capacitated Network Design with Uncertain Demand},
  author = {Riis, Morten and Andersen, Kim Allan},
  year = 2002,
  journal = {INFORMS Journal on Computing},
  volume = {14},
  number = {3},
  pages = {247--260},
  doi = {10.1287/ijoc.14.3.247.111},
  url = {https://doi.org/10.1287/ijoc.14.3.247.111}
}

@article{rockafellar2000optimization,
  title = {Optimization of Conditional Value-at-Risk},
  author = {Rockafellar, R. Tyrrell and Uryasev, Stanislav},
  year = 2000,
  journal = {The Journal of Risk},
  volume = {2},
  number = {3},
  pages = {21--41},
  issn = {14651211},
  doi = {10.21314/JOR.2000.038},
  url = {http://www.risk.net/journal-of-risk/technical-paper/2161159/optimization-conditional-value-risk}
}

@article{rujeerapaiboon2022scenario,
  title = {Scenario Reduction Revisited: Fundamental Limits and Guarantees},
  author = {Rujeerapaiboon, Napat and Schindler, Kilian and Kuhn, Daniel and Wiesemann, Wolfram},
  year = 2022,
  journal = {Mathematical Programming},
  volume = {191},
  number = {1},
  pages = {207--242},
  publisher = {{Springer Science and Business Media LLC}},
  issn = {0025-5610, 1436-4646},
  doi = {10.1007/s10107-018-1269-1},
  url = {https://doi.org/10.1007/s10107-018-1269-1}
}

@article{sadana2025survey,
  title = {A Survey of Contextual Optimization Methods for Decision-Making under Uncertainty},
  author = {Sadana, Utsav and Chenreddy, Abhilash and Delage, Erick and Forel, Alexandre and Frejinger, Emma and Vidal, Thibaut},
  year = 2025,
  journal = {European Journal of Operational Research},
  volume = {320},
  number = {2},
  pages = {271--289},
  doi = {10.1016/j.ejor.2024.03.020},
  url = {https://doi.org/10.1016/j.ejor.2024.03.020}
}

@article{stratigakos2022prescriptive,
  title = {Prescriptive {{Trees}} for {{Integrated Forecasting}} and {{Optimization Applied}} in {{Trading}} of {{Renewable Energy}}},
  author = {Stratigakos, Akylas and Camal, Simon and Michiorri, Andrea and Kariniotakis, Georges},
  year = 2022,
  month = nov,
  journal = {IEEE Transactions on Power Systems},
  volume = {37},
  number = {6},
  pages = {4696--4708},
  publisher = {{Institute of Electrical and Electronics Engineers (IEEE)}},
  issn = {0885-8950, 1558-0679},
  doi = {10.1109/TPWRS.2022.3152667},
  url = {https://doi.org/10.1109/tpwrs.2022.3152667}
}

@article{tropp2012userfriendly,
  title = {User-{{Friendly Tail Bounds}} for {{Sums}} of {{Random Matrices}}},
  author = {Tropp, Joel A.},
  year = 2012,
  journal = {Foundations of Computational Mathematics},
  volume = {12},
  number = {4},
  pages = {389--434},
  doi = {10.1007/s10208-011-9099-z},
  url = {https://doi.org/10.1007/s10208-011-9099-z}
}

@article{virtanen2020scipy,
  title = {{{SciPy}} 1.0: {{Fundamental}} Algorithms for Scientific Computing in Python},
  author = {Virtanen, Pauli and others},
  year = 2020,
  journal = {Nature Methods},
  volume = {17},
  number = {3},
  pages = {261--272},
  doi = {10.1038/s41592-019-0686-2}
}

@article{wiese2019open,
  title = {Open {{Power System Data}} -- {{Frictionless}} Data for Electricity System Modelling},
  author = {Wiese, Frauke and Schlecht, Ingmar and Bunke, Wolf-Dieter and Gerbaulet, Clemens and Hirth, Lion and Jahn, Martin and Kunz, Friedrich and Lorenz, Casimir and M{\"u}hlenpfordt, Jonathan and Reimann, Juliane and Schill, Wolf-Peter},
  year = 2019,
  journal = {Applied Energy},
  volume = {236},
  pages = {401--409},
  doi = {10.1016/j.apenergy.2018.11.097},
  url = {https://doi.org/10.1016/j.apenergy.2018.11.097}
}

@inproceedings{wilder2019melding,
  title = {Melding the {{Data-Decisions Pipeline}}: {{Decision-Focused Learning}} for {{Combinatorial Optimization}}},
  booktitle = {Proceedings of the {{AAAI Conference}} on {{Artificial Intelligence}}},
  author = {Wilder, Bryan and Dilkina, Bistra and Tambe, Milind},
  year = 2019,
  volume = {33},
  pages = {1658--1665},
  doi = {10.1609/aaai.v33i01.33011658},
  url = {https://doi.org/10.1609/aaai.v33i01.33011658}
}

@article{zahm2020gradient,
  title = {Gradient-{{Based Dimension Reduction}} of {{Multivariate Vector-Valued Functions}}},
  author = {Zahm, Olivier and Constantine, Paul G. and Prieur, Cl{\'e}mentine and Marzouk, Youssef M.},
  year = 2020,
  month = jan,
  journal = {SIAM Journal on Scientific Computing},
  volume = {42},
  number = {1},
  pages = {A534-A558},
  publisher = {Society for Industrial \& Applied Mathematics (SIAM)},
  issn = {1064-8275, 1095-7197},
  doi = {10.1137/18M1221837},
  url = {https://doi.org/10.1137/18m1221837}
}

@article{zhang2023optimized,
  title = {Optimized {{Scenario Reduction}}: {{Solving Large-Scale Stochastic Programs}} with {{Quality Guarantees}}},
  author = {Zhang, Wei and Wang, Kai and Jacquillat, Alexandre and Wang, Shuaian},
  year = 2023,
  month = jul,
  journal = {INFORMS Journal on Computing},
  volume = {35},
  number = {4},
  pages = {886--908},
  publisher = {{Institute for Operations Research and the Management Sciences (INFORMS)}},
  issn = {1091-9856, 1526-5528},
  doi = {10.1287/ijoc.2023.1295},
  url = {https://doi.org/10.1287/ijoc.2023.1295}
}

@article{zhang2026piecewiseaffine,
  title = {Data-{{Driven Piecewise Affine Decision Rules}} for {{Stochastic Programming}} with {{Covariate Information}}},
  author = {Zhang, Yiyang and Liu, Junyi and Zhao, Xiaobo},
  year = 2026,
  month = jun,
  journal = {Operations Research},
  volume = {74},
  number = {5},
  pages = {2629--2649},
  publisher = {{Institute for Operations Research and the Management Sciences (INFORMS)}},
  issn = {0030-364X, 1526-5463},
  doi = {10.1287/opre.2023.0175},
  url = {https://doi.org/10.1287/opre.2023.0175}
}

\ECSwitch
\renewcommand{\thepage}{EC.\arabic{page}}
\makeatletter
\renewcommand{\AAnormalsizeXI}{%
 \@setfontsize\normalsizeXI{11}{11.5}%
 \Dskips{0.45}{0.15}{0.1}{0.35}%
 \let\@listi\@listI
 \baselineskip12.7pt plus 0.3pt minus 0.2pt
 \@bls=12.7pt\@em=10pt\relax}
\normalsize
\setlength{\jot}{2pt}
\renewcommand{\theHsection}{EC.\arabic{section}}
\renewcommand{\theHsubsection}{EC.\arabic{section}.\arabic{subsection}}
\renewcommand{\theHequation}{EC.\arabic{equation}}
\@ifundefined{theHtable}{}{\renewcommand{\theHtable}{EC.\arabic{table}}}
\@ifundefined{theHfigure}{}{\renewcommand{\theHfigure}{EC.\arabic{figure}}}
\@ifundefined{theHtheorem}{}{\renewcommand{\theHtheorem}{EC.\arabic{theorem}}}
\@ifundefined{theHlemma}{}{\renewcommand{\theHlemma}{EC.\arabic{lemma}}}
\@ifundefined{theHproposition}{}{\renewcommand{\theHproposition}{EC.\arabic{proposition}}}
\@ifundefined{theHcorollary}{}{\renewcommand{\theHcorollary}{EC.\arabic{corollary}}}
\@ifundefined{theHassumption}{}{\renewcommand{\theHassumption}{EC.\arabic{assumption}}}
\def\theHdefinition{EC.\arabic{definition}}
\makeatother
\ECHead{E-Companion}

\begingroup
\small
\renewcommand{\arraystretch}{1.01}
\setlength{\tabcolsep}{4.5pt}
\begin{longtable}{@{}
>{\raggedright\arraybackslash}p{0.16\textwidth}
>{\raggedright\arraybackslash}p{0.31\textwidth}
>{\raggedright\arraybackslash}p{0.47\textwidth}@{}}
\caption{Summary of notation.}\label{tab:notationSummary}\\
\toprule
Object family & Notation & Role \\
\midrule
\endfirsthead
\multicolumn{3}{@{}l}{\textit{Table~\thetable\ continued}}\\
\toprule
Object family & Notation & Role \\
\midrule
\endhead
\multicolumn{3}{r@{}}{\textit{Continued on next page}}\\
\endfoot
\bottomrule
\endlastfoot

\multirow{10}{0.16\textwidth}{Data and dimensions}
& $x\in\cX\subseteq\RR^{d_x}$ & Observed or new context. \\
& $\xi\in\Xi\subseteq\RR^{d_\xi}$ & Uncertain parameter in uncertainty space $\Xi$. \\
& $\xi_i,\tilde{\xi}_k\in\Xi$ & Observed realization $i$ and scenario $k$, respectively. \\
& $I_{\rm tr}$ & Index set of training observations $\{(x_i,\xi_i)\}_{i\in I_{\rm tr}}$. \\
& $K_{\rm ref}$ & Size of a reference scenario support. \\
& $p_s^\star(x)$ & True conditional probability of scenario $s$. \\
& $z\in\cZ\subseteq\RR^{d_z}$ & Decision in feasible decision set $\cZ$. \\
& $P_x$ & Conditional distribution of $\xi$ given context $x$. \\
& $d_b$ & Selected dimension of chart $b$. \\
& $d_b^\star$ & Structural dimension of chart $b$. \\
\midrule

\multirow{7}{0.16\textwidth}{Charts}
& $B$ & Candidate number of charts. \\
& $\widehat B$ & Selected number of charts. \\
& $\tau_B(x)$ & Assignment of context $x$ to a chart. \\
& $\bar x_b$ & Centroid of chart $b$. \\
& $I_b\subseteq I_{\rm tr}$ & Training indices assigned to chart $b$. \\
& $\widehat D_x$ & Diagonal matrix of context sample scales. \\
& $\|h\|_x=\|\widehat D_x^{-1}h\|_2$ & Standardized norm for context differences. \\
\midrule

\multirow{8}{0.16\textwidth}{Representations}
& $T:\Xi\to\RR^{d_\xi}$ & Transformation of realizations. \\
& $g:\operatorname{dom}(g)\to\Xi$ & Decoder mapping coordinates to $\Xi$. \\
& $\widehat\mu_b$ & Mean of the transformed realizations in chart $b$. \\
& $\widehat D_b$ & Diagonal matrix of transformed sample scales. \\
& $\varepsilon_{\rm num}$ & Numerical lower bound on scales in $\widehat D_b$. \\
& $T_b(\xi)=\widehat D_b^{-1}\{T(\xi)-\widehat\mu_b\}$ & Standardized chart encoder. \\
& $g_b(y)=g(\widehat\mu_b+\widehat D_b y)$ & Chart decoder. \\
& $r_{y,b}\in(0,\infty]$ & Coordinate radius in the affine-gradient analysis. \\
\midrule

\multirow{4}{0.16\textwidth}{Gradient information}
& $\phi(z,\xi;x)=\nabla_zQ(z,\xi;x)$ & Scenario gradient. \\
& $\bar\phi(z;x)=\E_{P_x}[\phi(z,\xi;x)]$ & Conditional mean scenario gradient. \\
& $\widehat C_b\in\RR^{d_\xi\times d_z}$ & Empirical cross-covariance matrix. \\
& $C_b^\star\in\RR^{d_\xi\times d_z}$ & Ideal affine-gradient signal. \\
\midrule

\multirow{6}{0.16\textwidth}{Bases and coordinates}
& $\widehat V_{b,d}$ & First $d$ left singular directions of $\widehat C_b$. \\
& $\widehat V_b=\widehat V_{b,d_b}$ & Local basis for chart $b$. \\
& $V_b^{\rm rel}$ & Basis of the gradient-relevant space. \\
& $V_b^\star$ & Basis of the ideal covariance signal space. \\
& $u_{b,k}$ & Coordinates of generated scenario $k$ in chart $b$. \\
& $u^{\rm rel}$ & Coordinates in the gradient-relevant space. \\
\midrule

\multirow{7}{0.16\textwidth}{Local scenario model}
& $\tilde{\xi}_{b,k}=g_b(\widehat V_bu_{b,k})$ & Generated scenario $k$ in chart $b$. \\
& $p_b(x)\in\Delta_{K_b}$ & Context-dependent scenario probabilities. \\
& $K_b=d_b+1$ & Number of scenarios used by the algorithm. \\
& $a_b\in\RR^{K_b}$ & Affine-softmax intercepts. \\
& $\Gamma_b\in\RR^{K_b\times d_x}$ & Affine-softmax context coefficients. \\
& $\theta_b=(u_{b,1},\ldots,u_{b,K_b},a_b,\Gamma_b)$ & Parameters fitted for chart $b$. \\
& $P_{\rm Atlas}$ & Total free parameter dimension across charts. \\
\midrule

\multirow{6}{0.16\textwidth}{Loss and objectives}
& $\ell=f+Q$ & Operational loss (deterministic cost $f$ plus scenario cost $Q$). \\
& $R_0(z;x)=\E_{P_x}[\ell(z,\xi;x)]$ & Conditional operational risk. \\
& $R(z;P,x)$ & Regularized objective under distribution $P$ at context $x$. \\
& $\widehat R_b(z;x)$ & Regularized fitted objective for chart $b$. \\
& $\nabla_zR(z;P_x,x)$ & Population first-order mapping. \\
& $\nabla_z\widehat R_b(z;x)$ & Fitted first-order mapping. \\
\midrule

\multirow{13}{0.16\textwidth}{Decisions and bounds}
& $z_0(x)$ & Prespecified regularization center. \\
& $z(P_x;x)$ & Minimizer of $R(\cdot;P_x,x)$. \\
& $z(\widetilde P_x;x)$ & Decision under the finite scenario distribution $\widetilde P_x$. \\
& $z_b(x;\theta_b)$ & Decision under chart parameters $\theta_b$. \\
& $\widehat z_b(x)$ & Decision under fitted chart parameters. \\
& $\pi_{\rm Atlas}(x;\theta)$ & Atlas decision rule for parameters $\theta$. \\
& $\eta_{\rm reg}\geq0$ & Added quadratic regularization coefficient. \\
& $\eta>0$ & Local strong-convexity constant of $R(\cdot;P_x,x)$. \\
& $r_{\rm loc}\in(0,\infty]$ & Radius of $\mathcal B_{\cZ}(z(P_x;x),r_{\rm loc})$. \\
& $\cZ_b$ & Local decision set in chart $b$. \\
& $\varepsilon$ with role superscripts & Error quantities distinguished by role. \\
& $\bar\varepsilon_B^{\rm aff}$ & Sample-weighted average affine-gradient residual. \\
& $\overline\varepsilon_b^C(\delta)$ & Upper bound on $\varepsilon_b^C$. \\
\end{longtable}
\endgroup

\section{Proofs and Auxiliary Results}
\label{ec:proofs}

\subsection{Cross-covariance estimation}
\label{ec:finiteSampleRecovery}
\label{ec:proofFiniteSampleCrossCovariance}

\begin{theorem}[Matrix Bernstein inequality, {\citealp[Theorem~1.6]{tropp2012userfriendly}}]
\label{thm:ecMatrixBernsteinSource}
Let $I$ be a finite index set and let $Z_i$, $i\in I$, be independent, self-adjoint random matrices in
$\RR^{d_{\rm mat}\times d_{\rm mat}}$ satisfying $\E[Z_i]=0$ and
$\|Z_i\|_{\rm op}\leq L$ almost surely, where $L>0$. Then, for every $t>0$,
\begin{equation}
 \mathbb P\left(\left\|\sum_{i\in I}Z_i\right\|_{\rm op}\geq t\right)
 \leq 2d_{\rm mat}\exp\left\{-\frac{t^2}{2\|\sum_{i\in I}\E[Z_i^2]\|_{\rm op}+(2L/3)t}\right\}.
 \label{eq:ecMatrixBernsteinSource}
\end{equation}
\end{theorem}

\begin{lemma}
\label{lem:ecMatrixBernstein}
Let $I$ be a finite index set and let \(\mathcal F\) be a sigma-algebra.
Conditional on \(\mathcal F\), let $Z_i$, $i\in I$, be independent, symmetric random matrices in
$\RR^{d_{\rm mat}\times d_{\rm mat}}$ satisfying
\(\E[Z_i\mid\mathcal F]=0\).
Suppose the nonnegative \(\mathcal F\)-measurable quantities \(L\) and \(v\)
satisfy \(\|Z_i\|_{\rm op}\leq L\) and
\(\|\sum_{i\in I}\E[Z_i^2\mid\mathcal F]\|_{\rm op}\leq v\) almost surely.
Then, for every $u>0$,
\begin{equation}
 \mathbb P\left(
 \left\|\sum_{i\in I}Z_i\right\|_{\rm op}
 >
 \sqrt{2vu}+\frac{2Lu}{3}
 \,\middle|\,\mathcal F
 \right)
 \leq 2d_{\rm mat}\exp(-u).
 \label{eq:ecMatrixBernstein}
\end{equation}
\end{lemma}

\proof{Proof.}
Under the conditional distribution,
if \(L=0\), all \(Z_i\) are zero matrices, so
Eq.~\eqref{eq:ecMatrixBernstein} holds. Otherwise, the matrices satisfy the hypotheses of
Theorem~\ref{thm:ecMatrixBernsteinSource}, with norm bound \(L\) and
variance parameter at most \(v\).
Set \(t\) to \(\sqrt{2vu}+2Lu/3>0\) in
Theorem~\ref{thm:ecMatrixBernsteinSource}. Expanding the square gives
\(t^2-u\{2v+(2L/3)t\}=(2Lu/3)\sqrt{2vu}\geq0\).
Thus, \(t^2/\{2v+(2L/3)t\}\geq u\), and we also have \(\|\sum_{i\in I}\E[Z_i^2\mid\mathcal F]\|_{\rm op}\leq v\) almost surely. Hence,
 \(\frac{t^2}{2\|\sum_{i\in I}\E[Z_i^2\mid\mathcal F]\|_{\rm op}+(2L/3)t} \geq \frac{t^2}{2v+(2L/3)t} \geq u\). Substitution proves
Eq.~\eqref{eq:ecMatrixBernstein} almost surely. \Halmos
\endproof

For a candidate partition and chart \(b\), write
\(\bar\mu_{T,b}:=|I_b|^{-1}\sum_{i\in I_b}\mu_T(x_i)\).
\begin{definition}[Cross-covariance error terms]
For \(\delta\in(0,1)\), define
\begin{align}
 \varepsilon_b^{\rm loc}
 &:=\left\|\frac{1}{|I_b|}\sum_{i\in I_b}
 \left[\Sigma_T(x_i)+\{\mu_T(x_i)-\bar\mu_{T,b}\}
 \{\mu_T(x_i)-\bar\mu_{T,b}\}^{\top}\right]
 -\Sigma_T(\bar x_b)\right\|_{\rm op},
 \label{eq:chartMixtureBias}\\
 \varepsilon_b^{\rm samp}(\delta)
 &:=R_T^2\left(3\sqrt{\frac{2\log\{4(d_\xi+1)/\delta\}}{|I_b|}}
 +\frac{10\log\{4(d_\xi+1)/\delta\}}{3|I_b|}\right),
 \label{eq:chartSamplingRadius}\\
 \overline\varepsilon_b^C(\delta)
 &:=\overline H_b\|\widehat D_b^{-1}\|_{\rm op}^2
 \{\varepsilon_b^{\rm loc}+\varepsilon_b^{\rm samp}(\delta)\}
 +2R_T\|\widehat D_b^{-1}\|_{\rm op}\varepsilon_b^{\rm aff}.
 \label{eq:finiteSampleCrossCovarianceRadius}
\end{align}
\end{definition}

\begin{lemma}
\label{prop:finiteSampleCrossCovariance}
Fix a candidate partition satisfying Assumptions~\ref{ass:regularizedProblem},
\ref{ass:differentiableScenarioCosts}, and~\ref{ass:finiteSampleCovariance}
in Sections~\ref{sec:problem} and~\ref{subsec:chartTradeoff} and a chart \(b\) with \(|I_b|\geq2\).
Then, for every \(\delta\in(0,1)\), the following bound holds with probability
at least \(1-\delta\) conditional on the training contexts,
\begin{equation}
 \varepsilon_b^C
 =
 \bigl\|\widehat C_b-C_b^\star\bigr\|_{\rm op}
 \leq
 \overline\varepsilon_b^C(\delta).
 \label{eq:finiteSampleCrossCovarianceBound}
\end{equation}
\end{lemma}

\proof{Proof.}
All expectations in this proof are conditional on the training contexts $(x_i)_{i\in I_{\rm tr}}$.

\emph{Step 1. Relating cross-covariance error to covariance error.}
Let \(W_i:=T(\xi_i)\),
\(\widehat\Sigma_W:=|I_b|^{-1}\sum_{i\in I_b}
(W_i-\widehat\mu_b)(W_i-\widehat\mu_b)^\top\),
\(Y_i:=T_b(\xi_i)=\widehat D_b^{-1}(W_i-\widehat\mu_b)\),
\(\widehat\Sigma_Y:=|I_b|^{-1}\sum_{i\in I_b}Y_iY_i^\top\), and
\(H:=H_b(z_b^{\rm SAA},\bar x_b)\).
Since \(\widehat\mu_b=|I_b|^{-1}\sum_{i\in I_b}W_i\), we have
\(\sum_{i\in I_b}Y_i
=\widehat D_b^{-1}(\sum_{i\in I_b}W_i-|I_b|\widehat\mu_b)=0\).
Assumption~\ref{ass:finiteSampleCovariance}(\ref{cond:chartwiseSupportEnvelope}) gives
\(\|W_i\|_2\leq R_T\) and \(g_b(Y_i)=\xi_i\), and
Assumption~\ref{ass:finiteSampleCovariance}(\ref{cond:chartwiseAffineRepresentation}) with Eq.~\eqref{eq:targetAffineResponse} gives
\(\phi(z_b^{\rm SAA},\xi_i;\bar x_b)
=\phi(z_b^{\rm SAA},g_b(Y_i);\bar x_b)
=\beta_b(z_b^{\rm SAA},\bar x_b)+H(V_b^{\rm rel})^\top Y_i
+r_b(z_b^{\rm SAA},Y_i,\bar x_b)\).
Define \(R_b^{\rm aff}:=|I_b|^{-1}\sum_{i\in I_b}
Y_ir_b(z_b^{\rm SAA},Y_i,\bar x_b)^\top\).
Substitution into Eq.~\eqref{eq:gradientCrossCovariance} cancels the
intercept because \(\sum_{i\in I_b}Y_i=0\):
\begin{align*}
 \widehat C_b
 =\frac{1}{|I_b|}\sum_{i\in I_b}
 Y_i\phi(z_b^{\rm SAA},\xi_i;\bar x_b)^\top
 &=\left(\frac{1}{|I_b|}\sum_{i\in I_b}Y_i\right)
 \beta_b(z_b^{\rm SAA},\bar x_b)^\top
 +\widehat\Sigma_YV_b^{\rm rel}H^\top+R_b^{\rm aff} \\
 &=\widehat\Sigma_YV_b^{\rm rel}H^\top+R_b^{\rm aff}.
\end{align*}
Since \(\|\widehat\mu_b\|_2\leq|I_b|^{-1}\sum_{j\in I_b}\|W_j\|_2\leq R_T\),
the radius bound in Eq.~\eqref{eq:chartwiseRegularity} gives
\(\|Y_i\|_2=\|\widehat D_b^{-1}(W_i-\widehat\mu_b)\|_2
\leq\|\widehat D_b^{-1}\|_{\rm op}(\|W_i\|_2+\|\widehat\mu_b\|_2)
\leq2R_T\|\widehat D_b^{-1}\|_{\rm op}\leq r_{y,b}\) almost surely.
Hence Eq.~\eqref{eq:affineResidualBound}
gives \(\|r_b(z_b^{\rm SAA},Y_i,\bar x_b)\|_2\leq\varepsilon_b^{\rm aff}\)
almost surely, so the rank-one outer products satisfy
\(\|Y_ir_b(z_b^{\rm SAA},Y_i,\bar x_b)^\top\|_{\rm op}
=\|Y_i\|_2\|r_b(z_b^{\rm SAA},Y_i,\bar x_b)\|_2
\leq2R_T\|\widehat D_b^{-1}\|_{\rm op}\varepsilon_b^{\rm aff}\) almost surely,
and the triangle inequality gives
\(\|R_b^{\rm aff}\|_{\rm op}\leq2R_T\|\widehat D_b^{-1}\|_{\rm op}\varepsilon_b^{\rm aff}\).
Using Eq.~\eqref{eq:idealSignalDefinition}, we obtain
\begin{equation}
 \widehat C_b-C_b^\star
 =(\widehat\Sigma_Y-\Sigma_b)V_b^{\rm rel}H^\top+R_b^{\rm aff}.
 \label{eq:ecCrossCovarianceDecomposition}
\end{equation}
Thus,
\begin{equation}
 \|\widehat C_b-C_b^\star\|_{\rm op}
 \leq\|\widehat\Sigma_Y-\Sigma_b\|_{\rm op}
 \|V_b^{\rm rel}\|_{\rm op}\|H^\top\|_{\rm op}+\|R_b^{\rm aff}\|_{\rm op}
 \leq\overline H_b\|\widehat\Sigma_Y-\Sigma_b\|_{\rm op}
 +2R_T\|\widehat D_b^{-1}\|_{\rm op}\varepsilon_b^{\rm aff},
 \label{eq:ecCrossCovarianceReduction}
\end{equation}
where the first inequality follows from
Eq.~\eqref{eq:ecCrossCovarianceDecomposition}, and the second uses
the orthonormal columns of \(V_b^{\rm rel}\) and
\(\|H^\top\|_{\rm op}=\|H\|_{\rm op}\leq\overline H_b\) from
Eq.~\eqref{eq:affineResidualBound}.

Diagonality of \(\widehat D_b\) gives
\(\widehat\Sigma_Y=\widehat D_b^{-1}\widehat\Sigma_W\widehat D_b^{-1}\).
Since \(\Sigma_b=\widehat D_b^{-1}\Sigma_T(\bar x_b)\widehat D_b^{-1}\),
\begin{equation}
 \begin{aligned}
 \widehat\Sigma_Y-\Sigma_b
 &=\widehat D_b^{-1}\{\widehat\Sigma_W-\Sigma_T(\bar x_b)\}\widehat D_b^{-1},\\
 \|\widehat\Sigma_Y-\Sigma_b\|_{\rm op}
 &\leq\|\widehat D_b^{-1}\|_{\rm op}^2
 \|\widehat\Sigma_W-\Sigma_T(\bar x_b)\|_{\rm op}.
 \end{aligned}
 \label{eq:ecCovarianceScaling}
\end{equation}
It remains to bound \(\|\widehat\Sigma_W-\Sigma_T(\bar x_b)\|_{\rm op}\).

\emph{Step 2. Bounding the covariance error.}
Expanding the definition of \(\widehat\Sigma_W\) gives
\(\widehat\Sigma_W=|I_b|^{-1}\sum_{i\in I_b}W_iW_i^\top
-\widehat\mu_b\widehat\mu_b^\top\).
The identities
\(\E[W_iW_i^\top]=\Sigma_T(x_i)+\mu_T(x_i)\mu_T(x_i)^\top\)
and \(\sum_{i\in I_b}\mu_T(x_i)=|I_b|\bar\mu_{T,b}\) give
\[
 \begin{aligned}
 \frac{1}{|I_b|}\sum_{i\in I_b}\E[W_iW_i^\top]
 -\bar\mu_{T,b}\bar\mu_{T,b}^\top
 &=\frac{1}{|I_b|}\sum_{i\in I_b}
 \{\Sigma_T(x_i)+\mu_T(x_i)\mu_T(x_i)^\top\}
 -\bar\mu_{T,b}\bar\mu_{T,b}^\top\\
 &=\frac{1}{|I_b|}\sum_{i\in I_b}\Sigma_T(x_i)
 +\frac{1}{|I_b|}\sum_{i\in I_b}
 \{\mu_T(x_i)-\bar\mu_{T,b}\}\{\mu_T(x_i)-\bar\mu_{T,b}\}^\top.
 \end{aligned}
\]
Let \(Z_i:=W_iW_i^\top-\E[W_iW_i^\top]\).
Adding and subtracting the average second moment and the outer product
of \(\bar\mu_{T,b}\) gives
\begin{equation}
 \begin{aligned}
 \widehat\Sigma_W-\Sigma_T(\bar x_b)
 &=\frac{1}{|I_b|}\sum_{i\in I_b}
 \{W_iW_i^\top-\E[W_iW_i^\top]\}
 +\frac{1}{|I_b|}\sum_{i\in I_b}\E[W_iW_i^\top]
 -\widehat\mu_b\widehat\mu_b^\top-\Sigma_T(\bar x_b)\\
 &=\frac{1}{|I_b|}\sum_{i\in I_b}Z_i
 -\{\widehat\mu_b\widehat\mu_b^\top
 -\bar\mu_{T,b}\bar\mu_{T,b}^\top\}\\
 &\quad+\left[
 \frac{1}{|I_b|}\sum_{i\in I_b}\E[W_iW_i^\top]
 -\bar\mu_{T,b}\bar\mu_{T,b}^\top-\Sigma_T(\bar x_b)
 \right].
 \end{aligned}
 \label{eq:ecCovarianceErrorIdentity}
\end{equation}
Taking operator norms in Eq.~\eqref{eq:ecCovarianceErrorIdentity} and
applying the triangle inequality yields
\begin{equation}
 \begin{aligned}
 \|\widehat\Sigma_W-\Sigma_T(\bar x_b)\|_{\rm op}
 &\leq\frac{1}{|I_b|}\left\|\sum_{i\in I_b}Z_i\right\|_{\rm op}
 +\|\widehat\mu_b\widehat\mu_b^\top
 -\bar\mu_{T,b}\bar\mu_{T,b}^\top\|_{\rm op}\\
 &\quad+\left\|
 \frac{1}{|I_b|}\sum_{i\in I_b}\E[W_iW_i^\top]
 -\bar\mu_{T,b}\bar\mu_{T,b}^\top-\Sigma_T(\bar x_b)
 \right\|_{\rm op}\\
 &=\frac{1}{|I_b|}\left\|\sum_{i\in I_b}Z_i\right\|_{\rm op}
 +\|\widehat\mu_b\widehat\mu_b^\top
 -\bar\mu_{T,b}\bar\mu_{T,b}^\top\|_{\rm op}
 +\varepsilon_b^{\rm loc}.
 \end{aligned}
 \label{eq:ecRawCovarianceDecomposition}
\end{equation}
By Eqs.~\eqref{eq:ecCrossCovarianceReduction} and
\eqref{eq:ecCovarianceScaling}, it now suffices to bound the sum of
the two random terms in Eq.~\eqref{eq:ecRawCovarianceDecomposition}
by \(\varepsilon_b^{\rm samp}(\delta)\) with conditional probability
at least \(1-\delta\).

\emph{Step 3. Bounding the second moment error.}
Let \(t:=\log\{4(d_\xi+1)/\delta\}\).
For symmetric matrices, \(A\preceq B\) means that
\(B-A\) is positive semidefinite.
Throughout this e-companion, $\mathrm{Id}_s$ denotes the $s\times s$
identity matrix.
Assumption~\ref{ass:finiteSampleCovariance}(\ref{cond:chartwiseSupportEnvelope}) gives
\(\|W_i\|_2\leq R_T\).
Both \(W_iW_i^\top\) and \(\E[W_iW_i^\top]\) lie between \(0\) and
\(R_T^2\mathrm{Id}_{d_\xi}\), so
\(-R_T^2\mathrm{Id}_{d_\xi}\preceq Z_i\preceq R_T^2\mathrm{Id}_{d_\xi}\)
and \(\|Z_i\|_{\rm op}\leq R_T^2\).
Expanding the square and using conditional measurability of
\(\E[W_iW_i^\top]\) gives
\[
 \begin{aligned}
 \E[Z_i^2]
 &=\E\!\left[(W_iW_i^\top)^2
 -W_iW_i^\top\E[W_iW_i^\top]
 -\E[W_iW_i^\top]W_iW_i^\top
 +\{\E[W_iW_i^\top]\}^2\right]\\
 &=\E[\|W_i\|_2^2W_iW_i^\top]-\{\E[W_iW_i^\top]\}^2,\\
 0\preceq\E[Z_i^2]
 &\preceq\E[\|W_i\|_2^2W_iW_i^\top]
 \preceq R_T^2\E[W_iW_i^\top]\preceq R_T^4\mathrm{Id}_{d_\xi}.
 \end{aligned}
\]
Hence \(\|\sum_{i\in I_b}\E[Z_i^2]\|_{\rm op}\leq |I_b|R_T^4\).
Assumption~\ref{ass:finiteSampleCovariance}(\ref{cond:chartwiseIndependence}) gives conditional independence;
centering gives \(\E[Z_i]=0\).
Lemma~\ref{lem:ecMatrixBernstein}, conditional on
\(\mathcal F_{x,\rm tr}\), with \(L=R_T^2\), \(v=|I_b|R_T^4\), and
\(u=t\), gives
\begin{equation}
 \left\|\frac{1}{|I_b|}\sum_{i\in I_b}
 \{W_iW_i^\top-\E[W_iW_i^\top]\}\right\|_{\rm op}
 \leq R_T^2\sqrt{\frac{2t}{|I_b|}}+\frac{2R_T^2t}{3|I_b|},
 \label{eq:ecSecondMomentConcentration}
\end{equation}
with conditional failure probability at most
\(2d_\xi\exp(-t)\leq\delta/2\).

\emph{Step 4. Bounding the sample mean error.}
Let \(E_i:=W_i-\mu_T(x_i)\).
Lemma~\ref{lem:ecMatrixBernstein} bounds sums of symmetric matrices,
so we express the sample mean error as
\begin{equation}
 \|\widehat\mu_b-\bar\mu_{T,b}\|_2
 =\frac{1}{|I_b|}\left\|\sum_{i\in I_b}E_i\right\|_2
 =\frac{1}{|I_b|}\left\|\sum_{i\in I_b}
 \begin{pmatrix}0&E_i^\top\\E_i&0\end{pmatrix}\right\|_{\rm op}.
 \label{eq:ecMeanErrorDilation}
\end{equation}
Assumption~\ref{ass:finiteSampleCovariance}(\ref{cond:chartwiseIndependence}) gives conditional independence
of the block matrices. Since \(\E[W_i]=\mu_T(x_i)\), we have
\(\E[E_i]=0\), so the block matrices have conditional mean zero.
The bounds \(\|\mu_T(x_i)\|_2\leq\E[\|W_i\|_2]\leq R_T\)
and \(\|W_i\|_2\leq R_T\) give \(\|E_i\|_2\leq2R_T\).
Expanding the squared norm gives
\[
 \E[\|E_i\|_2^2]
 =\E[\|W_i\|_2^2]-2\mu_T(x_i)^\top\E[W_i]+\|\mu_T(x_i)\|_2^2
 =\E[\|W_i\|_2^2]-\|\mu_T(x_i)\|_2^2\leq R_T^2.
\]
Since \(E_iE_i^\top\preceq\|E_i\|_2^2\mathrm{Id}_{d_\xi}\), the
symmetric block matrices satisfy
\[
 \begin{aligned}
 \begin{pmatrix}0&E_i^\top\\E_i&0\end{pmatrix}^{\!2}
 &=\begin{pmatrix}\|E_i\|_2^2&0\\0&E_iE_i^\top\end{pmatrix},\\
 \sum_{i\in I_b}\E\!\left[
 \begin{pmatrix}0&E_i^\top\\E_i&0\end{pmatrix}^{\!2}\right]
 &=\begin{pmatrix}
 \displaystyle\sum_{i\in I_b}\E[\|E_i\|_2^2]&0\\
 0&\displaystyle\sum_{i\in I_b}\E[E_iE_i^\top]
 \end{pmatrix}
 \preceq |I_b|R_T^2\mathrm{Id}_{d_\xi+1}.
 \end{aligned}
\]
We also have
\(\left\|\begin{pmatrix}0&E_i^\top\\E_i&0\end{pmatrix}\right\|_{\rm op}
=\|E_i\|_2\leq2R_T\).
Applying Lemma~\ref{lem:ecMatrixBernstein} to the block matrices with
\(L=2R_T\), \(v=|I_b|R_T^2\), and \(u=t\), and using
Eq.~\eqref{eq:ecMeanErrorDilation} gives
\begin{equation}
 \|\widehat\mu_b-\bar\mu_{T,b}\|_2
 \leq R_T\sqrt{\frac{2t}{|I_b|}}+\frac{4R_Tt}{3|I_b|},
 \label{eq:ecMeanConcentration}
\end{equation}
with conditional failure probability at most
\(2(d_\xi+1)\exp(-t)=\delta/2\).

Both \(\|\widehat\mu_b\|_2\) and \(\|\bar\mu_{T,b}\|_2\) are bounded by
\(R_T\), so on the event in Eq.~\eqref{eq:ecMeanConcentration},
\[
 \begin{aligned}
 \widehat\mu_b\widehat\mu_b^\top-\bar\mu_{T,b}\bar\mu_{T,b}^\top
 &=(\widehat\mu_b-\bar\mu_{T,b})\widehat\mu_b^\top
  +\bar\mu_{T,b}(\widehat\mu_b-\bar\mu_{T,b})^\top,\\
 \|\widehat\mu_b\widehat\mu_b^\top-\bar\mu_{T,b}\bar\mu_{T,b}^\top\|_{\rm op}
 &\leq2R_T\|\widehat\mu_b-\bar\mu_{T,b}\|_2
 \leq2R_T^2\sqrt{\frac{2t}{|I_b|}}+\frac{8R_T^2t}{3|I_b|}.
 \end{aligned}
\]
\emph{Step 5. Combining the error bounds.}
By a union bound, Eqs.~\eqref{eq:ecSecondMomentConcentration} and
\eqref{eq:ecMeanConcentration} hold simultaneously with conditional probability
at least \(1-\delta\).
Substituting Eq.~\eqref{eq:ecSecondMomentConcentration} and the
preceding outer-product bound into
Eq.~\eqref{eq:ecRawCovarianceDecomposition}, and using
Eq.~\eqref{eq:chartSamplingRadius}, gives
\begin{equation}
 \|\widehat\Sigma_W-\Sigma_T(\bar x_b)\|_{\rm op}
 \leq\varepsilon_b^{\rm loc}
 +3\sqrt{2}R_T^2\sqrt{\frac{t}{|I_b|}}+\frac{10R_T^2t}{3|I_b|}
 =\varepsilon_b^{\rm loc}+\varepsilon_b^{\rm samp}(\delta).
 \label{eq:ecRawCovarianceConcentration}
\end{equation}
Eqs.~\eqref{eq:ecCrossCovarianceReduction},
\eqref{eq:ecCovarianceScaling}, and
\eqref{eq:ecRawCovarianceConcentration} now yield
\[
 \begin{aligned}
 \varepsilon_b^C
 =\|\widehat C_b-C_b^\star\|_{\rm op}
 &\leq\overline H_b\|\widehat D_b^{-1}\|_{\rm op}^2
 \|\widehat\Sigma_W-\Sigma_T(\bar x_b)\|_{\rm op}
 +2R_T\|\widehat D_b^{-1}\|_{\rm op}\varepsilon_b^{\rm aff}\\
 &\leq\overline H_b\|\widehat D_b^{-1}\|_{\rm op}^2
 \{\varepsilon_b^{\rm loc}+\varepsilon_b^{\rm samp}(\delta)\}
 +2R_T\|\widehat D_b^{-1}\|_{\rm op}\varepsilon_b^{\rm aff}
 =\overline\varepsilon_b^C(\delta).
 \end{aligned}
\]
\hfill\Halmos
\endproof

\subsection{Covariance localization}
\label{ec:spectralRefinements}
\label{ec:proofFiniteSampleCovarianceLocalization}

For each chart, define
\(\widehat J_b:=|I_b|^{-1}\sum_{i\in I_b}
\|x_i-\bar x_b\|_x^2\).

\begin{lemma}
\label{cor:finiteSampleCovarianceLocalization}
Suppose that there exist finite nonnegative constants $L_{\mu,b}$ and
$L_{\Sigma,b}$ such that, for every context assigned to chart $b$,
\(\bigl\|\mu_T(x)-\mu_T(\bar x_b)\bigr\|_2
\leq L_{\mu,b}\|x-\bar x_b\|_x\) and
\(\bigl\|\Sigma_T(x)-\Sigma_T(\bar x_b)\bigr\|_{\rm op}
\leq L_{\Sigma,b}\|x-\bar x_b\|_x\).
Then
\begin{equation}
 \varepsilon_b^{\rm loc}
 \leq
 L_{\Sigma,b}\sqrt{\widehat J_b}
 +
 L_{\mu,b}^2\widehat J_b.
 \label{eq:finiteSampleCovarianceLocalization}
\end{equation}
\end{lemma}
\proof{Proof.}

Let \(M_b:=|I_b|^{-1}\sum_{i\in I_b}
\{\mu_T(x_i)-\bar\mu_{T,b}\}
\{\mu_T(x_i)-\bar\mu_{T,b}\}^{\top}\).
Eq.~\eqref{eq:chartMixtureBias} gives
\[
 \begin{aligned}
 \varepsilon_b^{\rm loc}
 &=\left\|
 \frac{1}{|I_b|}\sum_{i\in I_b}
 \{\Sigma_T(x_i)-\Sigma_T(\bar x_b)\}+M_b
 \right\|_{\rm op}\\
 &\leq
 \left\|
 \frac{1}{|I_b|}\sum_{i\in I_b}
 \{\Sigma_T(x_i)-\Sigma_T(\bar x_b)\}
 \right\|_{\rm op}
 +\|M_b\|_{\rm op}\\
 &\leq
 \frac{L_{\Sigma,b}}{|I_b|}
 \sum_{i\in I_b}\|x_i-\bar x_b\|_x+\|M_b\|_{\rm op}\\
 &\leq
 L_{\Sigma,b}
 \left(\frac{1}{|I_b|}\sum_{i\in I_b}
 \|x_i-\bar x_b\|_x^2\right)^{1/2}
 +\|M_b\|_{\rm op}\\
 &=L_{\Sigma,b}\sqrt{\widehat J_b}+\|M_b\|_{\rm op}.
 \end{aligned}
\]

Since \(\sum_{i\in I_b}\{\mu_T(x_i)-\bar\mu_{T,b}\}=0\),
expanding the squared norms yields
\[
 \begin{aligned}
\sum_{i\in I_b}\|\mu_T(x_i)-\mu_T(\bar x_b)\|_2^2 &=
 \sum_{i\in I_b}\|\mu_T(x_i)-\bar\mu_{T,b}\|_2^2
 +|I_b|\|\bar\mu_{T,b}-\mu_T(\bar x_b)\|_2^2\\
 &\qquad+
 2\left\langle
 \sum_{i\in I_b}\{\mu_T(x_i)-\bar\mu_{T,b}\},
 \bar\mu_{T,b}-\mu_T(\bar x_b)
 \right\rangle\\
 &=
 \sum_{i\in I_b}\|\mu_T(x_i)-\bar\mu_{T,b}\|_2^2
 +|I_b|\|\bar\mu_{T,b}-\mu_T(\bar x_b)\|_2^2.
 \end{aligned}
\]
As an average of outer products, \(M_b\) is positive semidefinite.
Its operator norm is therefore at most its trace $\operatorname{tr}(M_b)$, and we then have
\[
 \begin{aligned}
 \|M_b\|_{\rm op}
 &\leq\operatorname{tr}(M_b)
 =\frac{1}{|I_b|}\sum_{i\in I_b}
 \|\mu_T(x_i)-\bar\mu_{T,b}\|_2^2\\
 &\leq
 \frac{1}{|I_b|}\sum_{i\in I_b}
 \|\mu_T(x_i)-\mu_T(\bar x_b)\|_2^2\\
 &\leq
 \frac{L_{\mu,b}^2}{|I_b|}\sum_{i\in I_b}
 \|x_i-\bar x_b\|_x^2
 =L_{\mu,b}^2\widehat J_b.
 \end{aligned}
\]
Combining the two bounds proves
Eq.~\eqref{eq:finiteSampleCovarianceLocalization}.
\Halmos
\endproof

\subsection{Proof of
Theorem~\ref{thm:atlasCrossCovarianceTradeoff}}
\label{ec:proofAtlasCrossCovarianceTradeoff}

\proof{Proof.}
Let \(t:=\log\{4(d_\xi+1)B/\delta\}\).
The sample size requirement in Theorem~\ref{thm:atlasCrossCovarianceTradeoff}
gives \(|I_b|\geq\max\{2,t\}\) for every chart.
Under Assumption~\ref{ass:finiteSampleCovariance},
Lemma~\ref{prop:finiteSampleCrossCovariance}, with conditional failure
probability at most \(\delta/B\) for each chart, gives
\begin{equation}
 \varepsilon_b^C
 \leq\overline\varepsilon_b^C(\delta/B)
 =\overline H_b\|\widehat D_b^{-1}\|_{\rm op}^2
 \{\varepsilon_b^{\rm loc}+\varepsilon_b^{\rm samp}(\delta/B)\}
 +2R_T\|\widehat D_b^{-1}\|_{\rm op}\varepsilon_b^{\rm aff}.
 \label{eq:ecChartwiseEstimationBound}
\end{equation}
To bound the weighted average over all charts, we need these bounds to
hold simultaneously. A union bound gives
\[
 \mathbb P\left(
 \bigcup_{b=1}^B
 \{\varepsilon_b^C>\overline\varepsilon_b^C(\delta/B)\}
 \,\middle|\,\mathcal F_{x,\rm tr}\right)
 \leq\sum_{b=1}^B
 \mathbb P\left(
 \varepsilon_b^C>\overline\varepsilon_b^C(\delta/B)
 \,\middle|\,\mathcal F_{x,\rm tr}\right)
 \leq\delta.
\]
Thus, Eq.~\eqref{eq:ecChartwiseEstimationBound} holds for every chart
with conditional probability at least \(1-\delta\).
We work on this simultaneous event for the remainder of the proof.

Since \(|I_b|\geq t\), we have \(0<t/|I_b|\leq1\) and
\(t/|I_b|\leq\sqrt{t/|I_b|}\).
Using the definition in Eq.~\eqref{eq:chartSamplingRadius} and
\(3\sqrt{2}+10/3<8\), we obtain
\begin{equation}
 \varepsilon_b^{\rm samp}(\delta/B)
 =R_T^2\left(3\sqrt{2}\sqrt{\frac{t}{|I_b|}}
 +\frac{10t}{3|I_b|}\right)
 \leq\left(3\sqrt{2}+\frac{10}{3}\right)R_T^2\sqrt{\frac{t}{|I_b|}}
 \leq8R_T^2\sqrt{\frac{t}{|I_b|}}.
 \label{eq:ecChartSamplingSimplification}
\end{equation}
Lemma~\ref{cor:finiteSampleCovarianceLocalization}, with the common
constants \(L_\Sigma\) and \(L_\mu\) in Eq.~\eqref{eq:chartwiseMomentBounds}, gives
\(\varepsilon_b^{\rm loc}\leq
L_\Sigma\sqrt{\widehat J_b}+L_\mu^2\widehat J_b\).
The scaling bound in Eq.~\eqref{eq:chartwiseRegularity} gives
\(\overline H_b\|\widehat D_b^{-1}\|_{\rm op}^2\leq A_{\rm sc}\) and
\(2R_T\|\widehat D_b^{-1}\|_{\rm op}\leq A_{\rm sc}\).
Substituting these bounds and
Eq.~\eqref{eq:ecChartSamplingSimplification} into
Eq.~\eqref{eq:ecChartwiseEstimationBound} yields, for every chart,
\begin{equation}
 \varepsilon_b^C\leq A_{\rm sc}
 \left\{L_\Sigma\sqrt{\widehat J_b}+L_\mu^2\widehat J_b
 +\varepsilon_b^{\rm aff}+8R_T^2\sqrt{\frac{t}{|I_b|}}\right\}.
 \label{eq:ecChartErrorBound}
\end{equation}

We now average Eq.~\eqref{eq:ecChartErrorBound} with weights
\(|I_b|/|I_{\rm tr}|\).
The partition and distortion definitions give
\(\sum_{b=1}^B|I_b|/|I_{\rm tr}|=1\) and
\(\sum_{b=1}^B|I_b|\widehat J_b/|I_{\rm tr}|=\widehat J_B\).
The two remaining weighted sums satisfy, by Cauchy--Schwarz,
\[
 \begin{aligned}
 \sum_{b=1}^B\frac{|I_b|}{|I_{\rm tr}|}\sqrt{\widehat J_b}
 &=\sum_{b=1}^B
 \sqrt{\frac{|I_b|}{|I_{\rm tr}|}}
 \sqrt{\frac{|I_b|\widehat J_b}{|I_{\rm tr}|}}
 \leq\left(\sum_{b=1}^B\frac{|I_b|}{|I_{\rm tr}|}\right)^{1/2}
 \left(\sum_{b=1}^B\frac{|I_b|\widehat J_b}{|I_{\rm tr}|}\right)^{1/2}
 =\sqrt{\widehat J_B},\\
 \sum_{b=1}^B\frac{|I_b|}{|I_{\rm tr}|}\sqrt{\frac{t}{|I_b|}}
 &=\sum_{b=1}^B
 \sqrt{\frac{|I_b|}{|I_{\rm tr}|}}
 \sqrt{\frac{t}{|I_{\rm tr}|}}
 \leq\left(\sum_{b=1}^B\frac{|I_b|}{|I_{\rm tr}|}\right)^{1/2}
 \left(\sum_{b=1}^B\frac{t}{|I_{\rm tr}|}\right)^{1/2}
 =\sqrt{\frac{Bt}{|I_{\rm tr}|}}.
 \end{aligned}
\]
Multiplying Eq.~\eqref{eq:ecChartErrorBound} by
\(|I_b|/|I_{\rm tr}|\), summing over \(b\), and using these identities
and inequalities gives
\[
 \begin{aligned}
 \sum_{b=1}^B\frac{|I_b|}{|I_{\rm tr}|}\varepsilon_b^C
 &\leq
 A_{\rm sc}\sum_{b=1}^B\frac{|I_b|}{|I_{\rm tr}|}
 \left\{
 L_\Sigma\sqrt{\widehat J_b}
 +L_\mu^2\widehat J_b
 +\varepsilon_b^{\rm aff}
 +8R_T^2\sqrt{\frac{t}{|I_b|}}
 \right\}\\
 &\leq
 A_{\rm sc}\left\{
 L_\Sigma\sqrt{\widehat J_B}
 +L_\mu^2\widehat J_B
 +\bar\varepsilon_B^{\rm aff}
 +8R_T^2\sqrt{\frac{Bt}{|I_{\rm tr}|}}
 \right\}.
 \end{aligned}
\]
Substituting \(t=\log\{4(d_\xi+1)B/\delta\}\) proves
Eq.~\eqref{eq:atlasCrossCovarianceTradeoff} with conditional probability
at least \(1-\delta\).
\Halmos
\endproof

\subsection{Ideal cross-covariance subspace}
\label{ec:idealCrossCovarianceSubspace}

\begin{lemma}
\label{lem:idealTiltedSpace}
Suppose Assumption~\ref{ass:localSignalNondegeneracy} holds for the
decomposition in Definition~\ref{def:affineGradientDecomposition}. Then
\begin{equation}
 \operatorname{range}(C_b^\star)
 =\operatorname{range}(\Sigma_bV_b^{\rm rel}),
 \qquad
 \operatorname{rank}(C_b^\star)=d_b^\star.
\label{eq:idealSignalRangeRank}
\end{equation}
Moreover, for every matrix
\(V_b^\star\in\RR^{d_\xi\times d_b^\star}\) with orthonormal columns spanning
\(\operatorname{range}(C_b^\star)\),
\((V_b^{\rm rel})^\top V_b^\star\) is nonsingular.
\end{lemma}

\proof{Proof.}

Let \(H:=H_b(z_b^{\rm SAA},\bar x_b)\).
Full column rank of \(H\) gives
\(a^\top H^\top Ha=\|Ha\|_2^2>0\) for every nonzero
\(a\in\RR^{d_b^\star}\).
Thus, \(H^\top H\) is nonsingular, and
\(H(H^\top H)^{-1}\) is a right inverse of \(H^\top\), since
\(H^\top H(H^\top H)^{-1}=\mathrm{Id}_{d_b^\star}\).
The identity \(C_b^\star=\Sigma_bV_b^{\rm rel}H^\top\) from
Eq.~\eqref{eq:idealSignalDefinition} and
\(C_b^\star H(H^\top H)^{-1}=\Sigma_bV_b^{\rm rel}\)
give both range inclusions
\[
 \begin{aligned}
 \operatorname{range}(C_b^\star)
 &=
 \operatorname{range}(\Sigma_bV_b^{\rm rel}H^\top)
 \subseteq\operatorname{range}(\Sigma_bV_b^{\rm rel}),\\
 \operatorname{range}(\Sigma_bV_b^{\rm rel})
 &=
 \operatorname{range}\!\left(
 C_b^\star H(H^\top H)^{-1}
 \right)
 \subseteq\operatorname{range}(C_b^\star).
 \end{aligned}
\]
If \(\Sigma_bV_b^{\rm rel}a=0\), then
\(a^\top(V_b^{\rm rel})^\top\Sigma_bV_b^{\rm rel}a
=(V_b^{\rm rel}a)^\top(\Sigma_bV_b^{\rm rel}a)=0\).
Positive definiteness of
\((V_b^{\rm rel})^\top\Sigma_bV_b^{\rm rel}\) forces \(a=0\).
Hence \(\Sigma_bV_b^{\rm rel}\) has full column rank \(d_b^\star\), and
both matrices have the same range of dimension \(d_b^\star\), proving
Eq.~\eqref{eq:idealSignalRangeRank}.

For any orthonormal basis \(V_b^\star\) of
\(\operatorname{range}(\Sigma_bV_b^{\rm rel})\), each column of
\(V_b^\star\) is a linear combination of the columns of
\(\Sigma_bV_b^{\rm rel}\). Hence we can write
\(V_b^\star=\Sigma_bV_b^{\rm rel}A_0\), where
\(A_0\in\RR^{d_b^\star\times d_b^\star}\).
If \(A_0a=0\), then \(V_b^\star a=0\), and orthonormality gives
\(\|a\|_2^2=\|V_b^\star a\|_2^2=0\).
Hence \(a\) must be the zero vector.
Thus, \(A_0\) is nonsingular.
The factorization
\((V_b^{\rm rel})^\top V_b^\star
=\{(V_b^{\rm rel})^\top\Sigma_bV_b^{\rm rel}\}A_0\)
is therefore a product of two nonsingular matrices, proving the
remaining claim.
\Halmos
\endproof

\subsection{\texorpdfstring{Spectral dimension recovery}{Spectral dimension recovery}}
\label{ec:proofSpectralDimensionRecovery}

\begingroup
To compare the selected dimension $d_b$ in Eq.~\eqref{eq:spectralWidthRule} with the
structural dimension $d_b^\star$, define the ideal retained fraction
\begin{equation}
 \alpha_b^\star(d):=
 \frac{\sum_{j=1}^d\sigma_j(C_b^\star)^2}{\|C_b^\star\|_F^2},
 \qquad d=0,\ldots,d_b^\star,
 \label{eq:idealSignalFraction}
\end{equation}
where an empty sum is zero.

\begin{proposition}[Spectral dimension recovery]
\label{thm:spectralDimensionRecovery}
Fix $\delta\in(0,1)$ and $\alpha_{\rm ret}\in(0,1)$. Suppose Assumption~\ref{ass:localSignalNondegeneracy}
holds almost surely.
Let a nonnegative $\overline\varepsilon_b^C(\delta)$ satisfy
\begin{equation}
 \mathbb P\left(
 \varepsilon_b^C\leq\overline\varepsilon_b^C(\delta)
 \;\middle|\;\mathcal F_{x,\rm tr}
 \right)\geq1-\delta.
 \label{eq:crossCovarianceBoundAbstract}
\end{equation}
Define
\begin{equation}
 \varepsilon_b^{\rm spec}(\delta):=
 \frac{\sqrt{\min(d_\xi,d_z)}\,\overline\varepsilon_b^C(\delta)}
   {\|C_b^\star\|_F}
 \left\{2+
 \frac{\sqrt{\min(d_\xi,d_z)}\,\overline\varepsilon_b^C(\delta)}
   {\|C_b^\star\|_F}\right\}.
\label{eq:spectralErrorRadius}
\end{equation}
Let $\mathcal E_b^{\rm ret}(\delta)$ be the event on which
$\varepsilon_b^{\rm spec}(\delta)<1$ and
\begin{equation}
 \frac{2\varepsilon_b^{\rm spec}(\delta)}
   {1-\varepsilon_b^{\rm spec}(\delta)}
 <\min\{\alpha_{\rm ret}-\alpha_b^\star(d_b^\star-1),
     1-\alpha_{\rm ret}\}.
\label{eq:spectralSeparationCondition}
\end{equation}
Then
$\mathbb P(\mathcal E_b^{\rm ret}(\delta)
 \cap\{d_b\neq d_b^\star\}\mid\mathcal F_{x,\rm tr})\leq\delta$.
\end{proposition}
\endgroup

\proof{Proof.}
By Eq.~\eqref{eq:spectralWidthRule}, \(d_b=d_b^\star\) follows if the empirical
retained fraction is below \(\alpha_{\rm ret}\) for every
\(1\leq d<d_b^\star\) and above \(\alpha_{\rm ret}\) at \(d_b^\star\).

Fix the training contexts and suppose that
\(\varepsilon_b^C\leq\overline\varepsilon_b^C(\delta)\) and
\(\mathcal E_b^{\rm ret}(\delta)\) both hold.
The Frobenius error satisfies
\(\|\widehat C_b-C_b^\star\|_F
\leq\sqrt{\min(d_\xi,d_z)}
\|\widehat C_b-C_b^\star\|_{\rm op}
\leq\sqrt{\min(d_\xi,d_z)}\,\overline\varepsilon_b^C(\delta)\).
Write \(S:=\|C_b^\star\|_F^2>0\),
\(\widehat S:=\|\widehat C_b\|_F^2\),
\(S_b(d):=\sum_{j=1}^d\sigma_j(C_b^\star)^2\), and
\(\widehat S_b(d):=\sum_{j=1}^d\sigma_j(\widehat C_b)^2\)
for \(0\leq d\leq\min(d_\xi,d_z)\), with empty sums equal to zero.
By Lemma~\ref{lem:idealTiltedSpace}, \(C_b^\star\) has rank
\(d_b^\star\), so \(\sigma_j(C_b^\star)=0\) for \(j>d_b^\star\).
Consequently,
\begin{equation}
 S_b(d_b^\star)
 =\sum_{j=1}^{d_b^\star}\sigma_j(C_b^\star)^2
 =\sum_{j=1}^{\min(d_\xi,d_z)}\sigma_j(C_b^\star)^2
 =\|C_b^\star\|_F^2=S.
 \label{eq:ecIdealTotalEnergy}
\end{equation}
The equality with \(\|C_b^\star\|_F^2\) uses the identity that the squared
Frobenius norm equals the sum of all squared singular values.

The difference between the sums of the first \(d\) squared singular values satisfies
\begin{equation}
 \begin{aligned}
 |\widehat S_b(d)-S_b(d)|
 &=\left|\sum_{j=1}^d\{\sigma_j(\widehat C_b)-\sigma_j(C_b^\star)\}\{\sigma_j(\widehat C_b)+\sigma_j(C_b^\star)\}\right|\\
 &\leq\left(\sum_{j=1}^d\{\sigma_j(\widehat C_b)-\sigma_j(C_b^\star)\}^2\right)^{1/2}
   \left(\sum_{j=1}^d\{\sigma_j(\widehat C_b)+\sigma_j(C_b^\star)\}^2\right)^{1/2}\\
 &\leq\|\widehat C_b-C_b^\star\|_F\{\sqrt{\widehat S_b(d)}+\sqrt{S_b(d)}\}
 \leq\|\widehat C_b-C_b^\star\|_F\{\|\widehat C_b\|_F+\|C_b^\star\|_F\} \quad\text{(Mirsky)}\\
 &\leq\|\widehat C_b-C_b^\star\|_F\{2\sqrt S+\|\widehat C_b-C_b^\star\|_F\}\\
 &\leq\sqrt{\min(d_\xi,d_z)}\,\overline\varepsilon_b^C(\delta)
 \{2\sqrt S+\sqrt{\min(d_\xi,d_z)}\,\overline\varepsilon_b^C(\delta)\}
 =S\varepsilon_b^{\rm spec}(\delta).
 \end{aligned}
 \label{eq:ecPartialEnergyBound}
\end{equation}

At \(d=\min(d_\xi,d_z)\), the definitions give \(S_b(d)=S\) and
\(\widehat S_b(d)=\widehat S\). Substituting \(d=\min(d_\xi,d_z)\) into
Eq.~\eqref{eq:ecPartialEnergyBound} therefore yields
\(|\widehat S-S|\leq S\varepsilon_b^{\rm spec}(\delta)\).
Since \(\varepsilon_b^{\rm spec}(\delta)<1\) on
\(\mathcal E_b^{\rm ret}(\delta)\),
\(\widehat S\geq S-|\widehat S-S|
\geq S\{1-\varepsilon_b^{\rm spec}(\delta)\}>0\).
For every integer \(0\leq d\leq\min(d_\xi,d_z)\),
\(0\leq S_b(d)\leq S\), so the retained fractions satisfy
\begin{equation}
 \begin{aligned}
 \left|
 \frac{\widehat S_b(d)}{\widehat S}-\frac{S_b(d)}{S}
 \right|
 &=
 \left|
 \frac{\widehat S_b(d)-S_b(d)}{\widehat S}
 +\frac{S_b(d)(S-\widehat S)}{\widehat S S}
 \right|\\
 &\leq
 \frac{|\widehat S_b(d)-S_b(d)|}{\widehat S}
 +\frac{S_b(d)|\widehat S-S|}{\widehat S S}\\
 &\leq
 \frac{2S\varepsilon_b^{\rm spec}(\delta)}{\widehat S}
 \leq\frac{2\varepsilon_b^{\rm spec}(\delta)}
   {1-\varepsilon_b^{\rm spec}(\delta)}.
 \end{aligned}
 \label{eq:ecCumulativeEnergyRatioBound}
\end{equation}

Eq.~\eqref{eq:ecCumulativeEnergyRatioBound} bounds the difference from above
and below, giving
\[
 \frac{S_b(d)}{S}
 -\frac{2\varepsilon_b^{\rm spec}(\delta)}{1-\varepsilon_b^{\rm spec}(\delta)}
 \leq\frac{\widehat S_b(d)}{\widehat S}
 \leq\frac{S_b(d)}{S}
 +\frac{2\varepsilon_b^{\rm spec}(\delta)}{1-\varepsilon_b^{\rm spec}(\delta)}
 \qquad\text{for every }0\leq d\leq\min(d_\xi,d_z).
\]
For \(1\leq d<d_b^\star\), monotonicity gives
\(S_b(d)/S=\alpha_b^\star(d)\leq\alpha_b^\star(d_b^\star-1)\).
At \(d=d_b^\star\), Eq.~\eqref{eq:ecIdealTotalEnergy} gives
\(S_b(d_b^\star)/S=1\).
On \(\mathcal E_b^{\rm ret}(\delta)\),
Eq.~\eqref{eq:spectralSeparationCondition} makes the error bound smaller than
both \(\alpha_{\rm ret}-\alpha_b^\star(d_b^\star-1)\) and
\(1-\alpha_{\rm ret}\). Hence,
\[
 \begin{aligned}
 \frac{\widehat S_b(d)}{\widehat S}
 &\leq
 \frac{S_b(d)}{S}
 +\frac{2\varepsilon_b^{\rm spec}(\delta)}
    {1-\varepsilon_b^{\rm spec}(\delta)}\\
 &\leq
 \alpha_b^\star(d_b^\star-1)
 +\frac{2\varepsilon_b^{\rm spec}(\delta)}
    {1-\varepsilon_b^{\rm spec}(\delta)}
 <\alpha_{\rm ret},
 \qquad 1\leq d<d_b^\star,\\
 \frac{\widehat S_b(d_b^\star)}{\widehat S}
 &\geq
 \frac{S_b(d_b^\star)}{S}
 -\frac{2\varepsilon_b^{\rm spec}(\delta)}{1-\varepsilon_b^{\rm spec}(\delta)}
 =1-\frac{2\varepsilon_b^{\rm spec}(\delta)}{1-\varepsilon_b^{\rm spec}(\delta)}
 >\alpha_{\rm ret}.
 \end{aligned}
\]
Thus, the rule in Eq.~\eqref{eq:spectralWidthRule} selects
\(d_b=d_b^\star\) on the intersection of the two events. Consequently,
\(\mathcal E_b^{\rm ret}(\delta)\cap\{d_b\ne d_b^\star\}
\subseteq
\{\varepsilon_b^C>\overline\varepsilon_b^C(\delta)\}\).
Eq.~\eqref{eq:crossCovarianceBoundAbstract} therefore gives
\[
 \mathbb P\left(
 \mathcal E_b^{\rm ret}(\delta)\cap\{d_b\ne d_b^\star\}
 \,\middle|\,\mathcal F_{x,\rm tr}\right)
 \leq\mathbb P\left(
 \varepsilon_b^C>\overline\varepsilon_b^C(\delta)
 \,\middle|\,\mathcal F_{x,\rm tr}\right)\leq\delta.\Halmos
\]
\endproof

\subsection{Perturbation of the empirical subspace}
\label{ec:empiricalSubspacePerturbation}

Fix chart \(b\). Whenever \(\varepsilon_b^C<\sigma_b\), define
\(\varepsilon_b^V:=\sqrt{2}\varepsilon_b^C/(\sigma_b-\varepsilon_b^C)\); since
\(\sigma_b>0\) and \(\kappa_b\geq1\), the condition
\(\varepsilon_b^C<\sigma_b/(1+\sqrt{2}\kappa_b)\) is equivalent to
\(\varepsilon_b^C<\sigma_b\) and \(\kappa_b\varepsilon_b^V<1\), with
\(\kappa_b/(1-\kappa_b\varepsilon_b^V)=\kappa_b(\sigma_b-\varepsilon_b^C)/\{\sigma_b-(1+\sqrt{2}\kappa_b)\varepsilon_b^C\}\).
With a slight abuse of notation, we omit the subscript \(b\) from
\(d_b^\star\),
\(V_b^{\rm rel}\), \(C_b^\star\), \(\widehat C_b\), \(V_b^\star\),
\(\widehat V_{b,d^\star}\), \(\sigma_b\), \(\kappa_b\),
\(\varepsilon_b^C\), and \(\varepsilon_b^V\)
in the following lemma and the proof of Lemma~\ref{lem:empiricalResponseLift}.

\begin{lemma}
\label{lem:empiricalResponseCompleteness}
Suppose Assumption~\ref{ass:localSignalNondegeneracy} holds for the
decomposition in Definition~\ref{def:affineGradientDecomposition}.
Let \(\sigma=\sigma_{d^\star}(C^\star)\) denote the smallest nonzero
singular value of \(C^\star\), and let \(\widehat V_{d^\star}\) contain
the leading \(d^\star\) left singular vectors of \(\widehat C\).
If \(\varepsilon^C<\sigma\) and \(\kappa\varepsilon^V<1\), then there exists
an orthogonal matrix \(O\in\RR^{d^\star\times d^\star}\) such that
\(\|\widehat V_{d^\star}O-V^\star\|_{\rm op}\leq\varepsilon^V\).
Moreover, \((V^{\rm rel})^\top\widehat V_{d^\star}\) is nonsingular, and
\begin{equation}
 \bigl\|\bigl((V^{\rm rel})^\top \widehat V_{d^\star}\bigr)^{-1}\bigr\|_{\rm op}
 \leq
 \frac{\kappa}
 {1-\kappa\varepsilon^V}.
 \label{eq:ecEmpiricalGeneratorInverse}
\end{equation}
\end{lemma}

\proof{Proof.}
Let \(A^\star:=(V^{\rm rel})^\top V^\star\).
Lemma~\ref{lem:idealTiltedSpace} gives that \(A^\star\) is nonsingular and
\(\|(A^\star)^{-1}\|_{\rm op}=\kappa\).
For any orthogonal matrix \(O\in\RR^{d^\star\times d^\star}\), let
\(E_V:=(V^{\rm rel})^\top(\widehat V_{d^\star}O-V^\star)\).
The projected matrix after right multiplication by \(O\) factors as
\[
 (V^{\rm rel})^\top\widehat V_{d^\star}O
 =A^\star+E_V
 =A^\star\{\mathrm{Id}_{d^\star}+(A^\star)^{-1}E_V\}.
\]
If \(\|E_V\|_{\rm op}\leq\varepsilon^V\), then
\(\|(A^\star)^{-1}E_V\|_{\rm op}
\leq\kappa\|E_V\|_{\rm op}\leq\kappa\varepsilon^V<1\).
For any \(v\in\RR^{d^\star}\) satisfying
\(\{\mathrm{Id}_{d^\star}+(A^\star)^{-1}E_V\}v=0\),
we have \(v=-(A^\star)^{-1}E_Vv\), so
\(\|v\|_2\leq\kappa\varepsilon^V\|v\|_2\).
Since \(\kappa\varepsilon^V<1\), this forces \(v=0\).
Thus, \(\mathrm{Id}_{d^\star}+(A^\star)^{-1}E_V\) is nonsingular.
Since \(A^\star\) is nonsingular, the above factorization also makes
\((V^{\rm rel})^\top\widehat V_{d^\star}O\) nonsingular.
The Neumann series converges in operator norm and gives
\[
 \begin{aligned}
 \{(V^{\rm rel})^\top\widehat V_{d^\star}O\}^{-1}
 &=\{\mathrm{Id}_{d^\star}+(A^\star)^{-1}E_V\}^{-1}(A^\star)^{-1}
 =\left(\sum_{k=0}^{\infty}\{-(A^\star)^{-1}E_V\}^k\right)(A^\star)^{-1},\\
 \bigl\|\{(V^{\rm rel})^\top\widehat V_{d^\star}O\}^{-1}\bigr\|_{\rm op}
 &\leq\kappa\sum_{k=0}^{\infty}\|(A^\star)^{-1}E_V\|_{\rm op}^k
 =\frac{\kappa}{1-\|(A^\star)^{-1}E_V\|_{\rm op}}
 \leq\frac{\kappa}{1-\kappa\varepsilon^V}.
 \end{aligned}
\]

Under \(\|E_V\|_{\rm op}\leq\varepsilon^V\),
\((V^{\rm rel})^\top\widehat V_{d^\star}
=\{(V^{\rm rel})^\top\widehat V_{d^\star}O\}O^\top\)
is a product of nonsingular matrices and is therefore nonsingular. Hence,
\[
 \begin{aligned}
 \bigl\|\{(V^{\rm rel})^\top\widehat V_{d^\star}\}^{-1}\bigr\|_{\rm op}
 &=\bigl\|O\{(V^{\rm rel})^\top\widehat V_{d^\star}O\}^{-1}\bigr\|_{\rm op}
 =\bigl\|\{(V^{\rm rel})^\top\widehat V_{d^\star}O\}^{-1}\bigr\|_{\rm op}
 \leq\frac{\kappa}{1-\kappa\varepsilon^V}.
 \end{aligned}
\]
The inverse bound in Eq.~\eqref{eq:ecEmpiricalGeneratorInverse} therefore holds
whenever \(\|E_V\|_{\rm op}\leq\varepsilon^V\).

Because \(V^{\rm rel}\) has orthonormal columns,
\[
 \|E_V\|_{\rm op}
 \leq\|(V^{\rm rel})^\top\|_{\rm op}
 \|\widehat V_{d^\star}O-V^\star\|_{\rm op}
 =\|\widehat V_{d^\star}O-V^\star\|_{\rm op}.
\]
The preceding argument establishes nonsingularity and
Eq.~\eqref{eq:ecEmpiricalGeneratorInverse} whenever an orthogonal matrix \(O\)
satisfies \(\|\widehat V_{d^\star}O-V^\star\|_{\rm op}\leq\varepsilon^V\).
It remains to construct such a matrix.

Let \(E_C:=\widehat C-C^\star\) and
\(P_\perp:=\mathrm{Id}_{d_\xi}-V^\star(V^\star)^\top\).
Lemma~\ref{lem:idealTiltedSpace} gives \(\operatorname{rank}(C^\star)=d^\star\)
and \(C^\star=V^\star(V^\star)^\top C^\star\).
Thus \(P_\perp C^\star=0\), since \(P_\perp V^\star=0\);
moreover, \(P_\perp^\top=P_\perp^2=P_\perp\) and
\(\|P_\perp\|_{\rm op}\leq1\).
Let \(\widehat W_{d^\star}\) contain the corresponding \(d^\star\)
right singular vectors of \(\widehat C\), and
write \(\widehat\Lambda_{d^\star}:=
\operatorname{Diag}(\sigma_1(\widehat C),\ldots,\sigma_{d^\star}(\widehat C))\),
where $\operatorname{Diag}(\cdot)$ forms a diagonal matrix from its entries.
Then \(\widehat C\widehat W_{d^\star}
=\widehat V_{d^\star}\widehat\Lambda_{d^\star}\), and Weyl's inequality for singular values gives
\[
 \begin{aligned}
 \sigma_{d^\star}(\widehat C)
 &\geq\sigma_{d^\star}(C^\star)-\|E_C\|_{\rm op}
 =\sigma-\varepsilon^C>0,\\
 \|\widehat\Lambda_{d^\star}^{-1}\|_{\rm op}
 &=\max_{1\leq j\leq d^\star}\frac{1}{\sigma_j(\widehat C)}
 =\frac{1}{\sigma_{d^\star}(\widehat C)}
 \leq\frac{1}{\sigma-\varepsilon^C}.
 \end{aligned}
\]
Projecting the singular-vector equation and using
\(\|\widehat W_{d^\star}\|_{\rm op}=1\) yields
\[
 \begin{aligned}
 P_\perp\widehat V_{d^\star}\widehat\Lambda_{d^\star}
 &=P_\perp(C^\star+E_C)\widehat W_{d^\star}
 =P_\perp E_C\widehat W_{d^\star},\\
 \|P_\perp\widehat V_{d^\star}\|_{\rm op}
 &=\|P_\perp E_C\widehat W_{d^\star}\widehat\Lambda_{d^\star}^{-1}\|_{\rm op}\\
 &\leq\|P_\perp\|_{\rm op}\|E_C\|_{\rm op}
 \|\widehat W_{d^\star}\|_{\rm op}\|\widehat\Lambda_{d^\star}^{-1}\|_{\rm op}
 \leq\frac{\varepsilon^C}{\sigma-\varepsilon^C}.
 \end{aligned}
\]
Let \(G:=(V^\star)^\top\widehat V_{d^\star}\).
For \(j=1,\ldots,d^\star\), the orthonormal columns of \(V^\star\) and \(\widehat V_{d^\star}\) give
\(0\leq\sigma_j(G)\leq\|G\|_{\rm op}\leq1\).
Write \(\theta_j:=\arccos(\sigma_j(G))\) for the principal angles,
ordered as \(0\leq\theta_1\leq\cdots\leq\theta_{d^\star}\leq\pi/2\), and let
\(\theta_{\max}:=\theta_{d^\star}\).
Take the singular value decomposition
\(G=O_1\operatorname{Diag}(\cos\theta_1,\ldots,\cos\theta_{d^\star})O_2^\top\),
where \(O_1,O_2\in\RR^{d^\star\times d^\star}\) are orthogonal.
Using \(P_\perp^\top P_\perp=P_\perp\), the Gram matrix of the projected basis is
\[
 (P_\perp\widehat V_{d^\star})^\top
 (P_\perp\widehat V_{d^\star})
 =\mathrm{Id}_{d^\star}-G^\top G
 =O_2\operatorname{Diag}(\sin^2\theta_1,\ldots,\sin^2\theta_{d^\star})O_2^\top.
\]
Its largest eigenvalue is \(\sin^2\theta_{\max}\), so
\(\|P_\perp\widehat V_{d^\star}\|_{\rm op}=\sin\theta_{\max}\).

It remains to convert the principal-angle bound into the required
bound on \(\|\widehat V_{d^\star}O-V^\star\|_{\rm op}\).
Choose \(O:=O_2O_1^\top\), which is orthogonal since
\(O_1\) and \(O_2\) are orthogonal.
Then \(GO=O^\top G^\top
=O_1\operatorname{Diag}(\cos\theta_1,\ldots,\cos\theta_{d^\star})O_1^\top\), and
\[
 \begin{aligned}
 (\widehat V_{d^\star}O-V^\star)^\top(\widehat V_{d^\star}O-V^\star)
 &=2\mathrm{Id}_{d^\star}-O^\top G^\top-GO\\
 &=4O_1\operatorname{Diag}
 \bigl(\sin^2(\theta_1/2),\ldots,\sin^2(\theta_{d^\star}/2)\bigr)O_1^\top.
 \end{aligned}
\]
Taking the square root of its largest eigenvalue and using
\(\cos(\theta_{\max}/2)\geq1/\sqrt{2}\) gives
\[
 \|\widehat V_{d^\star}O-V^\star\|_{\rm op}
 =2\sin(\theta_{\max}/2)
 =\frac{\sin\theta_{\max}}{\cos(\theta_{\max}/2)}
 \leq\sqrt{2}\sin\theta_{\max}
 \leq\frac{\sqrt{2}\varepsilon^C}{\sigma-\varepsilon^C}
 =\varepsilon^V.
\]

Thus, \(\|E_V\|_{\rm op}\leq\varepsilon^V\), so
\((V^{\rm rel})^\top\widehat V_{d^\star}\) is nonsingular and its inverse
satisfies Eq.~\eqref{eq:ecEmpiricalGeneratorInverse}, completing the proof.
\Halmos
\endproof

\subsection{\texorpdfstring{Representation of relevant projections}{Representation of relevant projections}}
\label{ec:proofCoordinateRecovery}

\begingroup
\begin{lemma}[Representation of relevant projections]
\label{lem:empiricalResponseLift}
Suppose Assumption~\ref{ass:localSignalNondegeneracy} holds for the
decomposition in Definition~\ref{def:affineGradientDecomposition}.
If $\varepsilon_b^C<\sigma_b/(1+\sqrt{2}\kappa_b)$, then
$(V_b^{\rm rel})^\top\widehat V_{b,d_b^\star}$ is invertible.
Every $u^{\rm rel}\in\RR^{d_b^\star}$ can therefore be represented in the
column space of $\widehat V_{b,d_b^\star}$ using the scenario coordinate vector
\begin{equation}
 u=\left((V_b^{\rm rel})^\top\widehat V_{b,d_b^\star}\right)^{-1}u^{\rm rel},
 \qquad
 \bigl\|u\bigr\|_2
 \leq
 \frac{\kappa_b(\sigma_b-\varepsilon_b^C)}
 {\sigma_b-(1+\sqrt{2}\kappa_b)\varepsilon_b^C}\bigl\|u^{\rm rel}\bigr\|_2.
\label{eq:empiricalResponseLift}
\end{equation}
\end{lemma}
\endgroup

\proof{Proof.}
Let \(u^{\rm rel}\in\RR^{d^\star}\) be arbitrary.
To represent \(u^{\rm rel}\), it suffices to find \(u\in\RR^{d^\star}\) satisfying
\((V^{\rm rel})^\top\widehat V_{d^\star}u=u^{\rm rel}\) and the coordinate bound
in Eq.~\eqref{eq:empiricalResponseLift}.
Under the stated perturbation conditions,
Lemma~\ref{lem:empiricalResponseCompleteness} gives that
\(A:=(V^{\rm rel})^\top\widehat V_{d^\star}\) is nonsingular
and satisfies Eq.~\eqref{eq:ecEmpiricalGeneratorInverse}.
Set
\begin{equation}
 u:=A^{-1}u^{\rm rel}
 =\{(V^{\rm rel})^\top\widehat V_{d^\star}\}^{-1}u^{\rm rel}.
 \label{eq:ecResponseCoordinateLift}
\end{equation}
Then
\((V^{\rm rel})^\top\widehat V_{d^\star}u
=AA^{-1}u^{\rm rel}=u^{\rm rel}\), and
\[
 \|u\|_2\leq\|A^{-1}\|_{\rm op}\|u^{\rm rel}\|_2
 \leq\frac{\kappa}{1-\kappa\varepsilon^V}\|u^{\rm rel}\|_2
 =\frac{\kappa(\sigma-\varepsilon^C)}{\sigma-(1+\sqrt{2}\kappa)\varepsilon^C}\|u^{\rm rel}\|_2.
\]
Since \(u^{\rm rel}\) was arbitrary, the proof is complete.
\Halmos
\endproof

\subsection{Gradient approximation error}
\label{ec:gradientApproximationError}

Fix chart \(b\). With a slight abuse of notation, we omit the chart index \(b\) from
\(d_b\), \(d_b^\star\), \(V_b^{\rm rel}\), \(\widehat V_{b,d_b^\star}\),
\(r_{y,b}\), \(g_b\), \(\beta_b\), \(H_b\), \(\overline H_b\),
\(r_b\), \(\mu_b^{\rm rel}\),
\(\varepsilon_b^{\rm mean}\), \(\varepsilon_b^{\rm aff}\), \(\varepsilon_b^V\),
\(\varepsilon_b^\circ\), \(\kappa_b\), \(a_b\), \(a_b^\circ\),
\(\Gamma_b\), \(\Gamma_b^\circ\), \(p_b\), \(p_b^\circ\), \(p_{b,k}^\circ\),
\(u_{b,k}^{\rm rel}\), \(u_{b,k}^\circ\), \(\tilde{\xi}_{b,k}^\circ\), and \(\tilde{\boldsymbol{\xi}}_b^\circ\)
in Lemma~\ref{lem:responseToOperator}
and the proof of Theorem~\ref{thm:tiltedResponseCapacity}.

\begin{lemma}
\label{lem:responseToOperator}
Suppose Assumptions~\ref{ass:differentiableScenarioCosts} and
\ref{ass:affineUniformBounds} hold for the decomposition in
Definition~\ref{def:affineGradientDecomposition}.
Suppose $(V^{\rm rel})^\top\widehat V_{d^\star}$ is invertible, and let
$u_1^{\rm rel},\ldots,u_K^{\rm rel}\in\RR^{d^\star}$, with scenario coordinate
vectors $u_k$ given by Eq.~\eqref{eq:ecResponseCoordinateLift}.
Assume $\|u_k\|_2\leq r_y$ for every $k=1,\ldots,K$.
Set $\tilde{\xi}_k:=g(\widehat V_{d^\star}u_k)$. For any probability model $p:\cX_b\to\Delta_K$, define
$\varepsilon^{\rm mean}(p)
:=\sup_{x\in\cX_b}
\bigl\|\mu^{\rm rel}(x)-\sum_{k=1}^Kp_k(x)u_k^{\rm rel}\bigr\|_2$
and let
$R_p(z;x):=\sum_{k=1}^Kp_k(x)\ell(z,\tilde{\xi}_k;x)
 +\eta_{\rm reg}\|z-z_0(x)\|_2^2/2$
be the induced regularized objective.
Then
\begin{equation}
 \sup_{x\in\cX_b,\,z\in\cZ_b}
 \bigl\|\nabla_z R(z;P_x,x)-\nabla_z R_p(z;x)\bigr\|_2
 \leq
 \overline H\varepsilon^{\rm mean}(p)
 +2\varepsilon^{\rm aff}.
 \label{eq:ecResponseOperatorBound}
\end{equation}
\end{lemma}

\proof{Proof.}
Fix \(z\in\cZ_b\) and \(x\in\cX_b\).
Assumption~\ref{ass:differentiableScenarioCosts} permits differentiation
under the conditional expectation in Eq.~\eqref{eq:populationObjectives}.
Using \(\ell=f+Q\) and \(\phi=\nabla_zQ\), we obtain
\[
 \begin{aligned}
 \nabla_zR(z;P_x,x)
 &=\E_{P_x}[\nabla_z\ell(z,\xi;x)]
  +\eta_{\rm reg}(z-z_0(x))\\
 &=\nabla_zf(z;x)+\E_{P_x}[\phi(z,\xi;x)]
  +\eta_{\rm reg}(z-z_0(x)).
 \end{aligned}
\]
The weights \(p_k(x)\) do not depend on \(z\), so the definition of \(R_p\) gives
\[
 \begin{aligned}
 \nabla_zR_p(z;x)
 &=\sum_{k=1}^Kp_k(x)\nabla_z\ell(z,\tilde{\xi}_k;x)
  +\eta_{\rm reg}(z-z_0(x))\\
 &=\left(\sum_{k=1}^Kp_k(x)\right)\nabla_zf(z;x)
  +\sum_{k=1}^Kp_k(x)\phi(z,\tilde{\xi}_k;x)
  +\eta_{\rm reg}(z-z_0(x))\\
 &=\nabla_zf(z;x)+\sum_{k=1}^Kp_k(x)\phi(z,\tilde{\xi}_k;x)
  +\eta_{\rm reg}(z-z_0(x)).
 \end{aligned}
\]
Subtracting the two expressions cancels \(\nabla_zf(z;x)\) and
\(\eta_{\rm reg}(z-z_0(x))\), yielding
\[
 \nabla_zR(z;P_x,x)-\nabla_zR_p(z;x)
 =\E_{P_x}[\phi(z,\xi;x)]
 -\sum_{k=1}^Kp_k(x)\phi(z,\tilde{\xi}_k;x).
\]
Setting $y=T_b(\xi)$ in
Eq.~\eqref{eq:targetAffineResponse}, taking conditional expectations, and using the definition of
\(\mu^{\rm rel}(x)\) gives
\begin{equation}
 \begin{aligned}
 \E_{P_x}[\phi(z,\xi;x)]
 &=\beta(z,x)+H(z,x)\mu^{\rm rel}(x)
 +\E_{P_x}[r(z,T_b(\xi),x)].
 \end{aligned}
 \label{eq:ecPopulationAffineMean}
\end{equation}
The feasible scenario coordinates $u_k$ satisfy
\((V^{\rm rel})^\top\widehat V_{d^\star}u_k=u_k^{\rm rel}\) by
Eq.~\eqref{eq:ecResponseCoordinateLift}.
Thus, Eq.~\eqref{eq:targetAffineResponse} and
\(\sum_{k=1}^Kp_k(x)=1\) give
\begin{equation}
 \begin{aligned}
 \sum_{k=1}^Kp_k(x)\phi(z,\tilde{\xi}_k;x)
 &=\beta(z,x)+H(z,x)\sum_{k=1}^Kp_k(x)u_k^{\rm rel}
 +\sum_{k=1}^Kp_k(x)
 r(z,\widehat V_{d^\star}u_k,x).
 \end{aligned}
 \label{eq:ecCompressedAffineMean}
\end{equation}

Substituting the two affine mean expressions therefore yields
\[
 \begin{aligned}
 &\nabla_zR(z;P_x,x)-\nabla_zR_p(z;x)\\
 &\quad=
 \E_{P_x}[\phi(z,\xi;x)]
 -\sum_{k=1}^Kp_k(x)\phi(z,\tilde{\xi}_k;x)\\
 &\quad=
 H(z,x)\left\{
 \mu^{\rm rel}(x)-\sum_{k=1}^Kp_k(x)u_k^{\rm rel}
 \right\}+\E_{P_x}[r(z,T_b(\xi),x)]
 -\sum_{k=1}^Kp_k(x)
 r(z,\widehat V_{d^\star}u_k,x).
 \end{aligned}
\]
Eq.~\eqref{eq:affineResidualBound} and the definition of
\(\varepsilon^{\rm mean}(p)\) bound the first term on the right by
\(\overline H\varepsilon^{\rm mean}(p)\) in norm.
To obtain Eq.~\eqref{eq:ecResponseOperatorBound}, it remains to bound
each residual average by \(\varepsilon^{\rm aff}\).
Eq.~\eqref{eq:affineResidualBound} bounds the conditional mean of the
residual for observed realizations directly, so
\(\|\E_{P_x}[r(z,T_b(\xi),x)]\|_2\leq\varepsilon^{\rm aff}\). Since
\(\|\widehat V_{d^\star}u_k\|_2=\|u_k\|_2\leq r_y\), the pointwise bound on
\(r\) in Eq.~\eqref{eq:affineResidualBound} and
\(p(x)\in\Delta_K\) give
\begin{equation*}
 \left\|\sum_{k=1}^Kp_k(x)r(z,\widehat V_{d^\star}u_k,x)\right\|_2
 \leq\sum_{k=1}^Kp_k(x)\|r(z,\widehat V_{d^\star}u_k,x)\|_2
 \leq\varepsilon^{\rm aff}\sum_{k=1}^Kp_k(x)
 =\varepsilon^{\rm aff}.
\end{equation*}

Combining the bounds above gives
\[
 \|\nabla_zR(z;P_x,x)-\nabla_zR_p(z;x)\|_2
 \leq\overline H\varepsilon^{\rm mean}(p)+2\varepsilon^{\rm aff}.
\]
Taking the supremum over \(x\in\cX_b\) and \(z\in\cZ_b\) proves
Eq.~\eqref{eq:ecResponseOperatorBound}.
\Halmos
\endproof

\subsection{Decision and objective error bounds}
\label{ec:decisionObjectiveErrorBounds}

\begin{lemma}
\label{lem:monotoneDecisionTransfer}
Suppose \(\cZ\) is nonempty, closed, and convex.
Let $\eta>0$ and $r_{\rm loc}\in(0,\infty]$, and let $R_1$ and $R_2$ be
convex objectives differentiable on a neighborhood of $\cZ$. Let $z_1\in\cZ$ satisfy
$0\in\nabla R_1(z_1)+N_{\cZ}(z_1)$, and suppose $R_1$ is $\eta$-strongly
convex on $\mathcal B_{\cZ}(z_1,r_{\rm loc}):=\{z\in\cZ:\|z-z_1\|_2\leq r_{\rm loc}\}$
(with $\mathcal B_{\cZ}(z_1,\infty):=\cZ$). Then $z_1$ is the unique minimizer
of $R_1$ on $\cZ$. Moreover, let $z_2\in\cZ$ satisfy
$0\in\nabla R_2(z_2)+N_{\cZ}(z_2)$.
If $\sup_{z\in\mathcal B_{\cZ}(z_1,r_{\rm loc})}\bigl\|\nabla R_1(z)-\nabla R_2(z)\bigr\|_2\leq\varepsilon$
for some $\varepsilon\geq0$ satisfying $\varepsilon<\eta r_{\rm loc}$, then
$z_2\in\mathcal B_{\cZ}(z_1,r_{\rm loc})$,
$\bigl\|z_2-z_1\bigr\|_2\leq\varepsilon/\eta$, and
$R_1(z_2)-R_1(z_1)\leq\varepsilon^2/(2\eta)$.
\end{lemma}

\proof{Proof.}
\begingroup
\emph{Step 1. Uniqueness of $z_1$ on $\cZ$.}
Choose $v_1\in N_{\cZ}(z_1)$ satisfying $\nabla R_1(z_1)+v_1=0$.
For any $z'\in\cZ$, convexity of $R_1$ on $\cZ$ and $\langle v_1,z'-z_1\rangle\leq0$ give
$R_1(z')\geq R_1(z_1)+\langle\nabla R_1(z_1),z'-z_1\rangle
=R_1(z_1)-\langle v_1,z'-z_1\rangle\geq R_1(z_1)$, so $z_1$ minimizes $R_1$ on $\cZ$.
Suppose $z'\in\cZ\setminus\{z_1\}$ also satisfies $R_1(z')\leq R_1(z_1)$.
Set $t:=\min\{1,r_{\rm loc}/\|z'-z_1\|_2\}\in(0,1]$ and
$z_t:=tz'+(1-t)z_1\in\cZ$. Then $0<\|z_t-z_1\|_2=t\|z'-z_1\|_2\leq r_{\rm loc}$,
so $z_t\in\mathcal B_{\cZ}(z_1,r_{\rm loc})\setminus\{z_1\}$.
Convexity of $R_1$ on $\cZ$ gives $R_1(z_t)\leq tR_1(z')+(1-t)R_1(z_1)\leq R_1(z_1)$.
Applying $\eta$-strong convexity on $\mathcal B_{\cZ}(z_1,r_{\rm loc})$ at
$\frac12 z_t+\frac12 z_1\in\mathcal B_{\cZ}(z_1,r_{\rm loc})$ yields
\[
 R_1\!\left(\frac{z_t+z_1}{2}\right)
 \leq\frac12 R_1(z_t)+\frac12 R_1(z_1)-\frac{\eta}{8}\|z_t-z_1\|_2^2
 <R_1(z_1),
\]
contradicting optimality of $z_1$ on $\cZ$. Hence $z_1$ is unique.

\emph{Step 2. Localization $z_2\in\mathcal B_{\cZ}(z_1,r_{\rm loc})$ under $\varepsilon<\eta r_{\rm loc}$.}
Choose $v_2\in N_{\cZ}(z_2)$ satisfying $\nabla R_2(z_2)+v_2=0$.
If $r_{\rm loc}=\infty$, then $\mathcal B_{\cZ}(z_1,\infty)=\cZ$ and
$z_2\in\mathcal B_{\cZ}(z_1,r_{\rm loc})$. Suppose $r_{\rm loc}<\infty$ and,
for contradiction, that $\|z_2-z_1\|_2>r_{\rm loc}$.
Set $s:=r_{\rm loc}/\|z_2-z_1\|_2\in(0,1)$ and $z_s:=z_1+s(z_2-z_1)$.
Since $\cZ$ is convex, $z_s\in\cZ$, and $\|z_s-z_1\|_2=s\|z_2-z_1\|_2=r_{\rm loc}$,
so $z_s\in\mathcal B_{\cZ}(z_1,r_{\rm loc})$.
Strong convexity of $R_1$ on $\mathcal B_{\cZ}(z_1,r_{\rm loc})$ between $z_1$ and $z_s$ gives
\[
 \begin{aligned}
 R_1(z_s)&\geq R_1(z_1)+\langle\nabla R_1(z_1),z_s-z_1\rangle+\frac{\eta}{2}\|z_s-z_1\|_2^2,\\
 R_1(z_1)&\geq R_1(z_s)-\langle\nabla R_1(z_s),z_s-z_1\rangle+\frac{\eta}{2}\|z_s-z_1\|_2^2.
 \end{aligned}
\]
Adding these two inequalities yields
$\langle\nabla R_1(z_s)-\nabla R_1(z_1),z_s-z_1\rangle\geq\eta\|z_s-z_1\|_2^2=\eta r_{\rm loc}^2$.
Since $z_s\in\cZ$ and $v_1\in N_{\cZ}(z_1)$, we have
$\langle\nabla R_1(z_1),z_s-z_1\rangle=-\langle v_1,z_s-z_1\rangle\geq0$, so
$\langle\nabla R_1(z_s),z_s-z_1\rangle\geq\eta r_{\rm loc}^2$.
Let $e_s:=\nabla R_1(z_s)-\nabla R_2(z_s)$. Because $z_s\in\mathcal B_{\cZ}(z_1,r_{\rm loc})$,
$\|e_s\|_2\leq\varepsilon$, and Cauchy--Schwarz gives
\[
 \begin{aligned}
 \langle\nabla R_2(z_s),z_s-z_1\rangle
 &=\langle\nabla R_1(z_s),z_s-z_1\rangle-\langle e_s,z_s-z_1\rangle\\
 &\geq\eta r_{\rm loc}^2-\|e_s\|_2\|z_s-z_1\|_2
 \geq(\eta r_{\rm loc}-\varepsilon)r_{\rm loc}>0.
 \end{aligned}
\]
Because $z_2-z_s=(1-s)(z_2-z_1)=\{(1-s)/s\}(z_s-z_1)$ with $(1-s)/s>0$,
multiplying by $(1-s)/s$ gives $\langle\nabla R_2(z_s),z_2-z_s\rangle>0$.
On the other hand, convexity of $R_2$ on $\cZ$ gives
$R_2(z_s)\geq R_2(z_2)+\langle\nabla R_2(z_2),z_s-z_2\rangle$ and
$R_2(z_2)\geq R_2(z_s)+\langle\nabla R_2(z_s),z_2-z_s\rangle$.
Adding these two inequalities and using $\nabla R_2(z_2)=-v_2$ with
$\langle v_2,z_s-z_2\rangle\leq0$ yields
\[
 0\leq\langle\nabla R_2(z_2)-\nabla R_2(z_s),z_2-z_s\rangle
 =\langle v_2,z_s-z_2\rangle-\langle\nabla R_2(z_s),z_2-z_s\rangle
 \leq-\langle\nabla R_2(z_s),z_2-z_s\rangle<0,
\]
a contradiction. Hence $\|z_2-z_1\|_2\leq r_{\rm loc}$, so
$z_2\in\mathcal B_{\cZ}(z_1,r_{\rm loc})$.

\emph{Step 3. Decision and objective bounds.}
\endgroup
Let \(\Delta z:=z_2-z_1\) and
\(e:=\nabla R_1(z_2)-\nabla R_2(z_2)\), so \(\|e\|_2\leq\varepsilon\) since $z_2\in\mathcal B_{\cZ}(z_1,r_{\rm loc})$.
If \(\Delta z=0\), both bounds hold. Otherwise, the distance bound is
equivalent to \(\eta\|\Delta z\|_2^2\leq\varepsilon\|\Delta z\|_2\).
It therefore suffices to establish this inequality.
With \(v_1\in N_{\cZ}(z_1)\) and \(v_2\in N_{\cZ}(z_2)\) from Steps 1 and 2,
the definition of the normal cone gives
\(\langle v_1,\Delta z\rangle\leq0\) and
\(\langle v_2,-\Delta z\rangle\leq0\).
Hence \(\langle v_2,\Delta z\rangle\geq0\) and
\(\langle v_2-v_1,\Delta z\rangle\geq0\).

Strong convexity at the two decisions gives
\[
 \begin{aligned}
 R_1(z_2)
 &\geq R_1(z_1)+\langle\nabla R_1(z_1),\Delta z\rangle
 +\frac{\eta}{2}\|\Delta z\|_2^2,\\
 R_1(z_1)
 &\geq R_1(z_2)-\langle\nabla R_1(z_2),\Delta z\rangle
 +\frac{\eta}{2}\|\Delta z\|_2^2.
 \end{aligned}
\]
Adding the two inequalities gives
\[
 \begin{aligned}
 \eta\|\Delta z\|_2^2
 &\leq\langle\nabla R_1(z_2)-\nabla R_1(z_1),\Delta z\rangle\\
 &\leq\langle\nabla R_1(z_2)-\nabla R_1(z_1)+v_2-v_1,\Delta z\rangle
 &&\text{(}\langle v_2-v_1,\Delta z\rangle\geq0\text{)}\\
 &=\langle\nabla R_1(z_2)-\nabla R_2(z_2),\Delta z\rangle
 &&\text{(}\nabla R_i(z_i)=-v_i,\ i=1,2\text{)}\\
 &=\langle e,\Delta z\rangle\\
 &\leq\|e\|_2\|\Delta z\|_2
\\
 &\leq\varepsilon\|\Delta z\|_2.
 \end{aligned}
\]
Division by \(\eta\|\Delta z\|_2>0\) gives
\(\|z_2-z_1\|_2\leq\varepsilon/\eta\).

It remains to bound \(R_1(z_2)-R_1(z_1)\).
Strong convexity at \(z_2\) gives
\[
 \begin{aligned}
 R_1(z_2)-R_1(z_1)
 &\leq\langle\nabla R_1(z_2),\Delta z\rangle
 -\frac{\eta}{2}\|\Delta z\|_2^2\\
 &=\langle e,\Delta z\rangle-\langle v_2,\Delta z\rangle
 -\frac{\eta}{2}\|\Delta z\|_2^2
 &&\text{(}\nabla R_1(z_2)=e-v_2\text{)}\\
 &\leq\langle e,\Delta z\rangle-\frac{\eta}{2}\|\Delta z\|_2^2
 &&\text{(}\langle v_2,\Delta z\rangle\geq0\text{)}\\
 &\leq\|e\|_2\|\Delta z\|_2-\frac{\eta}{2}\|\Delta z\|_2^2\\
 &\leq\varepsilon\|\Delta z\|_2-\frac{\eta}{2}\|\Delta z\|_2^2\\
 &=\frac{\varepsilon^2}{2\eta}
 -\frac{\eta}{2}
 \left(\|\Delta z\|_2-\frac{\varepsilon}{\eta}\right)^2\\
 &\leq\frac{\varepsilon^2}{2\eta}.\Halmos
 \end{aligned}
\]
\endproof

\subsection{Proof of Theorem~\ref{thm:tiltedResponseCapacity}}
\label{ec:proofTiltedResponseCapacity}

Every bounded subset of $\RR^{d_b^\star}$ lies in a simplex. If
$\|m\|_\infty\leq\overline m_b$ for every $m\in\mathcal M_b$ and some
$\overline m_b>0$, then $\mathcal M_b$ is contained in the simplex
$\{u\in\RR^{d_b^\star}:u_j\geq-\overline m_b\text{ for every }j,\,
\mathbf 1^\top u\leq d_b^\star\overline m_b\}$, whose vertices are
$-\overline m_b\mathbf 1$ and
$-\overline m_b\mathbf 1+2d_b^\star\overline m_be_j$ for
$j=1,\ldots,d_b^\star$, where $e_j$ is the $j$th standard basis vector.
Boundedness is therefore the only requirement on $\mathcal M_b$ for the
vertices used below to exist.

\proof{Proof.}
\emph{Part~(\ref{part:sharedSupportInvertibility}).}
Under Assumption~\ref{ass:localSignalNondegeneracy} and
\(\varepsilon^C<\sigma/(1+\sqrt{2}\kappa)\),
Lemma~\ref{lem:empiricalResponseLift} shows directly that
\((V^{\rm rel})^\top\widehat V_{d^\star}\) is invertible.

\emph{Part~(\ref{part:sharedSupportSufficiency}).}
Let \(K:=d^\star+1\) and fix \(\varepsilon>0\).
\emph{Step 1. Construct the support.}
For each \(k=1,\ldots,K\), set
\(u_k^\circ:=\{(V^{\rm rel})^\top\widehat V_{d^\star}\}^{-1}u_k^{\rm rel}\).
Then
\[
 (V^{\rm rel})^\top\widehat V_{d^\star}u_k^\circ
 =\{(V^{\rm rel})^\top\widehat V_{d^\star}\}
 \{(V^{\rm rel})^\top\widehat V_{d^\star}\}^{-1}u_k^{\rm rel}
 =u_k^{\rm rel}.
\]
Since \(\widehat V_{d^\star}\) has orthonormal columns and
\(\|u_k^\circ\|_2\leq r_y\) by hypothesis,
\(\|\widehat V_{d^\star}u_k^\circ\|_2=\|u_k^\circ\|_2\leq r_y\), so
\(\widehat V_{d^\star}u_k^\circ\in\operatorname{dom}(g)\) by
Assumption~\ref{ass:affineUniformBounds}.
Since \(g:\operatorname{dom}(g)\to\Xi\), we can therefore define
\[
 \tilde{\xi}_k^\circ
 :=g(\widehat V_{d^\star}u_k^\circ)
 =g\!\left(
 \widehat V_{d^\star}
 \{(V^{\rm rel})^\top\widehat V_{d^\star}\}^{-1}u_k^{\rm rel}
 \right)\in\Xi,
 \qquad k=1,\ldots,K.
\]
This defines one scenario tuple
\(\tilde{\boldsymbol{\xi}}^\circ:=(\tilde{\xi}_1^\circ,\ldots,\tilde{\xi}_K^\circ)\in\Xi^K\),
without reference to \(x\) or \(\varepsilon\).
For any probability model \(p:\cX_b\to\Delta_K\), the projection identities give
\[
 \sum_{k=1}^Kp_k(x)(V^{\rm rel})^\top\widehat V_{d^\star}u_k^\circ
 =\sum_{k=1}^Kp_k(x)u_k^{\rm rel},
 \qquad x\in\cX_b.
\]
Thus, the same \(K=d^\star+1\) scenarios satisfy the support requirements of
Lemma~\ref{lem:responseToOperator} for every \(x\in\cX_b\);
the probability model determines their weights and the resulting mean error.

\emph{Step 2. Choose the probability model.} We first show that the infimum in
Eq.~\eqref{eq:conditionalMeanApproximationError} is finite by evaluating one
model in this class. Setting \(a=0\) and \(\Gamma=0\) in
Eq.~\eqref{eq:localProbabilityModel} gives
\[
 p_k(x;0,0)=\frac{\exp(0)}{\sum_{j=1}^K\exp(0)}=\frac1K,
 \qquad k=1,\ldots,K,\quad x\in\cX_b.
\]
To bound the error of this uniform model, fix \(x\in\cX_b\).
Because the convex hull of the vertices contains
\(\mu^{\rm rel}(x)\), there exist
\(\lambda_k(x)\geq0\) with \(\sum_{k=1}^K\lambda_k(x)=1\) and
\(\mu^{\rm rel}(x)=\sum_{k=1}^K\lambda_k(x)u_k^{\rm rel}\).
Thus,
\[
 \begin{aligned}
 \|\mu^{\rm rel}(x)\|_2
 &=\left\|\sum_{k=1}^K\lambda_k(x)u_k^{\rm rel}\right\|_2
 \leq\sum_{k=1}^K\lambda_k(x)\|u_k^{\rm rel}\|_2\\
 &\leq
 \left(\max_{1\leq k\leq K}\|u_k^{\rm rel}\|_2\right)
 \sum_{k=1}^K\lambda_k(x)
 =\max_{1\leq k\leq K}\|u_k^{\rm rel}\|_2,\\
 \left\|
 \mu^{\rm rel}(x)-\frac1K\sum_{k=1}^Ku_k^{\rm rel}
 \right\|_2
 &\leq
 \|\mu^{\rm rel}(x)\|_2
 +\frac1K\sum_{k=1}^K\|u_k^{\rm rel}\|_2
 \leq2\max_{1\leq k\leq K}\|u_k^{\rm rel}\|_2.
 \end{aligned}
\]
The last bound holds for every \(x\in\cX_b\). Since the uniform model is included
in the infimum in Eq.~\eqref{eq:conditionalMeanApproximationError},
\[
 0\leq\varepsilon^{\rm mean}
 \leq\varepsilon^{\rm mean}\bigl(p(\cdot;0,0)\bigr)
 =\sup_{x\in\cX_b}
 \left\|\mu^{\rm rel}(x)-\frac1K\sum_{k=1}^Ku_k^{\rm rel}\right\|_2
 \leq2\max_{1\leq k\leq K}\|u_k^{\rm rel}\|_2<\infty.
\]
Suppose, for contradiction, that every finite pair
\(a\in\RR^K\), \(\Gamma\in\RR^{K\times d_x}\) has mean error greater than
\(\varepsilon^{\rm mean}+\varepsilon\).
Then \(\varepsilon^{\rm mean}+\varepsilon\) is a lower bound strictly larger than
\(\varepsilon^{\rm mean}\), contradicting the definition of the infimum.
Thus, there is a finite pair \((a^\circ,\Gamma^\circ)\) such that,
with \(p^\circ(x):=p(x;a^\circ,\Gamma^\circ)\),
\[
 \varepsilon^{\rm mean}(p^\circ)
 =\sup_{x\in\cX_b}
 \left\|\mu^{\rm rel}(x)-\sum_{k=1}^Kp_k^\circ(x)u_k^{\rm rel}\right\|_2
 \leq\varepsilon^{\rm mean}+\varepsilon.
\]

Under Assumptions~\ref{ass:differentiableScenarioCosts}
and~\ref{ass:affineUniformBounds} and \(\|u_k^\circ\|_2\leq r_y\), applying
Lemma~\ref{lem:responseToOperator} with \(u_k=u_k^\circ\),
\(\tilde{\xi}_k=\tilde{\xi}_k^\circ\), and \(p=p^\circ\) gives, for the resulting distribution \(P_x^\circ\),
\[
 \begin{aligned}
 \sup_{x\in\cX_b,\,z\in\cZ_b}
 \|\nabla_zR(z;P_x,x)-\nabla_zR(z;P_x^\circ,x)\|_2
 &\leq\overline H\varepsilon^{\rm mean}(p^\circ)+2\varepsilon^{\rm aff}\\
 &\leq\overline H(\varepsilon^{\rm mean}+\varepsilon)+2\varepsilon^{\rm aff}
 =\varepsilon^\circ(\varepsilon),
 \end{aligned}
\]
proving part~(\ref{part:sharedSupportSufficiency}).
By Eq.~\eqref{eq:ecEmpiricalGeneratorInverse}, a sufficient condition for
\(\|u_k^\circ\|_2\leq r_y\) for all \(k\) is
\(\kappa(1-\kappa\varepsilon^V)^{-1}\max_{1\leq k\leq K}\|u_k^{\rm rel}\|_2\leq r_y\).

\emph{Part~(\ref{part:sharedSupportMinimality}).}
By hypothesis, every $\mu_b^{\rm rel}(x)$ with $x\in\cX_b$ is a convex
combination of $u_1,\ldots,u_K$, so
$\mathcal M_b\subseteq\operatorname{aff}\{u_1,\ldots,u_K\}$. The affine hull
$\operatorname{aff}\{u_1,\ldots,u_K\}=u_1+\operatorname{span}\{u_2-u_1,\ldots,u_K-u_1\}$
is an affine set of dimension at most $K-1$. Because
$\operatorname{aff}(\mathcal M_b)$ is the smallest affine set containing
$\mathcal M_b$, it is contained in $\operatorname{aff}\{u_1,\ldots,u_K\}$,
so $\dim\operatorname{aff}(\mathcal M_b)\leq K-1$. If $\mathcal M_b$
has nonempty interior in $\RR^{d_b^\star}$, it contains a ball, whose
affine hull is $\RR^{d_b^\star}$; then
$\dim\operatorname{aff}(\mathcal M_b)=d_b^\star$ and $K\geq d_b^\star+1$.
\Halmos
\endproof

\subsection{\texorpdfstring{Local decision quality}{Local decision quality}}
\label{ec:proofLocalDecisionGuarantee}

\begingroup
\begin{lemma}[Local decision quality]
\label{thm:localDecisionGuarantee}
Suppose that Assumptions~\ref{ass:regularizedProblem},
\ref{ass:differentiableScenarioCosts}, and
\ref{ass:operatorContextRegularity} hold. For $x\in\cX_b$, set
$\varepsilon_b(x):=\varepsilon_b^{\rm ctr}+L_B\|x-\bar x_b\|_x$. If
$\varepsilon_b(x)<\eta r_{\rm loc}$, then
\begin{align}
 \bigl\|\widehat z_b(x)-z(P_x;x)\bigr\|_2
 &\leq \frac{\varepsilon_b(x)}{\eta},
 \label{eq:localDecisionDistance}\\
 R(\widehat z_b(x);P_x,x)-R(z(P_x;x);P_x,x)
 &\leq \frac{\varepsilon_b(x)^2}{2\eta}.
 \label{eq:localRegularizedGap}
\end{align}
\end{lemma}
\endgroup

\proof{Proof.}
Fix \(x\in\cX_b\) satisfying \(\varepsilon_b(x)<\eta r_{\rm loc}\).
We apply Lemma~\ref{lem:monotoneDecisionTransfer} with
\(R_1=R(\cdot;P_x,x)\), \(R_2=\widehat R_b(\cdot;x)\),
\(z_1=z(P_x;x)\), and \(z_2=\widehat z_b(x)\).
Assumption~\ref{ass:regularizedProblem} gives a nonempty, closed, convex
feasible set, convex objectives that attain their minima, and
\(\eta\)-strong convexity of \(R(\cdot;P_x,x)\) on \(\mathcal B_{\cZ}(z(P_x;x),r_{\rm loc})\).
Assumption~\ref{ass:differentiableScenarioCosts} gives differentiability.
The minimizers therefore satisfy
\[
 0\in\nabla_zR(z(P_x;x);P_x,x)+N_{\cZ}(z(P_x;x)),
 \qquad
 0\in\nabla_z\widehat R_b(\widehat z_b(x);x)
 +N_{\cZ}(\widehat z_b(x)).
\]
It remains to establish the gradient error bound $\varepsilon_b(x)$ required by
Lemma~\ref{lem:monotoneDecisionTransfer}.
\begingroup
For every \(z\in\mathcal B_{\cZ}(z(P_x;x),r_{\rm loc})\subseteq\cZ_b\),
adding and subtracting \(\nabla_zR(z;P_{\bar x_b},\bar x_b)-\nabla_z\widehat R_b(z;\bar x_b)\)
and applying Eqs.~\eqref{eq:centroidOperatorError} and~\eqref{eq:contextRegularity}
yield
\[
 \begin{aligned}
 &\|\nabla_zR(z;P_x,x)-\nabla_z\widehat R_b(z;x)\|_2\\
 &\leq\|\nabla_zR(z;P_{\bar x_b},\bar x_b)-\nabla_z\widehat R_b(z;\bar x_b)\|_2\\
 &\quad+\bigl\|\{\nabla_zR(z;P_x,x)-\nabla_z\widehat R_b(z;x)\}
 -\{\nabla_zR(z;P_{\bar x_b},\bar x_b)-\nabla_z\widehat R_b(z;\bar x_b)\}\bigr\|_2\\
 &\leq\varepsilon_b^{\rm ctr}+L_B\|x-\bar x_b\|_x=\varepsilon_b(x).
 \end{aligned}
\]
Taking the supremum over \(z\in\mathcal B_{\cZ}(z(P_x;x),r_{\rm loc})\) gives
\[
 \sup_{z\in\mathcal B_{\cZ}(z(P_x;x),r_{\rm loc})}
 \|\nabla_zR(z;P_x,x)-\nabla_z\widehat R_b(z;x)\|_2
 \leq\varepsilon_b(x)<\eta r_{\rm loc}.
\]
\endgroup
Lemma~\ref{lem:monotoneDecisionTransfer}, with error bound
\(\varepsilon_b(x)\), now yields
\[
 \|\widehat z_b(x)-z(P_x;x)\|_2
 \leq\frac{\varepsilon_b(x)}{\eta},
 \qquad
 R(\widehat z_b(x);P_x,x)-R(z(P_x;x);P_x,x)
 \leq\frac{\varepsilon_b(x)^2}{2\eta}.
\]
These are Eqs.~\eqref{eq:localDecisionDistance} and
\eqref{eq:localRegularizedGap}, completing the proof.
\Halmos
\endproof

\subsection{Proof of Theorem~\ref{thm:fixedAtlasRisk}}
\label{ec:fixedAtlasTheory}
\label{ec:unregularizedRisk}

\proof{Proof.}
Write \(\varepsilon(x):=\varepsilon_{\tau_B(x)}(x)
=\varepsilon_{\tau_B(x)}^{\rm ctr}+L_B\|x-\bar x_{\tau_B(x)}\|_x\).
Fix \(x\) with \(\varepsilon(x)<\eta r_{\rm loc}\) and \(\varepsilon(x)<\infty\). Let \(b=\tau_B(x)\).
To bound the expected excess operational risk, we first compare
\(R_0(\widehat z_b(x);x)\) with \(R_0(z;x)\) for an arbitrary
\(z\in\cZ\) satisfying \(\|z-z_0(x)\|_2\leq\operatorname{rad}_{\cZ,0}\).
The regularization penalty is nonnegative, so
\(R_0(\widehat z_b(x);x)\leq R(\widehat z_b(x);P_x,x)\).
Lemma~\ref{thm:localDecisionGuarantee} then gives
\[
 R_0(\widehat z_b(x);x)
 \leq R(z(P_x;x);P_x,x)+\frac{\varepsilon(x)^2}{2\eta}.
\]
Since \(z(P_x;x)\) minimizes \(R(\cdot;P_x,x)\) over \(\cZ\),
\(R(z(P_x;x);P_x,x)\leq R(z;P_x,x)\) for every \(z\in\cZ\). Hence,
\[
 \begin{aligned}
 R_0(\widehat z_b(x);x)
 &\leq R(z;P_x,x)+\frac{\varepsilon(x)^2}{2\eta}\\
 &=R_0(z;x)+\frac{\eta_{\rm reg}}2\|z-z_0(x)\|_2^2
 +\frac{\varepsilon(x)^2}{2\eta}.
 \end{aligned}
\]
Substituting \(\|z-z_0(x)\|_2\leq \operatorname{rad}_{\cZ,0}\) yields
\[
 R_0(\widehat z_b(x);x)
 \leq R_0(z;x)+\frac{\eta_{\rm reg}}2\operatorname{rad}_{\cZ,0}^2
 +\frac{\varepsilon(x)^2}{2\eta}.
\]
The two error terms do not depend on \(z\). Taking the infimum over
\(\{z\in\cZ:\|z-z_0(x)\|_2\leq\operatorname{rad}_{\cZ,0}\}\) therefore gives
\begin{equation}
 R_0(\widehat z_b(x);x)-R_0^\star(x)
 \leq
 \frac{\varepsilon(x)^2}{2\eta}
 +\frac{\eta_{\rm reg}}2\operatorname{rad}_{\cZ,0}^2.
 \label{eq:ecFixedAtlasPointwise}
\end{equation}
Since \(\varepsilon(x)<\eta r_{\rm loc}\) and \(\varepsilon(x)<\infty\) almost surely and
both operational risks are integrable, taking expectations over the deployment context \(x\) gives
\[
 \E[R_0(\widehat z_{\tau_B(x)}(x);x)]-\E[R_0^\star(x)]
 \leq\frac{\E[\varepsilon(x)^2]}{2\eta}
 +\frac{\eta_{\rm reg}}2\operatorname{rad}_{\cZ,0}^2.
\]
It remains to show that
\(\E[\varepsilon(x)^2]\leq
(\varepsilon_B^{\rm ctr}+L_B\sqrt{J_B})^2\).

By definition,
\(\varepsilon(x)=\varepsilon_{\tau_B(x)}^{\rm ctr}
+L_B\|x-\bar x_{\tau_B(x)}\|_x\),
\(\E[(\varepsilon_{\tau_B(x)}^{\rm ctr})^2]
=(\varepsilon_B^{\rm ctr})^2\), and
\(\E[\|x-\bar x_{\tau_B(x)}\|_x^2]=J_B\).
Both quantities are finite by hypothesis. The cross term satisfies
\[
 \begin{aligned}
 &\E\left[
 \varepsilon_{\tau_B(x)}^{\rm ctr}
 \|x-\bar x_{\tau_B(x)}\|_x
 \right]\\
 &\quad\leq
 \left\{\E[(\varepsilon_{\tau_B(x)}^{\rm ctr})^2]\right\}^{1/2}
 \left\{\E[\|x-\bar x_{\tau_B(x)}\|_x^2]\right\}^{1/2}\\
 &\quad=\varepsilon_B^{\rm ctr}\sqrt{J_B}.
 \end{aligned}
\]
Squaring \(\varepsilon(x)\) and taking expectations yields
\[
 \begin{aligned}
 \E[\varepsilon(x)^2]
 &=
 \E\left[
 \{\varepsilon_{\tau_B(x)}^{\rm ctr}
 +L_B\|x-\bar x_{\tau_B(x)}\|_x\}^2
 \right]\\
 &=(\varepsilon_B^{\rm ctr})^2
 +2L_B\E\left[
 \varepsilon_{\tau_B(x)}^{\rm ctr}
 \|x-\bar x_{\tau_B(x)}\|_x
 \right]
 +L_B^2J_B\\
 &\leq
 \left(\varepsilon_B^{\rm ctr}+L_B\sqrt{J_B}\right)^2.
 \end{aligned}
\]
The last inequality uses the cross-term bound above.
Substituting this second-moment bound into the expected risk bound gives
\[
 \E[R_0(\widehat z_{\tau_B(x)}(x);x)]-\E[R_0^\star(x)]
 \leq
 \frac{(\varepsilon_B^{\rm ctr}+L_B\sqrt{J_B})^2}{2\eta}
 +\frac{\eta_{\rm reg}}2\operatorname{rad}_{\cZ,0}^2,
\]
which proves Eq.~\eqref{eq:fixedAtlasRisk}.
\Halmos
\endproof

\begingroup
\subsection{Covering and parameter perturbation bounds}
\label{ec:stageTwoPerturbation}

\begin{lemma}
\label{lem:ecSubspaceCovering}
Fix chart $b$ with integers $d_b\geq1$ and $K_b\geq1$, and let
$P_b:=K_bd_b+(K_b-1)(d_x+1)$.
Suppose $\Theta_b$ is nonempty and every
$\theta_b=(u_{b,1},\ldots,u_{b,K_b},a_b,\Gamma_b)\in\Theta_b$ satisfies
$\mathbf 1^\top a_b=0$, $\mathbf 1^\top\Gamma_b=0$, and
$\|\theta_b\|_2\leq R_\Theta$ for some $R_\Theta>0$.
Then, for every $\epsilon>0$, there exists a finite subset
$\mathcal N_b\subseteq\Theta_b$ such that
$|\mathcal N_b|\leq(1+2R_\Theta/\epsilon)^{P_b}$ and
$\min_{\theta_b^\dagger\in\mathcal N_b}\|\theta_b-\theta_b^\dagger\|_2\leq\epsilon$
for every $\theta_b\in\Theta_b$.
\end{lemma}

\proof{Proof.}
Let $\mathcal S_b:=\{(u_{b,1},\ldots,u_{b,K_b},a_b,\Gamma_b):\mathbf 1^\top a_b=0,\,\mathbf 1^\top\Gamma_b=0\}$.
When $K_b=1$, $\mathbf 1^\top a_b=0$ and $\mathbf 1^\top\Gamma_b=0$ force
$a_b=0$ and $\Gamma_b=0$, so $\Psi_b(\theta_b):=u_{b,1}$ is a linear bijection
from $\mathcal S_b$ onto $\RR^{d_b}=\RR^{P_b}$ satisfying
$\|\Psi_b(\theta_b)\|_2=\|\theta_b\|_2$.
When $K_b\geq2$, the subspace $\{v\in\RR^{K_b}:\mathbf 1^\top v=0\}$ is the
orthogonal complement of $\operatorname{span}(\mathbf 1)$ in $\RR^{K_b}$ and
has dimension $K_b-1$. Choose $E_b\in\RR^{K_b\times(K_b-1)}$ with
$E_b^\top E_b=\mathrm{Id}_{K_b-1}$ and $\mathbf 1^\top E_b=0$.
Every $\theta_b\in\mathcal S_b$ is uniquely represented by
$\check a_b:=E_b^\top a_b\in\RR^{K_b-1}$ and
$\check\Gamma_b:=E_b^\top\Gamma_b\in\RR^{(K_b-1)\times d_x}$ via
$a_b=E_b\check a_b$ and $\Gamma_b=E_b\check\Gamma_b$, with
$\|a_b\|_2^2=\check a_b^\top E_b^\top E_b\check a_b=\|\check a_b\|_2^2$ and
$\|\Gamma_b\|_F^2=\operatorname{tr}(\check\Gamma_b^\top E_b^\top E_b\check\Gamma_b)=\|\check\Gamma_b\|_F^2$.
Stacking $(u_{b,1},\ldots,u_{b,K_b},\check a_b,\operatorname{vec}(\check\Gamma_b))$
therefore defines a linear isometry $\Psi_b:\mathcal S_b\to\RR^{P_b}$
satisfying $\|\Psi_b(\theta_b)-\Psi_b(\theta_b')\|_2=\|\theta_b-\theta_b'\|_2$.

Set $A:=\Psi_b(\Theta_b)\subseteq\{v\in\RR^{P_b}:\|v\|_2\leq R_\Theta\}$,
which is nonempty because $\Theta_b$ is nonempty.
For $v\in\RR^{P_b}$ and $r>0$, write
$\mathbb B_{P_b}(v,r):=\{w\in\RR^{P_b}:\|w-v\|_2<r\}$.
Call $S\subseteq A$ $\epsilon$-separated if $\|v-v'\|_2>\epsilon$ for all
distinct $v,v'\in S$. For any $\epsilon$-separated $S\subseteq A$, the balls
$\{\mathbb B_{P_b}(v,\epsilon/2):v\in S\}$ are pairwise disjoint, since if
$w\in \mathbb B_{P_b}(v,\epsilon/2)\cap \mathbb B_{P_b}(v',\epsilon/2)$,
then $\|v-v'\|_2\leq\|v-w\|_2+\|w-v'\|_2<\epsilon$, forcing $v=v'$.
Moreover, every $w\in \mathbb B_{P_b}(v,\epsilon/2)$ with $v\in S$ satisfies
$\|w\|_2\leq\|v\|_2+\|w-v\|_2<R_\Theta+\epsilon/2$, so
$\bigcup_{v\in S}\mathbb B_{P_b}(v,\epsilon/2)\subseteq \mathbb B_{P_b}(0,R_\Theta+\epsilon/2)$.
Let $\operatorname{vol}_{P_b}$ denote $P_b$-dimensional Lebesgue measure and
$\omega_{P_b}:=\operatorname{vol}_{P_b}(\mathbb B_{P_b}(0,1))>0$.
Translation invariance and homogeneity give
\[
 |S|\left(\frac{\epsilon}{2}\right)^{P_b}\omega_{P_b}
 =\sum_{v\in S}\operatorname{vol}_{P_b}\bigl(\mathbb B_{P_b}(v,\epsilon/2)\bigr)
 \leq\operatorname{vol}_{P_b}\bigl(\mathbb B_{P_b}(0,R_\Theta+\epsilon/2)\bigr)
 =\left(R_\Theta+\frac{\epsilon}{2}\right)^{P_b}\omega_{P_b},
\]
so $|S|\leq(1+2R_\Theta/\epsilon)^{P_b}<\infty$.
Because singletons in $A$ are $\epsilon$-separated and the cardinalities of
$\epsilon$-separated subsets of $A$ are bounded integers, there exists an
$\epsilon$-separated subset $S_{\max}\subseteq A$ of maximum cardinality.
For any $v\in A$, either $v\in S_{\max}$ or $S_{\max}\cup\{v\}$ is not
$\epsilon$-separated, which implies $\|v-v^\dagger\|_2\leq\epsilon$ for some
$v^\dagger\in S_{\max}$. Setting $\mathcal N_b:=\Psi_b^{-1}(S_{\max})\subseteq\Theta_b$
preserves both $|\mathcal N_b|=|S_{\max}|\leq(1+2R_\Theta/\epsilon)^{P_b}$
and $\min_{\theta_b^\dagger\in\mathcal N_b}\|\theta_b-\theta_b^\dagger\|_2\leq\epsilon$
by isometry of $\Psi_b$.
\Halmos
\endproof

\begin{lemma}[{Hoeffding's inequality}]
\label{lem:ecHoeffding}
Let $W_1,\ldots,W_n$ be independent real-valued random variables satisfying
$c\leq W_i\leq c+M_\ell$ almost surely for some
$c\in\RR$ and $M_\ell>0$. Then, for every $t>0$,
\begin{equation}
 \mathbb P\left(\left|\frac{1}{n}\sum_{i=1}^n(W_i-\E[W_i])\right|\geq t\right)
 \leq 2\exp\left(-\frac{2nt^2}{M_\ell^2}\right).
 \label{eq:ecHoeffding}
\end{equation}
\end{lemma}

\proof{Proof.}
Apply \citet[Theorem~2]{hoeffding1963probability} to $W_1,\ldots,W_n$ and to
$-W_1,\ldots,-W_n\in[-c-M_\ell,-c]$, and combine the two one-sided bounds by a
union bound. \Halmos
\endproof

\begin{lemma}[{Parameter Lipschitz continuity}]
\label{lem:ecChartParameterLipschitz}
Suppose Assumptions~\ref{ass:regularizedProblem} and
\ref{ass:differentiableScenarioCosts} hold, and that for chart $b$ every
$\theta_b\in\Theta_b$ satisfies $\widehat V_bu_{b,k}\in\operatorname{dom}(g_b)$
for $k=1,\ldots,K_b$. Suppose there exist finite constants
$\eta_{\rm fit}>0$, $R_x\geq0$, $G_\phi>0$, $L_\phi\geq0$, and $L_\ell>0$
such that, for every $x\in\cX_b$ and $\theta_b,\theta_b'\in\Theta_b$ with
$z=z_b(x;\theta_b)$ and $z'=z_b(x;\theta_b')$,
\begin{enumerate}[label=(\roman*),ref=\roman*]
\item\label{cond:ecCurvature}
$\widehat R_b(\cdot;x,\theta_b)$ is $\eta_{\rm fit}$-strongly convex on $\cZ$;
\item\label{cond:ecGradient}
$\|x-\bar x_b\|_x\leq R_x$, and for each $k=1,\ldots,K_b$, the scenario gradient at $z'$ satisfies\newline
$\|\phi(z',g_b(\widehat V_bu_{b,k});x)\|_2\leq G_\phi$ and
$\|\phi(z',g_b(\widehat V_bu_{b,k});x)-\phi(z',g_b(\widehat V_bu_{b,k}');x)\|_2
\leq L_\phi\|u_{b,k}-u_{b,k}'\|_2$;
\item\label{cond:ecLossLipschitz}
for every $\xi\in\Xi$,
$|\ell(z,\xi;x)-\ell(z',\xi;x)|\leq L_\ell\|z-z'\|_2$.
\end{enumerate}
Then, with
$L_\Theta:=L_\ell\sqrt{L_\phi^2+G_\phi^2(1+R_x^2)}/\eta_{\rm fit}>0$, every
$x\in\cX_b$, $\xi\in\Xi$, and $\theta_b,\theta_b'\in\Theta_b$ satisfy
\begin{equation}
 \bigl|\ell(z_b(x;\theta_b),\xi;x)-\ell(z_b(x;\theta_b'),\xi;x)\bigr|
 \leq L_\Theta\|\theta_b-\theta_b'\|_2.
 \label{eq:ecChartParameterLipschitz}
\end{equation}
\end{lemma}

\proof{Proof.}
Fix $x\in\cX_b$, $\xi\in\Xi$, and $\theta_b,\theta_b'\in\Theta_b$.
Write $z:=z_b(x;\theta_b)$,
$z':=z_b(x;\theta_b')$,
$R_1(\cdot):=\widehat R_b(\cdot;x,\theta_b)$, and
$R_2(\cdot):=\widehat R_b(\cdot;x,\theta_b')$.
By Assumption~\ref{ass:differentiableScenarioCosts} and condition~(\ref{cond:ecCurvature}), $R_1$ and $R_2$
are differentiable and convex on $\cZ$, and $R_1$ is $\eta_{\rm fit}$-strongly
convex on $\cZ$. First-order optimality of $z$ and $z'$ gives
$v_1\in N_{\cZ}(z)$ and $v_2\in N_{\cZ}(z')$ with
$\nabla R_1(z)+v_1=0$ and $\nabla R_2(z')+v_2=0$, so
$\langle\nabla R_1(z),z'-z\rangle=-\langle v_1,z'-z\rangle\geq0$ and
$-\langle\nabla R_2(z'),z'-z\rangle=\langle v_2,z'-z\rangle=-\langle v_2,z-z'\rangle\geq0$.
Strong convexity of $R_1$ between $z$ and $z'$ gives
$R_1(z')-R_1(z)\geq\langle\nabla R_1(z),z'-z\rangle+(\eta_{\rm fit}/2)\|z'-z\|_2^2$
and
$R_1(z)-R_1(z')\geq\langle\nabla R_1(z'),z-z'\rangle+(\eta_{\rm fit}/2)\|z'-z\|_2^2$.
Adding these inequalities and combining with the two normal-cone inequalities
yields
\[
 \begin{aligned}
 \eta_{\rm fit}\|z'-z\|_2^2
 &\leq\langle\nabla R_1(z')-\nabla R_1(z),z'-z\rangle\\
 &\leq\langle\nabla R_1(z')-\nabla R_2(z'),z'-z\rangle
 \leq\|\nabla R_1(z')-\nabla R_2(z')\|_2\|z'-z\|_2,
 \end{aligned}
\]
so $\|z-z'\|_2\leq\eta_{\rm fit}^{-1}\|\nabla R_1(z')-\nabla R_2(z')\|_2$.

We next bound $\|\nabla R_1(z')-\nabla R_2(z')\|_2$.
First, for $s(v):=\softmax(v)$ on $\RR^{K_b}$, its Jacobian
$J_s(v)=\operatorname{Diag}(s(v))-s(v)s(v)^\top$ satisfies, for any
$q=s(v)\in\Delta_{K_b}$ and $h\in\RR^{K_b}$ with $\bar h:=\sum_{j=1}^{K_b}q_jh_j$,
\[
 \|J_s(v)h\|_1
 =\sum_{k=1}^{K_b}q_k|h_k-\bar h|
 \leq\left(\sum_{k=1}^{K_b}q_k\right)^{1/2}
 \left(\sum_{k=1}^{K_b}q_k(h_k-\bar h)^2\right)^{1/2}
 =\left(\sum_{k=1}^{K_b}q_kh_k^2-\bar h^2\right)^{1/2}
 \leq\|h\|_2,
\]
since $0\leq q_k\leq1$. Integrating along $v'+t(v-v')$ for $t\in[0,1]$ gives
$\|s(v)-s(v')\|_1\leq\int_0^1\|J_s(v'+t(v-v'))(v-v')\|_1\,dt\leq\|v-v'\|_2$.
With $\check x_b:=\widehat D_x^{-1}(x-\bar x_b)$, $p:=p_b(x;a_b,\Gamma_b)$,
and $p':=p_b(x;a_b',\Gamma_b')$, condition~(\ref{cond:ecGradient})
gives $\|\check x_b\|_2\leq R_x$, so
\[
 \begin{aligned}
 \|p-p'\|_1
 &\leq\|(a_b-a_b')+(\Gamma_b-\Gamma_b')\check x_b\|_2
 \leq\|a_b-a_b'\|_2+R_x\|\Gamma_b-\Gamma_b'\|_F\\
 &\leq\sqrt{1+R_x^2}\,\bigl(\|a_b-a_b'\|_2^2+\|\Gamma_b-\Gamma_b'\|_F^2\bigr)^{1/2}.
 \end{aligned}
\]
Because $\nabla_z f(z';x)$ and $\eta_{\rm reg}\{z'-z_0(x)\}$ cancel between
$\nabla R_1(z')$ and $\nabla R_2(z')$, adding and subtracting
$\sum_{k=1}^{K_b}p_k'\phi(z',g_b(\widehat V_bu_{b,k});x)$ and applying
condition~(\ref{cond:ecGradient}) give
\[
 \begin{aligned}
 \|\nabla R_1(z')-\nabla R_2(z')\|_2
 &\leq\sum_{k=1}^{K_b}|p_k-p_k'|\,\|\phi(z',g_b(\widehat V_bu_{b,k});x)\|_2\\
 &\quad+\sum_{k=1}^{K_b}p_k'\,\|\phi(z',g_b(\widehat V_bu_{b,k});x)-\phi(z',g_b(\widehat V_bu_{b,k}');x)\|_2\\
 &\leq G_\phi\|p-p'\|_1+L_\phi\sum_{k=1}^{K_b}p_k'\|u_{b,k}-u_{b,k}'\|_2.
 \end{aligned}
\]
Since $p'\in\Delta_{K_b}$ satisfies $0\leq p_k'\leq1$ and $\sum_{k=1}^{K_b}p_k'=1$,
Cauchy--Schwarz gives
$\sum_{k=1}^{K_b}p_k'\|u_{b,k}-u_{b,k}'\|_2\leq(\sum_{k=1}^{K_b}p_k'\|u_{b,k}-u_{b,k}'\|_2^2)^{1/2}\leq(\sum_{k=1}^{K_b}\|u_{b,k}-u_{b,k}'\|_2^2)^{1/2}$.
Applying Cauchy--Schwarz to the two terms above yields
\[
 \|\nabla R_1(z')-\nabla R_2(z')\|_2
 \leq\sqrt{L_\phi^2+G_\phi^2(1+R_x^2)}\,\|\theta_b-\theta_b'\|_2.
\]
Combining this bound with $\|z-z'\|_2\leq\eta_{\rm fit}^{-1}\|\nabla R_1(z')-\nabla R_2(z')\|_2$
and condition~(\ref{cond:ecLossLipschitz}) gives
$|\ell(z,\xi;x)-\ell(z',\xi;x)|\leq L_\ell\|z-z'\|_2\leq L_\Theta\|\theta_b-\theta_b'\|_2$,
proving Eq.~\eqref{eq:ecChartParameterLipschitz}.
\Halmos
\endproof

\subsection{Proof of Theorem~\ref{thm:stageTwoGeneralization}}
\label{ec:proofStageTwoGeneralization}

\proof{Proof.}
\emph{Part~(\ref{part:uniformGeneralization}).}
For $\theta,\theta'\in\Theta_{\rm Atlas}$,
Assumption~\ref{ass:stageTwoRegularity}(\ref{cond:stageTwoLipschitz})
implies that $|\widehat{\mathcal L}_n(\theta)-\widehat{\mathcal L}_n(\theta')|$
and $|\mathcal L(\theta)-\mathcal L(\theta')|$ are both bounded by
$L_\Theta\max_{1\leq b\leq B}\|\theta_b-\theta_b'\|_2$ on the entire sample
space. Because $\Theta_{\rm Atlas}$ is a subset of a finite-dimensional
Euclidean space, it contains a countable dense subset $\Theta_{\rm cnt}$, so
$\sup_{\theta\in\Theta_{\rm Atlas}}|\widehat{\mathcal L}_n(\theta)-\mathcal L(\theta)|
=\sup_{\theta\in\Theta_{\rm cnt}}|\widehat{\mathcal L}_n(\theta)-\mathcal L(\theta)|$
is measurable.

Set $\epsilon:=M_\ell/(2L_\Theta\sqrt{n})>0$ and
$\kappa_n:=1+2R_\Theta/\epsilon=1+4R_\Theta L_\Theta\sqrt{n}/M_\ell>1$.
Under Assumption~\ref{ass:stageTwoRegularity}(\ref{cond:stageTwoDomain}), for
each chart $b=1,\ldots,B$, Lemma~\ref{lem:ecSubspaceCovering} gives an
$\epsilon$-net $\mathcal N_b\subseteq\Theta_b$ of cardinality
$|\mathcal N_b|\leq\kappa_n^{P_b}$, where
$P_b:=K_bd_b+(K_b-1)(d_x+1)$.
Form the product net
$\mathcal N_{\rm Atlas}:=\mathcal N_1\times\cdots\times\mathcal N_B\subseteq\Theta_{\rm Atlas}$,
whose cardinality satisfies
\[
 |\mathcal N_{\rm Atlas}|
 =\prod_{b=1}^B|\mathcal N_b|
 \leq\prod_{b=1}^B\kappa_n^{P_b}
 =\kappa_n^{P_{\rm Atlas}}.
\]
Set
$t:=M_\ell\sqrt{\{P_{\rm Atlas}\log\kappa_n+\log(2/\delta)\}/(2n)}>0$.
For each fixed $\theta^\dagger\in\mathcal N_{\rm Atlas}$,
Assumption~\ref{ass:stageTwoRegularity}(\ref{cond:stageTwoLoss}) gives an
interval $[c_{\theta^\dagger},c_{\theta^\dagger}+M_\ell]$ containing
$\ell(\pi_{\rm Atlas}(x;\theta^\dagger),\xi;x)$ for all $(x,\xi)\in\cX\times\Xi$.
Thus $W_i(\theta^\dagger):=\ell(\pi_{\rm Atlas}(x_i;\theta^\dagger),\xi_i;x_i)$,
$i=1,\ldots,n$, are independent with mean $\mathcal L(\theta^\dagger)$ and
lie in $[c_{\theta^\dagger},c_{\theta^\dagger}+M_\ell]$.
Lemma~\ref{lem:ecHoeffding} and a union bound over
$\mathcal N_{\rm Atlas}$ therefore give
\[
 \mathbb P\left(
 \max_{\theta^\dagger\in\mathcal N_{\rm Atlas}}
 \bigl|\widehat{\mathcal L}_n(\theta^\dagger)-\mathcal L(\theta^\dagger)\bigr|\geq t
 \right)
 \leq 2|\mathcal N_{\rm Atlas}|\exp\left(-\frac{2nt^2}{M_\ell^2}\right)
 \leq 2\kappa_n^{P_{\rm Atlas}}
 \cdot\frac{\delta}{2}\kappa_n^{-P_{\rm Atlas}}
 =\delta.
\]
Work on the complementary event
$\{\max_{\theta^\dagger\in\mathcal N_{\rm Atlas}}|\widehat{\mathcal L}_n(\theta^\dagger)-\mathcal L(\theta^\dagger)|<t\}$,
which has probability at least $1-\delta$.
Given any $\theta=(\theta_1,\ldots,\theta_B)\in\Theta_{\rm Atlas}$, select
$\theta^\dagger=(\theta_1^\dagger,\ldots,\theta_B^\dagger)\in\mathcal N_{\rm Atlas}$
satisfying $\|\theta_b-\theta_b^\dagger\|_2\leq\epsilon$ for every $b=1,\ldots,B$.
For each $(x,\xi)\in\cX\times\Xi$ with $b=\tau_B(x)$,
Assumption~\ref{ass:stageTwoRegularity}(\ref{cond:stageTwoLipschitz}) gives
$|\ell(z_b(x;\theta_b),\xi;x)-\ell(z_b(x;\theta_b^\dagger),\xi;x)|\leq L_\Theta\epsilon=M_\ell/(2\sqrt{n})$.
Averaging over $i=1,\ldots,n$ and taking expectations over $(x,\xi)$ yield
$|\widehat{\mathcal L}_n(\theta)-\widehat{\mathcal L}_n(\theta^\dagger)|\leq M_\ell/(2\sqrt{n})$
and $|\mathcal L(\theta)-\mathcal L(\theta^\dagger)|\leq M_\ell/(2\sqrt{n})$.
The triangle inequality then gives
\[
 \bigl|\widehat{\mathcal L}_n(\theta)-\mathcal L(\theta)\bigr|
 \leq\frac{M_\ell}{\sqrt{n}}+t
 =\varepsilon_n^{\rm est}(\delta),
\]
and taking the supremum over $\theta\in\Theta_{\rm Atlas}$ proves
Eq.~\eqref{eq:stageTwoUniformDeviation}.

\emph{Part~(\ref{part:oracleInequality}).}
On the same event, let
$\widehat\theta=(\widehat\theta_1,\ldots,\widehat\theta_B)\in\Theta_{\rm Atlas}$
satisfy the condition in
Theorem~\ref{thm:stageTwoGeneralization}(\ref{part:oracleInequality}), and fix an arbitrary
$\theta=(\theta_1,\ldots,\theta_B)\in\Theta_{\rm Atlas}$.
For every chart $b$ with $|I_{b,n}|\geq1$, the chartwise condition gives
$|I_{b,n}|^{-1}\sum_{i\in I_{b,n}}\ell(z_b(x_i;\widehat\theta_b),\xi_i;x_i)+\Omega_b(\widehat\theta_b)
\leq|I_{b,n}|^{-1}\sum_{i\in I_{b,n}}\ell(z_b(x_i;\theta_b),\xi_i;x_i)+\Omega_b(\theta_b)+\varepsilon_{\rm opt}$.
Multiplying by $|I_{b,n}|/n$ gives
$n^{-1}\sum_{i\in I_{b,n}}\ell(z_b(x_i;\widehat\theta_b),\xi_i;x_i)+(|I_{b,n}|/n)\Omega_b(\widehat\theta_b)
\leq n^{-1}\sum_{i\in I_{b,n}}\ell(z_b(x_i;\theta_b),\xi_i;x_i)+(|I_{b,n}|/n)\Omega_b(\theta_b)+(|I_{b,n}|/n)\varepsilon_{\rm opt}$,
which holds identically as $0\leq0$ when $I_{b,n}=\emptyset$.
Summing over $b=1,\ldots,B$ and using $\Omega_b(\widehat\theta_b)\geq0$,
$\sum_{b=1}^B(|I_{b,n}|/n)\Omega_b(\theta_b)\leq\{\max_{1\leq b\leq B}\Omega_b(\theta_b)\}\sum_{b=1}^B(|I_{b,n}|/n)=\max_{1\leq b\leq B}\Omega_b(\theta_b)$,
and $\sum_{b=1}^B|I_{b,n}|=n$
gives $\widehat{\mathcal L}_n(\widehat\theta)\leq\widehat{\mathcal L}_n(\theta)+\max_{1\leq b\leq B}\Omega_b(\theta_b)+\varepsilon_{\rm opt}$.
Applying Eq.~\eqref{eq:stageTwoUniformDeviation} at $\widehat\theta$ and $\theta$
then yields
\[
 \begin{aligned}
 \mathcal L(\widehat\theta)
 &\leq\widehat{\mathcal L}_n(\widehat\theta)+\varepsilon_n^{\rm est}(\delta)
 \leq\widehat{\mathcal L}_n(\theta)+\max_{1\leq b\leq B}\Omega_b(\theta_b)+\varepsilon_n^{\rm est}(\delta)+\varepsilon_{\rm opt}\\
 &\leq\mathcal L(\theta)+\max_{1\leq b\leq B}\Omega_b(\theta_b)+2\varepsilon_n^{\rm est}(\delta)+\varepsilon_{\rm opt}.
 \end{aligned}
\]
Taking the infimum over $\theta\in\Theta_{\rm Atlas}$ proves
Eq.~\eqref{eq:stageTwoOracleInequality}.
\Halmos
\endproof

\subsection{Proof of Corollary~\ref{cor:endToEndRisk}}
\label{ec:proofEndToEndRisk}

\proof{Proof.}
Writing $p_b^\circ(x)=p_b(x;\tilde a_b,\tilde\Gamma_b)$ in
Theorem~\ref{thm:tiltedResponseCapacity}(\ref{part:sharedSupportSufficiency}),
set $a_b^\circ:=\tilde a_b-\mathbf 1(\mathbf 1^\top\tilde a_b)/K_b$ and
$\Gamma_b^\circ:=\tilde\Gamma_b-\mathbf 1(\mathbf 1^\top\tilde\Gamma_b)/K_b$.
Then $\mathbf 1^\top a_b^\circ=0$ and $\mathbf 1^\top\Gamma_b^\circ=0$, and
subtracting a common value from all logits at each $x$ preserves the
probabilities in Eq.~\eqref{eq:localProbabilityModel}.
Because $d_b=d_b^\star$, the basis $\widehat V_b$ equals
$\widehat V_{b,d_b^\star}$, and the objective in
Eq.~\eqref{eq:localChartProblem} at $\theta_b^\circ$ is
$R_b^\circ(\cdot;x):=R(\cdot;P_{b,x}^\circ,x)$. Set
$\theta^\circ:=(\theta_1^\circ,\ldots,\theta_B^\circ)$, which belongs to
$\Theta_{\rm Atlas}$ by hypothesis.

Fix $x\in\cX$ and let $b=\tau_B(x)$. Write
$z_b^\circ(x):=z_b(x;\theta_b^\circ)$ for the minimizer of
$R_b^\circ(\cdot;x)$. By
Theorem~\ref{thm:tiltedResponseCapacity}(\ref{part:sharedSupportSufficiency}),
$\sup_{z\in\cZ_b}\|\nabla_z R(z;P_x,x)-\nabla_z R_b^\circ(z;x)\|_2
\leq\varepsilon_b^\circ(\varepsilon)$, and $\cZ_b$ contains
$\mathcal B_{\cZ}(z(P_x;x),r_{\rm loc})$.
Assumptions~\ref{ass:regularizedProblem} and~\ref{ass:differentiableScenarioCosts}
give differentiable convex objectives and local strong convexity of
$R(\cdot;P_x,x)$ on this ball, and the minimizers satisfy the first-order
optimality conditions in Lemma~\ref{lem:monotoneDecisionTransfer}. Since
$\varepsilon_b^\circ(\varepsilon)<\eta r_{\rm loc}$, that lemma gives
$R(z_b^\circ(x);P_x,x)-R(z(P_x;x);P_x,x)\leq
\{\varepsilon_b^\circ(\varepsilon)\}^2/(2\eta)$.

The regularization penalty is nonnegative, so $R_0(z_b^\circ(x);x)\leq R(z_b^\circ(x);P_x,x)$.
By the preceding bound and optimality of $z(P_x;x)$, for every $z\in\cZ$
with $\|z-z_0(x)\|_2\leq\operatorname{rad}_{\cZ,0}$,
\[
 \begin{aligned}
 R_0(z_b^\circ(x);x)
 &\leq R(z(P_x;x);P_x,x)+\frac{\{\varepsilon_b^\circ(\varepsilon)\}^2}{2\eta}
 \leq R(z;P_x,x)+\frac{\{\varepsilon_b^\circ(\varepsilon)\}^2}{2\eta}\\
 &=R_0(z;x)+\frac{\eta_{\rm reg}}{2}\|z-z_0(x)\|_2^2
 +\frac{\{\varepsilon_b^\circ(\varepsilon)\}^2}{2\eta}
 \leq R_0(z;x)+\frac{\eta_{\rm reg}\operatorname{rad}_{\cZ,0}^2}{2}
 +\frac{\{\varepsilon_b^\circ(\varepsilon)\}^2}{2\eta}.
 \end{aligned}
\]
Taking the infimum over these $z$ and using the definition of
$\operatorname{rad}_{\cZ,0}$ gives
\[
 R_0(\pi_{\rm Atlas}(x;\theta^\circ);x)-R_0^\star(x)
 \leq
 \max_{1\leq b\leq B}\frac{\{\varepsilon_b^\circ(\varepsilon)\}^2}{2\eta}
 +\frac{\eta_{\rm reg}\operatorname{rad}_{\cZ,0}^2}{2}.
\]
The loss $\ell(\pi_{\rm Atlas}(x;\theta^\circ),\xi;x)$ takes values in a bounded
interval of length $M_\ell$ by
Assumption~\ref{ass:stageTwoRegularity}(\ref{cond:stageTwoLoss}), and
conditioning on $x$ gives
$\mathcal L(\theta^\circ)=\E[R_0(\pi_{\rm Atlas}(x;\theta^\circ);x)]$.
Taking expectations in the previous display therefore yields
$\mathcal L(\theta^\circ)-\E[R_0^\star(x)]
\leq\max_b\{\varepsilon_b^\circ(\varepsilon)\}^2/(2\eta)
+\eta_{\rm reg}\operatorname{rad}_{\cZ,0}^2/2$.
On the event of Theorem~\ref{thm:stageTwoGeneralization},
Eq.~\eqref{eq:stageTwoOracleInequality} with $\theta=\theta^\circ$ gives
$\mathcal L(\widehat\theta)\leq\mathcal L(\theta^\circ)
+\max_{1\leq b\leq B}\Omega_b(\theta_b^\circ)+2\varepsilon_n^{\rm est}(\delta)
+\varepsilon_{\rm opt}$. Combining the last two inequalities proves
Eq.~\eqref{eq:endToEndRisk}.
\Halmos
\endproof
\endgroup

\section[Models and Experimental Settings]{Models and Experimental Settings}
\label{ec:testProblems}

In this section, we give data generation details and
application-specific parameter choices for the computational study in
Section~\ref{sec:computation}.

\subsection{Smooth capacity allocation}
\label{ec:capacityProblem}

We construct $A$ by drawing its entries independently from a gamma
distribution with shape \(0.65\) and unit scale, then normalizing each row.
We set $C=5d_z$ and $\tau=0.35$. For
$j=1,\ldots,d_\xi$, we set
$c_j^u=3+0.4\{1+\sin(0.31(j-1))\}$ and
$c_j^o=0.7+0.15\{1+\cos(0.23(j-1))\}$.

We generate four demand states indexed by
\(v_r\in\{(-1,-1),(-1,1),(1,-1),(1,1)\}\).
We draw a \(d_z\times2\) matrix with independent standard normal entries,
subtract each column mean, and let $V_{\rm res}$ contain the orthonormal
columns of its QR factorization. We fix this exogenous resource basis before
fitting any method. Define
\(z_{\rm req}^{(r)}=5d_z\,\operatorname{softmax}(V_{\rm res}v_r)\).
The reference scenario set contains \(K_{\rm ref}/4\) scenarios
from each state, and \(r(s)\) denotes the state of scenario \(s\).
With independent $E_{sj}\sim\mathcal N(0,1)$, the scenario components are
\begin{equation}
 \tilde{\xi}_j^{(s)}
 =
 (A z_{\rm req}^{(r(s))})_j
 \exp\!\left\{0.035E_{sj}-\tfrac12(0.035)^2\right\}.
 \label{eq:capacityCatalog}
\end{equation}
The support size study in Section~\ref{subsubsec:supportSizeStudy} uses 500 scenarios.

Contexts are independent draws from $\mathcal N(0,\mathrm{Id}_8)$. Let
\(V_{\rm ctx}\in\RR^{8\times2}\) contain the orthonormal columns of the QR
factorization of an independent standard Gaussian matrix. The conditional
probabilities are uniform within each state, with state probabilities given
by an affine softmax
\begin{equation}
 p_s^\star(x)
 =
 \frac{4}{K_{\rm ref}}
 \frac{\exp\{x^\top V_{\rm ctx}v_{r(s)}\}}
 {\sum_{r=1}^{4}\exp\{x^\top V_{\rm ctx}v_r\}}.
 \label{eq:capacityConditionalProbability}
\end{equation}

CAP transforms positive demands by $[T(\xi)]_r=\log(\xi_r)$.
Let $\underline\xi_r$ and $\overline\xi_r$ be the minimum
and maximum values of component $r$ among the fixed reference scenarios.
The decoder on $\RR^{d_\xi}$ is
\begin{equation}
 [g(y)]_r=\min\{\overline\xi_r,\max\{\underline\xi_r,\exp(y_r)\}\}.
\label{eq:capacityDecoder}
\end{equation}

Differentiating Eq.~\eqref{eq:capacityComputationalProblem} gives
the scenario gradient used to construct $\widehat C_b$. Here
$\operatorname{logit}^{-1}(a):=(1+\exp(-a))^{-1}$ and $\odot$ denotes
componentwise multiplication.
\begin{equation}
 \phi_{\rm CAP}(z,\xi)
 =
 \frac1{d_\xi} A^\top
 \left[
 (c^u+c^o)\odot
 \operatorname{logit}^{-1}\!\left(\frac{Az-\xi}{0.35}\right)-c^u
 \right].
 \label{eq:capacityGradient}
\end{equation}

Since $h_\tau''(a)=\exp(-a/\tau)/[\tau\{1+\exp(-a/\tau)\}^2]>0$ is
continuous on the compact product of $\cZ_{\rm CAP}(C)$ and
$\prod_{r=1}^{d_\xi}[\underline\xi_r,\overline\xi_r]$, the diagonal entries of
$D(z,\xi):=\operatorname{Diag}((c_j^o+c_j^u)h_\tau''((Az)_j-\xi_j))$ admit a
strictly positive minimum $\underline h>0$, and full column rank of $A$ yields
$\nabla_z^2\ell_{\rm CAP}(z,\xi)=d_\xi^{-1}A^\top D(z,\xi)A\succeq
(\underline h\,\sigma_{\min}(A)^2/d_\xi)\mathrm{Id}_{d_z}$. Every CAP objective
is therefore strongly convex with
$\eta=\underline h\,\sigma_{\min}(A)^2/d_\xi>0$, so CAP sets
$\eta_{{\rm reg},{\rm CAP}}=0$ (with $z_0=5\mathbf 1_{d_z}$).

\subsection{Smoothed conditional value-at-risk portfolio optimization}
\label{ec:portfolioProblem}

We set $r_{\max}=0.25$, $\alpha=0.10$,
$\lambda_\mu=0.25$, $\lambda_w=0.01$, and $\tau=0.015$.

The master scenario set contains 1,000 scenarios for
128 assets. Let
$H_s\sim\mathcal N(0,\mathrm{Id}_6)$ denote the factor vector and
$C_s\sim\operatorname{Bernoulli}(0.12)$ the crash indicator for scenario
$s$. For asset $j$, let $a_j$ denote a normalized Gaussian loading,
$r_j\sim\operatorname{Uniform}[0.5,1.5]$ its crash exposure, and
$\mu_j^0\sim\mathcal N(0.004,0.0015^2)$ its return intercept. Given
independent $E_{sj}\sim\mathcal N(0,1)$, the return is
\begin{equation}
 \tilde{\xi}_j^{(s)}
 =
 0.245\tanh\!\left(
 \frac{\mu_j^0+0.032H_s^\top a_j+0.008E_{sj}
 -0.075C_s r_j}{0.25}
 \right).
 \label{eq:portfolioCatalog}
\end{equation}

The settings for $K_{\rm ref}$ and $d_\xi$ use nested prefixes. Contexts follow
$\mathcal N(0,\mathrm{Id}_8)$. We independently draw
$G_1^P,G_2^P,G_3^P\in\RR^{8\times7}$ with independent standard normal entries and
$\omega^P\in[-\pi,\pi]^7$ with independent uniform entries.
We construct the context map from the features
\begin{equation}
\begin{aligned}
 h_{\rm lin}^{P}(x)&=\frac{x^\top G_1^P}{\sqrt8},\qquad
 h_{\rm per}^{P}(x)=\sin\!\left(\frac{x^\top G_2^P}{\sqrt8}+\omega^P\right),\\
 h_{\rm quad}^{P}(x)&=\tanh\!\left(
 \frac{(x\odot x-\mathbf 1_8)^\top G_3^P}{\sqrt8}\right).
\end{aligned}
\label{eq:portfolioContextFeatures}
\end{equation}
Using these features, we set
\begin{equation}
\begin{aligned}
 \psi_P^0(x)
 &=\tanh\!\left\{
 h_{\rm lin}^{P}(x)
 +0.45[0.65h_{\rm per}^{P}(x)+0.35h_{\rm quad}^{P}(x)]
 \right\},\\
 \psi_P(x)
 &=\left(\psi_{P,1:6}^0(x),
 \frac{1+\psi_{P,7}^0(x)}{2}\right).
\end{aligned}
\label{eq:portfolioContextMap}
\end{equation}
For $v_s=(H_s,C_s)$ and independent
$\delta_s\sim\mathcal N(0,0.10^2)$, the conditional probabilities are
\begin{equation}
 p_s^\star(x)
 =
 \frac{
 \exp\{-\|\psi_P(x)-v_s\|_2^2/(2\cdot1.35^2)+\delta_s\}
 }{
 \sum_{r=1}^{K_{\rm ref}}
 \exp\{-\|\psi_P(x)-v_r\|_2^2/(2\cdot1.35^2)+\delta_r\}
 }.
 \label{eq:portfolioConditionalProbability}
\end{equation}

PORT uses $[T(\xi)]_r=\operatorname{atanh}(\xi_r/0.25)$ and
$[g(y)]_r=0.25\tanh(y_r)$, with
$\operatorname{dom}(g)=\RR^{d_\xi}$.

In Eq.~\eqref{eq:portfolioComputationalProblem},
\(f((w,t))=0.01\|w\|_2^2/2\), and \(Q\) contains the remaining terms,
including the linear threshold term. Writing
$\omega_{\rm tail}(w,t,\xi)=
\operatorname{logit}^{-1}\{(-\xi^\top w-t)/0.015\}$,
we differentiate \(Q\) to obtain the scenario gradient used to construct
$\widehat C_b$.
\begin{equation}
 \phi_{\rm PORT}((w,t),\xi)
 =
 \left(-[0.25+10\omega_{\rm tail}(w,t,\xi)]\xi,\
 1-10\omega_{\rm tail}(w,t,\xi)\right).
 \label{eq:portfolioGradient}
\end{equation}
The leverage penalty is part of the operational loss and gives
positive curvature in the directions of the portfolio weights, but no
uniform positive bound has been established in the threshold direction.
We therefore add a quadratic term with prespecified regularization center
\(z_0=(d_\xi^{-1}\mathbf 1_{d_\xi},0)\)
and coefficient \(\eta_{{\rm reg},{\rm PORT}}=10^{-6}\)
only in the solution mapping.
In both PORT and NET, the operational objective is convex, so
the added quadratic term makes each finite-distribution objective in
Eq.~\eqref{eq:finiteDecisionLayer} strictly convex on $\cZ$. Moreover, on
any bounded local neighborhood $\mathcal B_{\cZ}(z(P_x;x),r_{\rm loc})$,
the leverage penalty $\lambda_w\|w\|_2^2/2$ and smooth tail curvature
$h_\tau''(-\xi^\top w-t)>0$ under $P_x$ in PORT, and expected
quadratic congestion curvature on active links in NET, allow the local
strong-convexity constant $\eta$ of $R(\cdot;P_x,x)$ in
Eq.~\eqref{eq:strongConvexity} to exceed $\eta_{\rm reg}$
substantially, with $\eta\geq\eta_{\rm reg}>0$ always
providing a global fallback when $r_{\rm loc}=\infty$.

\subsection{Multicommodity network reservation}
\label{ec:reservationProblem}

The network has node set \(\{0,1,2,3\}^2\) and both
orientations of every
horizontal and vertical grid edge.
In each replication, we lexicographically order all ordered
origin--destination pairs whose row and column coordinates both differ,
apply a fixed permutation for that replication, and retain the first 128
pairs as the master commodity list. For each commodity, one path moves
horizontally and then vertically; the other moves vertically and then
horizontally.
We set \(c=\mathbf 1_{48}\), \(\lambda=10\), and \(\rho=0.05\).
We generate 1,000 reference scenarios with 128 demand components
and use nested prefixes for the settings of \(K_{\rm ref}\) and \(d_\xi\).

Define
\(b_j=0.04+(j-1)(0.06-0.04)/127\) and
\(g_j=1+\{(j-1)\bmod4\}\). For each scenario, we draw
\(H_s\sim\mathcal N(0,\mathrm{Id}_6)\), \(E_{sj}\sim\mathcal N(0,1)\),
\(R_s\sim\operatorname{Bernoulli}(0.12)\), and
\(T_s\) uniformly from \(\{1,2,3,4\}\). Define
\begin{equation}
 M=
 \begin{pmatrix}
 1&1&0&0\\0&0&1&1\\1&0&1&0\\0&1&0&1
 \end{pmatrix}.
\label{eq:reservationRegimeMatrix}
\end{equation}
Using the regime matrix in Eq.~\eqref{eq:reservationRegimeMatrix}, we
generate the master scenarios as
\begin{equation}
 \tilde{\xi}_j^{(s)}
 =
 b_j\exp\!\left\{
 0.35H_{s1}+0.40H_{s,1+g_j}+0.20E_{sj}
 +0.90R_sM_{T_s,g_j}
 -\frac{0.35^2+0.40^2+0.20^2}{2}
 \right\}.
\label{eq:reservationCatalog}
\end{equation}
We compute normalized total demand from all 128 components of
each reference scenario, even when the setting retains fewer commodities,
by
\begin{equation}
 \bar\xi_s=\frac1{128}\sum_{j=1}^{128}\frac{\tilde{\xi}_j^{(s)}}{b_j},
 \qquad s=1,\ldots,1{,}000.
\label{eq:reservationNormalizedDemand}
\end{equation}
Within the current \(K_{\rm ref}\)-scenario prefix, we order the scenarios by
the normalized demand in Eq.~\eqref{eq:reservationNormalizedDemand}, breaking
an exact tie by scenario index, and let
\(\mathcal L\) and \(\mathcal H\) contain the lower and upper
\(K_{\rm ref}/2\) scenarios, respectively. Contexts follow
\(\mathcal N(0,\mathrm{Id}_8)\). Let \(u=8^{-1/2}\mathbf 1_8\) and
\(\varpi(x)=\operatorname{logit}^{-1}(3u^\top x)\). The conditional distribution is
\begin{equation}
 p_s^\star(x)
 =
 \begin{cases}
 2\{1-\varpi(x)\}/K_{\rm ref}, & s\in\mathcal L,\\
 2\varpi(x)/K_{\rm ref}, & s\in\mathcal H.
 \end{cases}
 \label{eq:reservationConditionalProbability}
\end{equation}
The conditional distribution depends only on \(u^\top x\), and the seven
orthogonal context directions do not affect it.

NET transforms positive demands by $[T(\xi)]_r=\log(\xi_r)$
and reconstructs them by $[g(y)]_r=\exp(y_r)$,
with $\operatorname{dom}(g)=\RR^{d_\xi}$.

For fixed reservation \(z\) and demand \(\xi\), the routing
regularizer in Eq.~\eqref{eq:reservationComputationalProblem} makes the
recourse decision \(q(z,\xi)\) unique. Danskin's theorem then gives the
scenario gradient used to construct $\widehat C_b$.
\begin{equation}
 \phi_{\rm NET}(z,\xi)
 =-2\lambda[L(q(z,\xi),\xi)-z]_+,
 \label{eq:reservationGradient}
\end{equation}
where $[a]_+:=\max\{a,0\}$ is applied componentwise.
The reservation cost gradient \(c\) belongs to the deterministic cost \(f\)
and is therefore not part of the scenario gradient \(\phi=\nabla_zQ\).
NET can be flat
in reservation directions whose capacities are inactive.
The added quadratic term has prespecified regularization center
\(z_0=0.1\cdot\mathbf 1_{48}\) and coefficient
\(\eta_{{\rm reg},{\rm NET}}=10^{-6}\), ensuring strong convexity.

\subsection{Rolling interconnected grid capacity-coverage protocol}
\label{ec:grid64Protocol}

We extract hourly actual electricity load, wind generation, and solar
generation across $R=8$ interconnected Continental European control areas
(\texttt{DE\_50hertz}, \texttt{DE\_amprion}, \texttt{DE\_tennet},
\texttt{DE\_transnetbw}, \texttt{FR}, \texttt{AT}, \texttt{DK\_1}, and
\texttt{IT\_NORD}) from the Open Power System Data compilation
\citep{wiese2019open} from January 1,
2015, through September 30, 2020 ($2{,}100$ Coordinated Universal Time (UTC) operational days). Each zone
$r\in\{1,\ldots,8\}$ is scaled by its 2015--2016 median hourly load
$m_r^{\rm load}>0$ ($9.83$, $20.73$, $17.34$, $6.89$, $52.65$, $6.90$, $2.17$,
and $17.40$~GW), with missing or meter-dropout readings imputed causally from
hour $t-24$. Averaging hourly per-unit residual net load
$(\mathrm{load}-\mathrm{wind}-\mathrm{solar})/m_r^{\rm load}$ within $H=8$
three-hour UTC blocks $h\in\{1,\ldots,8\}$ defines the
$d_\xi=64$-dimensional daily realization $\xi_t\in(-2.2,2.2)^{64}$, yielding
$2{,}098$ complete observations $(x_t,\xi_t)$ after the two-day feature lag
under $[T(\xi)]_j=\operatorname{atanh}(\xi_j/2.2)$ and
$[g(y)]_j=2.2\tanh(y_j)$ on $\operatorname{dom}(g)=\RR^{64}$.

Let $\bar\xi_{s,r}^{\rm am}$ denote the mean per-unit residual net load of
zone $r$ over hours 05:00--08:00 UTC on day $s$ (ahead of the 12:00 Central
European day-ahead gate closure), $\mathcal R_{\rm DE}:=\{1,2,3,4\}$ the four
German zones, $\mathcal R_{\rm oth}:=\{6,7,8\}$,
$\omega_t:=2\pi\,\mathrm{doy}(t)/365.25$, and $w_t^{\rm wknd}\in\{0,1\}$ the
weekend indicator. For operational day $t$, the context vector $x_t\in\RR^{16}$ is
\begin{equation}
 \begin{aligned}
 x_t=\Biggl(
 &\bar\xi_{t-1,1}^{\rm am},\ldots,\bar\xi_{t-1,8}^{\rm am},\;
  \frac{1}{4}\sum_{r\in\mathcal R_{\rm DE}}\xi_{t-2,(r,7)},\;
  \xi_{t-2,(5,7)},\;
  \frac{1}{3}\sum_{r\in\mathcal R_{\rm oth}}\xi_{t-2,(r,7)},\\
 &\sin\omega_t,\;\cos\omega_t,\;\sin(2\omega_t),\;\cos(2\omega_t),\;
  w_t^{\rm wknd}
 \Biggr)^\top.
 \end{aligned}
 \label{eq:grid64ContextVector}
\end{equation}

Let $A_{\rm grid}\in\{0,1\}^{8\times 8}$ denote the symmetric adjacency
matrix of the $13$ undirected corridors among the eight zones
(\texttt{DE\_50hertz}--\texttt{DE\_tennet},
\texttt{DE\_50hertz}--\texttt{DK\_1},
\texttt{DE\_amprion}--\texttt{DE\_tennet},
\texttt{DE\_amprion}--\texttt{DE\_transnetbw},
\texttt{DE\_amprion}--\texttt{FR},
\texttt{DE\_amprion}--\texttt{AT},
\texttt{DE\_tennet}--\texttt{DE\_transnetbw},
\texttt{DE\_tennet}--\texttt{AT},
\texttt{DE\_tennet}--\texttt{DK\_1},
\texttt{DE\_transnetbw}--\texttt{FR},
\texttt{DE\_transnetbw}--\texttt{AT},
\texttt{FR}--\texttt{IT\_NORD}, and
\texttt{AT}--\texttt{IT\_NORD})
and $D_{\rm deg}:=\operatorname{Diag}(A_{\rm grid}\mathbf 1_8)$. With sharing
fraction $\gamma_{\rm cpl}=0.20$,
$M_{\rm grid}:=(1-\gamma_{\rm cpl})\mathrm{Id}_8+\gamma_{\rm cpl}A_{\rm grid}D_{\rm deg}^{-1}$
maps zonal per-unit commitments $z\in\cZ_{\rm GRID}:=[0,2.5]^8$ to effective per-unit
coverage $[M_{\rm grid}z]_r$. Because each zone is normalized by $m_r^{\rm load}$,
$M_{\rm grid}$ operates in normalized zonal coordinates and does not represent
physical gigawatt transfers.
With $c_0=2.5$, $\kappa_{\rm sc}=70$, block
weights $w_{\rm blk}=(0.10,0.10,1.85,1.35,0.80,1.60,2.10,0.50)$,
$c_h^u:=c_0\kappa_{\rm sc}w_{{\rm blk},h}$, and $\tau=0.028$, the operational
loss is
\begin{equation}
 \ell_{\rm GRID}(z,\xi)
 :=c_0\mathbf 1_8^\top z
 +\frac{1}{8}\sum_{r=1}^8\sum_{h=1}^8
 c_h^u\,h_\tau\!\bigl(\xi_{r,h}-[M_{\rm grid}z]_r\bigr),
 \label{eq:grid64DecisionProblem}
\end{equation}
augmented during solving by $(\eta_{\rm reg}/2)\|z-z_0\|_2^2$ around
$z_0=\mathbf 1_8$ with $\eta_{\rm reg}=10^{-6}$.

For each quarterly window $q\in\{\text{2018Q1},\ldots,\text{2020Q3}\}$
($N_{\rm oos}=1{,}004$ days), validation uses the preceding 364 days and
training uses the 728 days prior to validation. Atlas selects
$\widehat B\in\{1,2,4\}$ and $d_b$ via five-fold assessment and Eq.~\eqref{eq:spectralWidthRule}
($\alpha_{\rm ret}=0.9999$, $n_{\min}=50$,
$c_{\min}^{\rm par}=1.5$), sets
$K_b=d_b+1$, and runs $1{,}600$ Adam updates at learning rate $0.025$;
ICEO and ICEO-red use learning rate $0.001$, and all remaining method and
validation rules follow Section~\ref{ec:methodTrainingRules}. Let
$\ell_{m,1},\ldots,\ell_{m,N_{\rm oos}}$ denote the concatenated daily
out-of-sample losses of method $m$ with mean $\bar\ell_m$, centered series
$\tilde\ell_{m,t}:=\ell_{m,t}-\bar\ell_m$, and sample autocovariances
$\hat\gamma_{m,l}:=N_{\rm oos}^{-1}\sum_{t=l+1}^{N_{\rm oos}}\tilde\ell_{m,t}\tilde\ell_{m,t-l}$.
With a two-week bandwidth ($L_{\rm NW}=14$), we compute the
\citet{newey1987simple} standard error
\begin{equation}
 \mathrm{SE}_{\rm NW}(\bar\ell_m)
 :=\left[
 \frac{1}{N_{\rm oos}}\left\{
 \hat\gamma_{m,0}
 +2\sum_{l=1}^{L_{\rm NW}}\left(1-\frac{l}{L_{\rm NW}+1}\right)\hat\gamma_{m,l}
 \right\}_+
 \right]^{1/2},
 \label{eq:grid64NeweyWest}
\end{equation}
and apply Eq.~\eqref{eq:grid64NeweyWest} analogously to paired daily loss
differences $\ell_{{\rm Atlas},t}-\ell_{m,t}$. Table~\ref{tab:grid64EarlyDetails}
reports the detailed first-period results (over 2018Q1--2020Q3, the
prescriptiveness coefficient of \citealt{bertsimas2020predictive} relative to
SAA and hindsight optimum $20.76$ is $0.396$ for Atlas versus $0.390$ for ICEO).

\begin{table}[htbp]
\centering
\caption{Rolling day-ahead European grid capacity coverage across 11 quarterly windows (2018Q1--2020Q3, 1,004 days).}
\label{tab:grid64EarlyDetails}
\small
\renewcommand{\arraystretch}{1.03}
\setlength{\tabcolsep}{2pt}
\begin{tabular*}{\textwidth}{@{\extracolsep{\fill}}l*{8}{r}@{}}
\toprule
 & & \multicolumn{3}{c}{Out-of-sample operational loss $\downarrow$} & & \multicolumn{3}{c}{Deploy wall time $\downarrow$} \\
\cmidrule(lr){3-5}\cmidrule(lr){7-9}
Method & Scenarios & 11-qtr med. & Pre-2020 & Pooled (SE) & Build\,/s $\downarrow$ & 1 ctx\,/ms & 128, 1c\,/s & 128, 32c\,/ms \\
\midrule
Atlas & 8.4 & \textbf{24.05} & \textbf{25.65} & \textbf{25.10} (0.32) & 87.40 & 19.58 & 1.19 & 65.02 \\
ICEO & 728.0 & 24.16 & 25.73 & 25.15 (0.35) & 355.72 & 23.90 & 4.62 & 263.48 \\
ICEO-red & 8.4 & 26.65 & 27.72 & 26.95 (0.35) & 81.78 & 18.94 & \textbf{0.92} & \textbf{50.29} \\
SOF & 106.4 & 24.38 & 25.85 & 25.38 (0.30) & 10.47 & 13.39 & 2.13 & 181.04 \\
kNN & 19.6 & 24.14 & 26.02 & 25.51 (0.31) & 2.78 & \textbf{9.18} & 1.02 & 68.05 \\
kNN-red & 8.3 & 24.70 & 26.49 & 25.85 (0.37) & 60.12 & 13.20 & 1.05 & 77.39 \\
SAA & 728.0 & 27.59 & 28.22 & 27.96 (0.17) & $<0.01$ & 24.25 & 3.09 & 132.85 \\
\bottomrule
\end{tabular*}
\par\vspace{2pt}
\begin{minipage}{\textwidth}
\renewcommand{\baselinestretch}{1.02}
\footnotesize
\textit{Notes.} Scenarios average over the 11 windows;
11-qtr med.\ is the median quarterly loss. Pre-2020 (730 days) and Pooled
(1,004 days) give daily means, with 14-day Newey--West standard errors
(SE) reported only in the Pooled column.
Build\,/s and deployment times are measured on the 2019Q1 window
($K_1=8$ versus 728 scenarios) without competing workloads.
Deployment times are medians over five runs. The 32-core batch timing
for kNN-red uses parallel clustering across 32 cores.
\end{minipage}
\end{table}

\subsection{Evaluation on later data}
\label{ec:grid64Later}

For 2020Q4--2024Q4, we extract hourly load, wind, and solar generation for the
same eight control areas from the ENTSO-E Transparency Platform
\citep{hirth2018entsoe} from October 1, 2017, through December 31, 2024,
keeping the per-unit bases
$m_r^{\rm load}$, imputation rule, context vector in
Eq.~\eqref{eq:grid64ContextVector}, loss in Eq.~\eqref{eq:grid64DecisionProblem},
and training rules of Section~\ref{ec:grid64Protocol} fixed in advance. Each
window trains on 728 days and validates on 364 days from this compilation.
Table~\ref{tab:grid64Later} reports the results; Atlas selects one chart in all
17 windows ($K_1=8$ through 2022Q4 and $K_1=9$ from 2023Q1).

\begin{table}[htbp]
\centering
\caption{Rolling day-ahead European grid capacity coverage on later data (17 quarterly windows, 2020Q4--2024Q4, 1,553 days).}
\label{tab:grid64Later}
\small
\renewcommand{\arraystretch}{1.03}
\setlength{\tabcolsep}{4pt}
\begin{tabular*}{\textwidth}{@{\extracolsep{\fill}}l*{6}{r}@{}}
\toprule
Method & Scenarios & 17-qtr med. & Pooled (SE) & Presc. & Atlas $-$ method (SE) & Quarters \\
\midrule
Atlas & 8.5 & 25.14 & 24.62 (0.28) & 0.401 &  &  \\
ICEO & 728.0 & \textbf{24.83} & \textbf{24.30} (0.28) & \textbf{0.443} & $+0.32$ (0.08) & 3 \\
ICEO-red & 8.5 & 27.18 & 27.22 (0.55) & 0.055 & $-2.61$ (0.34) & 16 \\
SOF & 117.4 & 26.24 & 25.10 (0.27) & 0.337 & $-0.48$ (0.10) & 14 \\
kNN & 17.9 & 26.00 & 25.07 (0.29) & 0.342 & $-0.45$ (0.10) & 14 \\
kNN-red & 8.5 & 25.79 & 25.28 (0.35) & 0.314 & $-0.66$ (0.14) & 13 \\
SAA & 728.0 & 27.36 & 27.64 (0.17) & 0.000 & $-3.03$ (0.18) & 17 \\
\bottomrule
\end{tabular*}
\par\vspace{2pt}\begin{minipage}{\textwidth}\renewcommand{\baselinestretch}{1.02}\footnotesize
\textit{Notes.} Scenarios average over the 17 windows; 17-qtr med.\ is the median quarterly loss. Pooled (SE) gives daily means over 1,553 days with 14-day Newey--West standard errors. Presc. is prescriptiveness (hindsight-optimal loss $20.10$). Atlas $-$ method is the mean paired daily difference (with 14-day SE), and Quarters counts windows where Atlas has lower loss.
\end{minipage}
\end{table}

\section{Computational Protocol}
\label{ec:protocol}

The source code, data-generation and preprocessing scripts, and
workflows for reproducing all experiments will be released as open source
upon publication.

\subsection{Sampling and evaluation}
\label{ec:samplingEvaluation}
\label{ec:symmetricApplicationPanel}

For each setting within a replication, all methods use the same
training and validation observations and are evaluated on the same
contexts and realizations.
We generate the reference scenario sets, context maps, and
training, validation, and evaluation samples using separate random streams
within each replication. PORT and NET also vary their physical problem
parameters across replications. CAP keeps its service technology fixed
and varies the resource demand states, contextual directions, and noise
within each state. The main comparison uses 10 paired replications in CAP
and PORT (random seeds 260911--260920) and three in NET (seeds 260911--260913)
for every method, with $|I_{\rm tr}|=12{,}800$, $|I_{\rm val}|=1{,}000$, and
$|I_{\rm eval}|=200$. The synthetic construction-timing, sample-size, and chart-sensitivity
studies (Figure~\ref{fig:chartSensitivity}) use replications 1--3
(seeds 260911--260913). Synthetic deployment timing uses replication 1
(seed 260911), while the controlled basis study
(Table~\ref{tab:controlledBasisMeans}) uses seeds 241011--241013.
Let \(\{x_i\}_{i\in I_{\rm eval}}\) denote the evaluation contexts,
independent of the training and validation samples. For each synthetic
application with reference support \(\{\tilde{\xi}^{(s)}\}_{s=1}^{K_{\rm ref}}\) and known conditional
probabilities \(p_s^\star(x_i)\), we evaluate the exact conditional operational
risk of any decision rule $\pi(\cdot)$ at $x_i$ by
\begin{equation}
 R_0(\pi(x_i);x_i)
 =
 \sum_{s=1}^{K_{\rm ref}}p_s^\star(x_i)\ell(\pi(x_i),\tilde{\xi}^{(s)};x_i),
 \label{eq:exactCatalogEvaluationRisk}
\end{equation}
and report the mean of Eq.~\eqref{eq:exactCatalogEvaluationRisk} across
$i\in I_{\rm eval}$.
For NET, we evaluate the operational loss by solving the recourse
problem in Eq.~\eqref{eq:reservationComputationalProblem} to convergence
with the reservation decision fixed.

\subsection[Regularization and structure selection]{Regularization and structure selection}
\label{ec:regularizationSelection}
\label{ec:commonStructureProtocol}

Across all fitted methods, settings, and replications, downstream solves
use $\eta_{\rm reg}=0$ in CAP and $\eta_{\rm reg}=10^{-6}$ in PORT, NET, and
the grid study, fixed in advance; $\eta_{\rm reg}$ enters only the downstream
solves and their derivatives, while every reported loss and risk is
unregularized. For CAP, PORT, NET, and the European grid study, Atlas
considers $B\in\{1,2,4\}$. The chart sensitivity study in
Section~\ref{subsec:chartMechanismResults} and its standalone check in
Table~\ref{tab:ecRobustnessChecks} use $B\in\{1,2,4,8,16,32\}$.
We use K-means with $20$ restarts and tolerance $10^{-8}$ for context charts
and $10$ restarts for scenario initialization. Each candidate is evaluated
using five internal folds per replication, with $\alpha_{\rm ret}=0.9999$
(Table~\ref{tab:ecAlphaSensitivity}).

A candidate is eligible if every chart contains
at least $n_{\min}=50$ observations, achieves the retention target within its
effective rank (counting singular values satisfying
$\sigma_j(\widehat C_b)>10^{-12}\sigma_1(\widehat C_b)$),
has at least $K_b$ distinct projected realizations
$\widehat V_b^\top T_b(\xi_i)$ (within tolerance $10^{-10}$), and has at least $c_{\min}^{\rm par}=3$
($1.5$ on the grid, owing to its shorter 728-day training windows) fitting observations per nominal free parameter
$P_b:=K_bd_b+(K_b-1)(d_x+1)$ across folds.
We solve Eq.~\eqref{eq:screeningProbabilityFit} using L-BFGS-B,
initialized from empirical cluster frequencies, with
$\mathrm{ftol}=10^{-12}$, $\mathrm{gtol}=10^{-8}$, and
$\mathrm{maxiter}=2{,}000$. We use an intercept-only fallback when needed
and reject candidates with empty clusters in assessment folds.

After fixing the structure, we jointly train all $\widehat B$ charts
for $1{,}600$ synchronized Adam steps. At each step, we sample uniformly
within each $I_b$ and distribute a global minibatch of $64$ observations
as evenly as possible across charts, giving
$\lfloor 64/\widehat B\rfloor$ or $\lceil 64/\widehat B\rceil$
observations per chart. We clip the gradient norm at $10.0$ for each
chart and project $(a_b,\Gamma_b)$ onto the zero-sum constraints after
each update. Finally, we select one common optimization
iteration from $\{0,400,800,1200,1600\}$ across all charts using
validation operational loss.

\subsection{Training rules by method}
\label{ec:methodTrainingRules}

ICEO standardizes contexts and maps them through a two-hidden-layer ReLU
network of width 32 to $K_{\rm ref}$ softmax logits over the reference support,
minimizing empirical operational loss. Atlas uses penalties
$\lambda_u=\lambda_\Gamma=10^{-5}$ and learning rates of $0.001$ in CAP,
PORT, and NET and $0.025$ on the European grid. The grid learning rate was
selected using the initial development windows and fixed before evaluating
the 17 later quarters.

We solve all downstream problems to a KKT residual tolerance of $10^{-6}$.
CAP uses SciPy SLSQP ($\mathrm{ftol}=10^{-13}$, at most $1{,}000$ iterations),
PORT uses Clarabel 0.11.1 ($\mathrm{cone\_tol}=10^{-9}$) with SLSQP when needed,
and NET uses Gurobi 13.0.2 with barrier and primal-simplex methods, followed
by up to 50 projected-gradient refinement steps. The grid study uses
projected Newton with Armijo line search and L-BFGS-B when needed
($\mathrm{ftol}=10^{-15}$, $\mathrm{gtol}=10^{-12}$, at most $500$ iterations).
All observed net loads lie within the decoder domain $(-2.2,2.2)$:
the ranges are $[-1.69,1.71]$ over 2,098 OPSD days and
$[-1.41,1.65]$ over 2,647 preprocessed ENTSO-E days.

For CAP, PORT, and the grid study, we compute training gradients by
implicit differentiation of the downstream KKT conditions.
Let $F(z;\tilde\xi,p)$ denote the regularized downstream objective,
$z^\star$ its minimizer, and $g_z:=\nabla_z\ell(z^\star,\xi;x)$
the gradient of the realized operational loss. Let
$H:=\nabla_{zz}^2F(z^\star;\tilde\xi,p)$ and let $G$ contain the
equality constraints and active bound constraints.
At a regular solution, we obtain the adjoint vector $v$ from
\[
\begin{bmatrix}
H & G^\top\\
G & 0
\end{bmatrix}
\begin{bmatrix}
v\\
\mu
\end{bmatrix}
=
\begin{bmatrix}
g_z\\
0
\end{bmatrix}.
\]
The parameter derivatives then follow from
$-\nabla_{(\tilde\xi,p)}
(\nabla_zF(z^\star;\tilde\xi,p)^\top v)$.
For active-set identification, SLSQP uses a bound tolerance of $10^{-8}$,
while Clarabel uses a primal tolerance of $10^{-6}$ and requires
the corresponding multiplier to exceed $10^{-6}$.
PORT has no active box constraints.
On the grid, bounds with multipliers exceeding $10^{-5}$ are treated
as active, while weakly active bounds are treated as free.
The clamped exponential decoder in CAP uses zero derivatives outside
its clipping bounds.

For NET, we express routing flows as fractions
$r_{k,m}:=q_{k,m}/\xi_{k,m}\in[0,1]$, so the routing constraints
do not depend on the learned scenarios. We check strict
complementarity in the original routing variables at tolerance $10^{-8}$,
using 80-digit \texttt{mpmath} refinement when needed.
Eliminating the independent routing blocks through a Schur complement
gives a positive-definite system for the reservation variables under
$\eta_{\rm reg}=10^{-6}$. The outer loss gradients follow from the
envelope theorem, holding the optimized routing flows fixed.

We include all solver refinement and alternative solve methods in the
reported wall times. If the required KKT conditions or derivative systems
cannot be verified, the run terminates without excluding observations
or omitting updates. All reported experiments completed successfully
without excluding any prescribed evaluation observations.
We also compare the analytic derivatives with central finite differences
at smooth test points. The median relative errors for scenario and
probability derivatives are $7.1\times10^{-6}$ and $2.7\times10^{-6}$
in CAP, $2.0\times10^{-4}$ and $9.1\times10^{-5}$ in PORT,
$2.7\times10^{-7}$ and $4.0\times10^{-8}$ in NET, and
$3.0\times10^{-8}$ and $6.5\times10^{-10}$ on the grid.
For PORT, the comparison uses $\mathrm{cone\_tol}=10^{-12}$
and finite-difference step size $10^{-3}$.

SOF uses the constrained approximate-risk
criterion of \citet{kallus2023stochastic} with 100 bootstrap trees, no honesty
split, minimum leaf size 10, and minimum child fraction $0.2$, with these
settings fixed without validation tuning. For kNN,
validation selects $k\in\{2^j:j=0,\ldots,\lfloor\log_2|I_{\rm tr}|\rfloor\}\cup\{|I_{\rm tr}|\}$
under standardized Euclidean context distances. Let
$K_{\rm match}:=\max_{1\leq b\leq\widehat B}K_b$ denote the maximum Atlas chart
support size. In standardized $T$-coordinates, ICEO-red partitions training
realizations into $K_{\rm match}$ K-means clusters
\citep{rujeerapaiboon2022scenario}, decodes the centroids through $g$, and fits
the same two-hidden-layer network over those $K_{\rm match}$ scenarios;
kNN-red independently selects $k$ by validation, applying weighted
K-means to each local empirical distribution in standardized $T$-coordinates.
We use three K-means initializations and at most 100 iterations, with the
cluster centroids decoded through $g$. When $k\leq K_{\rm match}$,
we retain the original empirical distribution without clustering.
For 32-core batch deployment, clustering is performed in parallel
across workers.

CondGen uses the same two-hidden-layer ReLU network of width 32 as ICEO
to generate context-dependent scenario offsets in standardized
$T$-coordinates, initialized at K-means centers. For $K_{\rm match}>1$,
the network also produces $K_{\rm match}$ softmax logits.
We center the output-layer weights and biases after each update,
removing 33 redundant parameters from the reported parameter counts.
CondGen minimizes the same operational loss using Adam with learning
rate $0.001$, 1,600 updates, and the same validation-based selection rule.

\subsection[Performance measures and timing]{Performance measures and timing}
\label{ec:qualityEndpoints}
\label{ec:timingProtocol}

For method $h$, let $\pi_h$ denote its decision
rule and $z_i^{\rm ref}$ an unregularized optimal decision under the known
conditional distribution at evaluation context $x_i$. We report relative excess risk as
\begin{equation}
 \mathrm{RER}_h
 =
 \frac{\sum_i\{R_0(\pi_h(x_i);x_i)-R_0(z_i^{\rm ref};x_i)\}}
 {\left|\sum_iR_0(z_i^{\rm ref};x_i)\right|}.
\label{eq:relativeExcessRate}
\end{equation}

All construction and deployment times in
Tables~\ref{tab:capPortConstructionCost},
\ref{tab:deploymentTiming}, \ref{tab:grid64Rolling},
and~\ref{tab:grid64EarlyDetails} are measured without concurrent
experiments on a dedicated group of 32 physical cores.
Construction time includes structure selection, model fitting, and
validation, but excludes final test evaluation and diagnostic evaluations
of the full training loss for all methods. For ICEO in NET replication 2,
we sum the times of both segments of a run resumed after active-set
reclassification. We report elapsed wall time and total CPU time across
the main process and workers, with CPU time measured at the operating
system's $1/64$~s resolution.

Deployment time reports five-repetition medians after two warm-up
runs from distribution prediction through decision assembly.
For single-core deployment of 128 contexts in NET, baselines
requiring more than $200$~ms per context use three timed runs after
one warm-up run. Synthetic deployment timings use replication 1
(seed 260911), with the first evaluation context for single-context
timing and the first 128 contexts for batch timing. Grid deployment
timings use 2019Q1, with the first test day for single-context timing
and the preceding 128 validation days for batch timing, since the test
quarter contains only 90 days.

\subsection{Chart experiment guided by theory}
\label{ec:chartMechanismProtocol}

For the chart problem in Section~\ref{subsec:chartMechanismResults}, we use
the uncertainty space and feasible decision set \([-4,4]^{10}\) and set
\(\eta_{\rm reg}=0\).
The Hessian in \(z\) is the identity.
The transformation and decoder are the identity on this box.
We reject supports containing scenarios outside it, so the
weighted mean is both the unconstrained optimum and a feasible optimum of
the box-constrained problem.

For each replication, the QR factorization of a standard Gaussian
\(10\times10\) matrix provides mutually orthogonal two-column blocks
\(V_{\rm c},V_1^{\rm aux},\ldots,V_4^{\rm aux}\).
The operating mode \(R\) is uniform on
\(\{1,2,3,4\}\). We use the scenario construction in Section~\ref{subsec:chartMechanismResults}. The four modes give 16 corners, with only four distinct scenarios at zero angle.
We reuse the random matrix, mode draws, contexts, and uniforms
for sampling scenario indices across angles, so only the mode-specific
directions change.

Let \(h_r\) enumerate \((-1,-1),(-1,1),(1,-1),(1,1)\), and set
\(c_r=2(h_{r1},h_{r1},h_{r1},h_{r1},h_{r2},h_{r2})\).
We independently draw \(c_e\) uniformly on \([-1/4,1/4]^6\) and \(s\) uniformly
on \([-0.8,0.8]^2\). The observed context is
\(x=(c_R+c_e,t)\), where \(t=\operatorname{atanh}(s)\) componentwise.
The context space \(\cX\) is the union of four products, each
consisting of a mode box and the two-dimensional domain of \(t\).
We do not supply the true mode labels during fitting.
Only the four corners of the observed mode have positive mass. Relative excess risk in Eq.~\eqref{eq:relativeExcessRate} uses the exact conditional mean in Eq.~\eqref{eq:chartMechanismMain}. Each value of $B$ in Figure~\ref{fig:chartSensitivity} is fitted separately.

\subsection{Ablation and sensitivity protocols}
\label{ec:additionalAblations}

The experiments in Section~\ref{subsec:ablationStudies} use the following
protocols. For the controlled local-basis study in
Section~\ref{subsec:controlledBasisResults}
(Table~\ref{tab:controlledBasisMeans}), we use three replications
(seeds 241011--241013). In each replication, we generate orthonormal
bases $[V_{\rm rel},V_{\rm noise}]\in\RR^{200\times12}$ using QR
decomposition. We construct five task groups at the vertices of a
four-dimensional regular simplex with radius $2$, scaled by $1.5$ to
norm $3$. Each group contains ten centered Gaussian microstates,
scaled to coordinate RMS $0.1$, crossed with ten centered and whitened
nuisance states in $\RR^8$, giving $K_{\rm ref}=500$ scenarios.
We draw contexts $x\in\RR^2$ from four unit-variance Gaussian clusters
centered at $5(\pm1,\pm1)$. Conditional group probabilities follow a
zero-intercept softmax model using the context relative to its nearest
cluster center, with five slope vectors
$1.4(\cos\tfrac{2\pi g}{5},\sin\tfrac{2\pi g}{5})$.
The $100$ scenarios within each group have equal conditional probabilities.
For each $\sigma\in\{1/4,4\}$, we generate
$|I_{\rm tr}|=1{,}600$, $|I_{\rm val}|=200$, and
$|I_{\rm eval}|=200$ observations. Both basis methods use
$B=4$, $d_b=4$, $K_b=5$, $\eta_{\rm reg}=0$, and
$\lambda_u=\lambda_\Gamma=10^{-5}$, with $80$ Adam iterations
at learning rate $0.01$ and evaluations at iterations
$\{0,20,40,80\}$.

The basis comparisons, reported in Table~\ref{tab:ecBasisAblation},
retain the chart partition and support sizes $(\widehat B,\{K_b\})$ selected by
Atlas on each replication or quarterly window. PCA uses the leading $d_b$
principal components of $\{T_b(\xi_i):i\in I_b\}$, and Full sets
$\widehat V_b=\mathrm{Id}_{d_\xi}$. Each retention threshold $\alpha_{\rm ret}$ reruns
Algorithm~\ref{alg:kktAtlas} from structure selection through joint
model fitting.

Each training sample size uses a nested prefix of the 12,800 training
observations, with the 1,000-observation validation sample and 200-context
evaluation panel fixed. For each sample size, kNN-red selects $k$
independently by validation (Table~\ref{tab:ecSampleScaling}).

We also compare construction times for Atlas and alternative scenario
bases under the same 32-core timing protocol
(Section~\ref{ec:timingProtocol}). For Atlas, we separately measure
structure selection (Step~1) and parameter fitting (Step~2), which
includes initialization, $1{,}600$ Adam iterations, and validation.
Over replications 1--3, the mean Step~1, Step~2, and total wall times
are $41.13$, $65.08$, and $106.22$ seconds in CAP, $41.11$, $48.20$,
and $89.32$ seconds in PORT, and $366.46$, $1{,}440.36$, and
$1{,}806.82$ seconds in NET. The NET results use the selected
chart counts $\widehat B=(4,2,4)$ from
Table~\ref{tab:capPortStructures} and agree with the construction
times in Table~\ref{tab:capPortConstructionCost}.

For comparison, we fit PCA and full-space models using the same
selected chart structures. Their respective Step~2 wall times are
$64.08$ and $64.15$ seconds in CAP, $47.25$ and $47.31$ seconds in PORT,
and $1{,}544.42$ and $1{,}462.91$ seconds in NET.
On the 2019Q1 grid window, Atlas requires $2.57$ seconds for Step~1
and $87.96$ seconds for Step~2, giving a total of $90.53$ seconds,
compared with $87.40$ seconds in the separate run reported in
Table~\ref{tab:grid64Rolling}. The corresponding fixed-structure
Step~2 times are $77.05$ seconds for PCA and $77.59$ seconds for
full-space fitting. Thus, across the four applications, Step~2
fitting times are comparable across the three bases. On the grid window,
the KKT fit takes about $14\%$ longer than the PCA fit.

\begin{table}[!htbp]
\centering
\caption{Local basis ablation across synthetic applications and rolling European grid windows.}
\label{tab:ecBasisAblation}
\footnotesize
\renewcommand{\arraystretch}{1.00}
\setlength{\tabcolsep}{0pt}
\begin{tabular*}{\textwidth}{@{\extracolsep{\fill}}llrrrr@{}}
\toprule
Study & Basis variant & Parameters & Metric (SE) $\downarrow$ & \shortstack[r]{KKT $-$ variant\\(SE)} & \shortstack[r]{KKT\\lower} \\
\midrule
\multirow{3}{*}{CAP (10 reps.)} & Cross-covariance (KKT) & 39--78 & \textbf{0.07} (0.01) &  &  \\
 & Standardized PCA & 39--78 & \textbf{0.07} (0.01) & $0.00$ (0.00) & 4/10 \\
 & Full space ($\RR^{64}$) & 283--566 & \textbf{0.07} (0.01) & $0.00$ (0.01) & 5/10 \\
\midrule
\multirow{3}{*}{PORT (10 reps.)} & Cross-covariance (KKT) & 144--459 & 2.94 (0.17) &  &  \\
 & Standardized PCA & 144--459 & 3.17 (0.17) & $-0.23$ (0.08) & 10/10 \\
 & Full space ($\RR^{32}$) & 360--729 & \textbf{2.75} (0.18) & $+0.18$ (0.05) & 1/10 \\
\midrule
\multirow{3}{*}{NET (3 reps.)} & Cross-covariance (KKT) & 371--800 & \textbf{0.12} (0.01) &  &  \\
 & Standardized PCA & 371--800 & 0.13 (0.01) & $-0.01$ (0.00) & 3/3 \\
 & Full space ($\RR^{16}$) & 507--1,064 & \textbf{0.12} (0.01) & $0.00$ (0.00) & 2/3 \\
\midrule
\multirow{3}{*}{Grid 2018Q1--2020Q3 (11 qtrs)} & Cross-covariance (KKT) & 175--208 & \textbf{25.10} (0.32) &  &  \\
 & Standardized PCA & 175--208 & 25.29 & $-0.18$ (0.12) & 9/11 \\
 & Full space ($\RR^{64}$) & 631--712 & 25.38 & $-0.28$ (0.15) & 8/11 \\
\midrule
\multirow{3}{*}{Grid 2020Q4--2024Q4 (17 qtrs)} & Cross-covariance (KKT) & 175--208 & \textbf{24.62} (0.28) &  &  \\
 & Standardized PCA & 175--208 & 25.16 & $-0.54$ (0.24) & 11/17 \\
 & Full space ($\RR^{64}$) & 631--712 & 24.91 & $-0.30$ (0.16) & 10/17 \\
\midrule
\multirow{3}{*}{Grid 2018Q1--2024Q4 (28 qtrs)} & Cross-covariance (KKT) & 175--208 & \textbf{24.81} &  &  \\
 & Standardized PCA & 175--208 & 25.21 & $-0.40$ (0.15) & 20/28 \\
 & Full space ($\RR^{64}$) & 631--712 & 25.10 & $-0.29$ (0.11) & 18/28 \\
\bottomrule
\end{tabular*}
\par\vspace{2pt}\begin{minipage}{\textwidth}\renewcommand{\baselinestretch}{1.02}\footnotesize
\textit{Notes.} Metric is mean relative excess risk (\%) with cross-replication SE on CAP, PORT, and NET, and pooled daily loss on the grid (PCA/Full grid rows report pooled means; 14-day Newey--West SEs for grid differences). KKT $-$ variant is the KKT minus row-variant metric; KKT lower is the fraction of replications or quarters where KKT has strictly lower risk or loss.
\end{minipage}
\end{table}

\begingroup
\setlength{\intextsep}{6pt}
\setlength{\floatsep}{6pt}
\subsection{Additional comparisons and robustness diagnostics}
\label{ec:round7Diagnostics}

Tables~\ref{tab:ecCondGenAll}--\ref{tab:ecHeadMatched} give the detailed results for Section~\ref{subsubsec:round7Robustness}.

\begin{table}[!htbp]
\centering

\caption{Shared-support KKT Atlas versus context-varying neural scenario generator (CondGen) at matched $K_{\rm match}$.}
\label{tab:ecCondGenAll}
\footnotesize
\renewcommand{\arraystretch}{1.01}
\setlength{\tabcolsep}{1.5pt}
\begin{tabular*}{\textwidth}{@{\extracolsep{\fill}}lrrrrrrr@{}}
\toprule
Application & \multicolumn{2}{c}{Step~2 Params $P$} & \multicolumn{2}{c}{Metric (SE) $\downarrow$} & \shortstack[r]{Atlas $-$ CondGen\\(paired SE)} & \multicolumn{2}{c}{128 ctx, 32c\,/s $\downarrow$} \\
\cmidrule(lr){2-3}\cmidrule(lr){4-5}\cmidrule(lr){7-8}
 & Atlas & CondGen & Atlas & CondGen & & Atlas & CondGen \\
\midrule
CAP (10 reps.) & 39--78 & 9,891 & \textbf{0.07} (0.01) & 0.49 (0.03) & $-0.42$ (0.03) & \textbf{0.04} & \textbf{0.04} \\
PORT (10 reps.) & 144--459 & 11,112--20,913 & \textbf{2.94} (0.17) & 5.17 (0.23) & $-2.23$ (0.20) & \textbf{0.03} & \textbf{0.03} \\
NET (3 reps.; default) & 371--800 & 6,921--7,482 & 0.12 (0.01) & \textbf{0.10} (0.00) & $+0.02$ (0.01) & \textbf{0.21} & 0.22 \\
NET (3 reps.; $B=1$) & 171--200 & 6,921--7,482 & \textbf{0.07} (0.01) & 0.10 (0.00) & $-0.03$ (0.01) & \textbf{0.21} & 0.22 \\
Grid 18Q1--20Q3 (11 qtrs) & 175--208 & 18,727--20,872 & 25.10 (0.32) & \textbf{24.98} (0.32) & $+0.12$ (0.13) & 0.07 & \textbf{0.04} \\
Grid 20Q4--24Q4 (17 qtrs) & 175--208 & 18,727--20,872 & 24.62 (0.28) & \textbf{24.22} (0.26) & $+0.40$ (0.08) & 0.07 & \textbf{0.04} \\
\bottomrule
\end{tabular*}
\par\vspace{2pt}\begin{minipage}{\textwidth}\renewcommand{\baselinestretch}{1.00}\footnotesize
\textit{Notes.} Metrics and SEs follow Table~\ref{tab:ecBasisAblation}. CondGen reuses the default Atlas budget (NET replication~2: 11 scenarios versus 10 for single-chart Atlas). Both grid timings come from 2019Q1. CondGen build times in CAP, PORT, and NET are 50.48, 43.26, and 630.59 seconds.
\end{minipage}
\end{table}

\begin{table}[!htbp]
\centering

\caption{Atlas sample-splitting, single-chart, probability-head, and multi-anchor diagnostics.}
\label{tab:ecAtlasDiagnostics}
\footnotesize
\renewcommand{\arraystretch}{1.01}
\setlength{\tabcolsep}{3pt}
\begin{tabular*}{\textwidth}{@{\extracolsep{\fill}}lrrrrr@{}}
\toprule
Application & \shortstack[r]{Atlas\\(default)} & \shortstack[r]{Atlas-split\\(50/50)} & \shortstack[r]{Atlas\\($B=1$)} & \shortstack[r]{Atlas-MLP-prob} & \shortstack[r]{5-anchor Atlas\\(min cos.)} \\
\midrule
CAP (10 reps.) & 0.07 (0.01) & 0.10 (0.01) & 0.07 (0.01) & 0.51 (0.07) & 0.07 (0.9998) \\
PORT (10 reps.) & 2.94 (0.17) & 3.43 (0.26) & 2.94 (0.17) & 3.50 (0.34) & 2.94 (0.9815) \\
NET (3 reps.) & 0.12 (0.01) & 0.19 (0.03) & 0.07 (0.01) & 0.24 (0.03) & 0.12 (0.9968) \\
Grid 18Q1--20Q3 (11 qtrs) & 25.10 (0.32) & 26.17 (0.48) & 25.10 (0.32) & 25.31 (0.33) & 25.24 (0.8335) \\
Grid 20Q4--24Q4 (17 qtrs) & 24.62 (0.28) & 27.51 (0.79) & 24.62 (0.28) & 24.85 (0.29) & 24.91 (0.7251) \\
\bottomrule
\end{tabular*}
\par\vspace{2pt}\begin{minipage}{\textwidth}\renewcommand{\baselinestretch}{1.00}\footnotesize
\textit{Notes.} Atlas-split uses disjoint halves for structure selection and fitting (chronological 364/364 days on the grid). The usual screening rule applies to the first half. Atlas-MLP-prob retains the charts and fits a two-hidden-layer ReLU head with validation-selected learning rate in $\{0.001,0.003\}$. Paired SEs against Atlas are 0.07/0.21/0.04 in CAP/PORT/NET and 0.14/0.10 on the early/later grid. 5-anchor Atlas averages the cross-covariance matrices evaluated at the SAA decision for the full chart and at the SAA decision for each of its four K-means subgroups. Parentheses in that column give the worst-case minimum principal cosine with the default basis. Differences from the default (Atlas minus five-anchor) are close to zero in the synthetic problems (paired SE 0.00--0.01) and are $-0.13$ (SE 0.09) and $-0.30$ (SE 0.14) in the early and later grid periods.
\end{minipage}
\end{table}

\begin{table}[htbp]
\centering

\caption{Additional robustness checks.}
\label{tab:ecRobustnessChecks}
\scriptsize
\renewcommand{\arraystretch}{1.03}
\setlength{\tabcolsep}{0pt}
\begin{tabular*}{\textwidth}{@{\extracolsep{\fill}}llrr@{}}
\toprule
Diagnostic & Setting & Metric & Reference \\
\midrule
\multirow{2}{*}{ICEO 3,200-step (CAP, 3 reps.)} & $\mathrm{lr}=10^{-3}$ (1,600 / 3,200) & 0.38\% / 0.35\% & Atlas 0.09\% \\
 & $\mathrm{lr}=3\!\times\!10^{-3}$ / $10^{-2}$ (3,200) & 0.66\% / 1.06\% & Atlas 0.09\% \\
\multirow{2}{*}{ICEO 3,200-step (PORT, 3 reps.)} & $\mathrm{lr}=10^{-3}$ (1,600 / 3,200) & 3.08\% / 3.19\% & Atlas 2.88\% \\
 & $\mathrm{lr}=3\!\times\!10^{-3}$ / $10^{-2}$ (3,200) & 4.10\% / 5.79\% & Atlas 2.88\% \\
Grid 11-qtr ICEO 3-LR sweep & $\mathrm{lr}\in\{3\!\times\!10^{-4},10^{-3},3\!\times\!10^{-3}\}$ & 25.08 / 25.15 / 25.22 & Val 25.11 (Atlas 25.10) \\
\midrule
Atlas iter.\ 0 (init.) & CAP/PORT/NET; Grid 11 qtrs & 0.05/33.27/0.90\%; 31.71 & 0.07/2.94/0.12\%; 25.10 \\
\midrule
Auto-$\widehat B$ diag.\ (\S\ref{ec:chartMechanismProtocol}) & min.\ op.\ loss ($0^\circ/45^\circ/90^\circ$) & $\widehat B=16$ (1.65/1.64/1.58\%) & Standalone; rejects $B=32$ \\
\midrule
Spectral gap $\sigma_{d_b}/\sigma_{d_b+1}$ & CAP / PORT / NET median & 13.77 / 1.05 / 1.09 & Energy $\geq 0.9999$ \\
\midrule
\multirow{2}{*}{Grid HAC $t$ ($L_{\rm NW}=7..28$)} & Early 11 qtrs: SOF / kNN-red & $-2.62$..$-2.12$ / $-5.51$..$-4.86$ & ICEO: $-0.38$..$-0.34$ \\
 & Later 17 qtrs: SOF / kNN-red & $-5.03$..$-4.58$ / $-4.86$..$-4.15$ & ICEO: $+4.10$..$+4.25$ \\
\midrule
kNN-red 32c (ser.$\to$par., s) & CAP / PORT / NET & 0.41$\to$0.21; 1.00$\to$0.24; 0.68$\to$0.40 & Atlas 0.04/0.03/0.21 \\
\bottomrule
\end{tabular*}
\end{table}

\begin{table}[htbp]
\centering

\caption{Head-matched affine-softmax ICEO comparisons.}
\label{tab:ecHeadMatched}
\footnotesize
\renewcommand{\arraystretch}{1.03}
\setlength{\tabcolsep}{0pt}
\begin{tabular*}{\textwidth}{@{\extracolsep{\fill}}llrrrr@{}}
\toprule
Application & Method (affine head) & Params $P$ & Metric (SE) $\downarrow$ & \shortstack[r]{Atlas $-$ Method\\(paired SE)} & \shortstack[r]{Atlas\\lower} \\
\midrule
\multirow{2}{*}{CAP (10 reps., $\mathrm{lr}=10^{-3}$)} & ICEO-affine ($K=500$) & 4,491 & 0.51 (0.02) & $-0.44$ (0.03) & 10/10 \\
 & ICEO-red-affine ($K=4$) & 27 & 0.56 (0.04) & $-0.49$ (0.04) & 10/10 \\
\multirow{2}{*}{CAP (reps.\ 1--3, val-sel.\ LR)} & ICEO-affine ($K=500$) & 4,491 & 0.35 (0.07) & $-0.26$ (0.06) & 3/3 \\
 & ICEO-red-affine ($K=4$) & 27 & 0.32 (0.11) & $-0.23$ (0.09) & 3/3 \\
\midrule
\multirow{2}{*}{PORT (10 reps., $\mathrm{lr}=10^{-3}$)} & ICEO-affine ($K=500$) & 4,491 & 3.71 (0.14) & $-0.77$ (0.06) & 10/10 \\
 & ICEO-red-affine ($K=9$--$18$) & 72--153 & 8.60 (3.14) & $-5.66$ (3.02) & 10/10 \\
\multirow{2}{*}{PORT (reps.\ 1--3, val-sel.\ LR)} & ICEO-affine ($K=500$) & 4,491 & 4.66 (0.78) & $-1.78$ (0.56) & 3/3 \\
 & ICEO-red-affine ($K=14$--$18$) & 117--153 & 5.70 (0.73) & $-2.82$ (0.58) & 3/3 \\
\midrule
NET (3 reps., $\mathrm{lr}=10^{-3}$) & ICEO-red-affine ($K=10$--$11$) & 81--90 & 0.22 (0.03) & $-0.10$ (0.03) & 3/3 \\
NET (3 reps., val-sel.\ LR) & ICEO-red-affine ($K=10$--$11$) & 81--90 & 0.14 (0.00) & $-0.02$ (0.01) & 3/3 \\
\midrule
\multirow{2}{*}{Grid 18Q1--20Q3 (11 qtrs, val-sel.\ LR)} & ICEO-affine ($K=728$) & 12,359 & 25.63 (0.36) & $-0.53$ (0.17) & 10/11 \\
 & ICEO-red-affine ($K=8$--$9$) & 119--136 & 27.09 (0.36) & $-1.98$ (0.17) & 11/11 \\
\bottomrule
\end{tabular*}
\par\vspace{2pt}\begin{minipage}{\textwidth}\renewcommand{\baselinestretch}{1.02}\footnotesize
\textit{Notes.} Both ICEO variants use Atlas's affine-softmax head. Validation selects from $\{0.001,0.01,0.025\}$ in CAP/PORT/NET and $\{0.001,0.025\}$ on the grid, with other hyperparameters at benchmark defaults. NET covers reduced-support ICEO only, and the grid covers the first 11 quarters.
\end{minipage}
\end{table}
\endgroup

\end{document}